\documentclass{amsart}
\usepackage{amssymb}
\usepackage{eucal}
\usepackage{amsfonts}
\usepackage{epsfig}
\usepackage[frame,ps,dvips,matrix,arrow,curve,rotate]{xy}
\let\oldcite\cite                                  
\newtheorem{thm}{Theorem}[section]
\newtheorem{cor}[thm]{Corollary}
\newtheorem{lem}[thm]{Lemma}

\newtheorem{prop}[thm]{Proposition}
\theoremstyle{definition}
\newtheorem{defn}[thm]{Definition}
\theoremstyle{remark}
\newtheorem{rem}[thm]{Remark}
\numberwithin{equation}{section}
\newtheorem{ex}[thm]{Example}
\theoremstyle{remark}

\let\:=\colon
\newcommand{\x}{\times}

\newcommand{\im}{\operatorname{im}}
\newcommand{\Ob}{\operatorname{Ob}}
\newcommand{\Mor}{\operatorname{Mor}}

\newcommand{\Hom}{\operatorname{Hom}}

\newcommand{\id}{\operatorname{id}}

\newcommand{\Tors}{\operatorname{Tors}}
\newcommand{\Glue}{\operatorname{Glue}}

\def\smashedlongrightarrow{\setbox0=\hbox{$\longrightarrow$}\ht0=1pt\box0}
\def\risom{\buildrel\sim\over{\smashedlongrightarrow}}

\def\smashedst{\setbox0=\hbox{$\rightrightarrows$}\ht0=4pt\box0}
\newcommand{\sst}[1]{\stackrel{#1}{\smashedst}}

\def\twomorphism{\setbox0=\hbox{$\Rightarrow$}\ht0=4pt\box0}
\newcommand{\twomor}[1]{\stackrel{#1}{\twomorphism}}

\newcommand{\Ra}{\Rightarrow}

\newcommand{\lra}{\longrightarrow}
\newcommand{\llra}[1]{\stackrel{#1}{\lra}}
\newcommand{\rsa}{\rightsquigarrow}

\newcommand{\hra}{\hookrightarrow}

\newcommand{\pX}{\pi_1(\X,x)}

\newcommand{\pY}{\pi_1(\Y,y)}

\newcommand{\bbZ}{\mathbb{Z}}
\newcommand{\bbR}{\mathbb{R}}
\newcommand{\bbC}{\mathbb{C}}
\newcommand{\bbN}{\mathbb{N}}

\newcommand{\bbQ}{\mathbb{Q}}

\newcommand{\al}{\alpha}
\newcommand{\be}{\beta}
\newcommand{\ep}{\epsilon}
\newcommand{\De}{\Delta}

\newcommand{\vf}{\varphi}

\newcommand{\X}{\mathcal{X}}
\newcommand{\Y}{\mathcal{Y}}
\newcommand{\Z}{\mathcal{Z}}
\newcommand{\PP}{\mathcal{P}}
\newcommand{\U}{\mathcal{U}}
\newcommand{\V}{\mathcal{V}}
\newcommand{\W}{\mathcal{W}}

\newcommand{\LL}{\mathcal{L}}

\newcommand{\C}{\mathsf{C}}
\newcommand{\B}{\mathcal{B}}
\newcommand{\A}{\mathcal{A}}
\newcommand{\OO}{\mathcal{O}}

\newcommand{\T}{\mathcal{T}}
\newcommand{\I}{\mathcal{I}}

\newcommand{\Xm}{X_{mod}}
\newcommand{\Ym}{Y_{mod}}
\newcommand{\pim}{\pi_{mod}}

\newcommand{\IX}{\mathcal{I}_{\X}}
\newcommand{\IY}{\mathcal{I}_{\Y}}

\newcommand{\ox}{\omega_x}
\newcommand{\oy}{\omega_y}

\begin{document}

\title[Foundations of topological stacks I]
{Foundations of topological stacks I}%
\author{Behrang Noohi}

\begin{abstract}
This is the first in a series of papers devoted to foundations
of topological stacks.

We begin developing a homotopy theory for topological stacks
along the lines of classical homotopy theory of topological
spaces. In this paper we go as far as introducing the homotopy
groups and establishing their basic properties. We also develop
a Galois theory of covering spaces for a (locally connected
semilocally 1-connected) topological stack.  Built into the
Galois theory is a method for determining the stacky structure
(i.e., inertia groups) of covering stacks. As a consequence,
we get for free a characterization of topological stacks that
are quotients of topological spaces by discrete group actions.
For example, this give a handy characterization of good orbifolds.

{\em Orbifolds}, {\em graphs of groups}, and {\em complexes of
groups} are examples of topological (Deligne-Mumford)  stacks. We
also show that any algebraic stack (of finite type over $\bbC$)
gives rise to a topological stack. We also prove a Riemann
Existence Theorem for stacks. In particular, the algebraic
fundamental group of an algebraic stack over $\mathbb{C}$ is
isomorphic to the profinite completion of the fundamental group of
its underlying topological stack.

The next paper in the series concerns function stacks (in
particular loop stacks) and fibrations of topological stacks.
This is the first in a series of papers devoted to foundations of
topological stacks.
%
%
%
%
%
\end{abstract}

\maketitle%
\newpage
\tableofcontents%
\newpage


\section{Introduction}{\label{S:introduction}}

The theory of algebraic stacks was invented by Deligne and Mumford
(largely based on ideas of Grothendieck)  and extended later on by
Artin, and has ever since become an indispensable part of
algebraic geometry.

Topological stacks, however, have only appeared in heuristic
arguments here and there, and no solid  theory for them exists in
the literature.\footnote{There is a paper by D. Metzler
\oldcite{Metzler} which has a bit of overlap (in some definitions)
with this paper, but has a different focus.} Since topological
stacks carry the underlying ``topology'' of their ubiquitous
algebraic counterparts, need for a systematic development for the
homotopy theory (and, for that matter, measure theory,
differential geometry and so on) seems   inevitable. Furthermore,
it is expected that (a substantial portion  of) the  theories of
topological and differential groupoids (e.g. works of Weinstein,
Moerdijk and others) to find their natural home in the theory of
topological stacks.

I should point out that certain special classes of  topological
stack have already been around for decades, under secret names of
{\em graphs of groups} (Serre, Bass), and {\em orbifolds}
(Thurston). The theory developed in this paper brings together
these two theories under the general umbrella of topological
stacks.

My goal has been to set up a machinery
 that is general enough so it can be applied in variety of situations.
The slogan is that, having chosen the correct definitions,
topological stacks can be treated pretty much the same way as
 topological spaces and,
modulo keeping track of certain extra structures that are bundled
with them (I will try to explain this below), we can do with
topological stack what we do with usual topological spaces (e.g.
homotopy and (co)homology theories, differential geometry, measure
theory and so on).

\vspace{0.1in} \noindent $\blacktriangleright$ Where do
topological stacks come from? \vspace{0.1in}

  Topological
stacks  provides a unified way to treat {\em equivariant}
problems,  as if they were non-equivariant:
 in most situations,
one can formulate an equivariant problem about a $G$-space $X$ as
a problem about the quotient stack $[X/G]$. This point of view
also provide a better functorial grasp on the problems, especially
in  situations where there is a change of groups involved.

The theory of topological stacks is a natural framework to study
topological groupoids, and the homotopy theory developed here
enables us to apply the machinery of algebraic topology in the
study of topological groupoids. Morita invariant phenomena belong
naturally to the realm of topological stacks.

Also, as in algebraic geometry, topological stacks can be useful
in the study of moduli problems (e.g. classifying spaces etc.).
Quotient constructions, especially in the presence of fixed
points, are more naturally performed in the framework of stacks.
The pathological behavior of quotient spaces, as well as the
enormous loss of information, is well taken care of if we use the
stacky approach.

Topological stacks can also appear as ``underlying spaces'' of
algebraic (or differential, analytic etc.)  stacks,  much the same
way that   topological space appear as underlying spaces of
algebraic varieties (or manifolds, complex manifolds etc.). The
homotopy theoretical properties of the underlying topological
stack is bound to play an important role in the study of the
original algebraic stack.

\vspace{0.1in} \noindent $\blacktriangleright$ How to think about
topological stacks? \vspace{0.1in}

Philosophically speaking, if one day it was declared that  {\em
set} theory, as a foundation for mathematics, should be thrown
away and be replaced by {\em groupoid} theory, then we would have
to replace topological spaces, complex manifolds, schemes and so
on, by topological stacks, analytic stacks, algebraic stacks and
so on!

Intuitively speaking, one way to think of a topological stack $\X$
is to think of a ``topological space'' in which a point is no
longer a single point, but a cluster of points that are equivalent
to each other.\footnote{To make our simplistic description closer
to reality, we should consider also $W$-valued points for $W$ an
arbitrary topological space. These are  maps from $W$ to $\X$.}
There are ``equivalences'' between pairs of points in such a
cluster, {\em possibly more than one for each pair}, and we would
like to keep track of these equivalences; that is why we build
them into the structure of our topological stack. If we chose to
actually identify all the points in each cluster, we end up with
an honest topological space $\Xm$, called the {\em coarse moduli
space} of $\X$. If we make a choice of a point $x$ in a cluster,
it comes with an additional structure: the group of equivalences
from $x$ to itself. This is called the {\em inertia group} of $x$,
and is denoted by $I_x$. Different points in the same cluster have
isomorphic inertia groups. A very rough intuitive picture for $\X$
would be then to think of it as the collection of $I_x$ bundled
together along the topological space $\Xm$.

Of course, this picture is extremely handicapped, and one should
be very cautious as there could occur quite a lot of pathological
phenomena that are overlooked in this simplistic description. One
is advised to work out a whole lot of examples (especially
pathological ones) so as to adjust the old fashioned topological
intuition to this situation.

An essential part of the intuition required to understand a stack
is based on the following meta-mathematical principle: never
identify equivalent (in a loose sense of the term) objects, only
remember that they {\em can} be identified, and remember
 {\em the ways} they can be identified (i.e. remember all
the equivalences between them).

Let me give a very simple example. Take a (discrete) set $X$, and
let $G$ be a discrete group acting on it. We would like to define
the quotient of this action. The old fashioned way to do this is
to look at the set $X/G$ of orbits. That is,  a given pair of
points  $x$ and $y$ in $X$ are identified, if there is a $g \in G$
sending $x$ to $y$. This construction, however, violates the
principle mentioned in the previous paragraph. So let us modify
the construction a little bit. What we want, instead of
identifying $x$ and $y$, is to remember that they {\em can} be
identified via $g$. To do so, we draw an arrow from $x$ to $y$,
and put a label $g$ on it. So now we have a set $X$ with a
collection of arrows between them, labeled by elements of $G$.
 It is easy to see that
there is a natural way to compose the arrows. So we have actually
constructed a category whose set of objects is $X$. This category
is indeed a groupoid, and is sometimes referred to as the {\em
translation groupoid} of the action of $G$ on $X$. We think of
this groupoid as the quotient of the action of $G$ on $X$, and
denote it by $[X/G]$. This is a baby example of a {\em quotient
stack}.

Of course this is a discrete example and may seem not so
interesting, but recall that Grothendieck tells us that any
``space'' $X$ (e.g. topological space, differential manifold,
scheme and so on) is just a collection of {\em sets}, coexisting
in a compatible way, i.e. in the form of a sheaf. For instance,
the information carried by a topological space $X$ is completely
captured by the sheaf $\mathsf{Top} \to \mathsf{Set}$ it
represents (Yoneda lemma). So, after all, everything boils down to
sets. What if we want to have {\em groupoids} instead of sets (we
just saw how groupoids turned up in the construction of our
discrete quotient stack)?

The groupoid version of  a sheaf is  what is called a {\em stack}.
So, a stack to a (discrete) groupoid is what a sheaf is to a set.
The axioms of a stack are jazzed up versions of the axioms of
sheaves (we have one more axiom actually), and they are sometimes
called {\em descent conditions}. They are to be thought of as
local-to-global conditions. The following examples explains what
descent conditions mean. Assume $\X$ is a stack over the category
$\mathsf{Top}$ of topological space, and let $W$ be a topological
space. Suppose we are given an open covering $\{W_i\}_{i \in I}$
of $W$, and for each $i$ we are given a map $f_i \: W_i \to \X$.
Assume over each double intersection  $W_{ij}$, $f_i$ can be
identified with $f_j$ (but, by the principle mentioned above, we
will not say that $f_i$ is equal to $f_j$ on $W_{ij}$), and we
record this identification by giving it a name $\vf_{ji} \: f_i
\twomor{} f_j$. To make things as compatible as possible, we
require that $\vf_{ij}$ is the inverse of $\vf_{ji}$, and that
over a triple intersection $W_{ijk}$  we have
$\vf_{ij}\circ\vf_{jk}=\vf_{ik}$ (cocycle condition). Having made
all these provisions,  it is natural to expect that the maps $f_i$
can be glued along the identification $\vf_{ij}$ to form a global
map $f \: W \to \X$, and we want that to happen in an essentially
unique way. This is what is called the {\em descent condition} or
{\em stack condition}.

Of course, not every sheaf on $\mathsf{Top}$ is ``topological
enough'' (i.e. comes from a topological space), so one should not
expect that any old stack to be ``topological enough'' either.
That is, simply imposing stack conditions is not enough for the
purpose of doing topology. {\em Topological} stacks are stacks
over $\mathsf{Top}$ which satisfy  some extra conditions. This
conditions appear in the form of existence of what we call a {\em
chart} for the stack. We may perhaps also require that certain
axioms  be satisfied by this chart (these will be discussed in
more detail in the text). Having imposed these conditions, one can
{\em pretend} that the stack came from a ``topological object'',
and the simplistic picture I gave a few paragraphs back is a way
of visualizing that topological object.

\vspace{0.1in}

\noindent$\blacktriangleright$ Structure of the paper
\vspace{0.1in}

The paper is divided into two parts. The first part is more on the
formal side and is devoted to  setting up the basic definitions
and constructions related to  stacks over $\mathsf{Top}$. No
homotopy theory appears in Part I. In the Part II,  we begin
developing the homotopy theory. To keep the paper within bounds,
we will not go that far into homotopy theory, and leave this task
to the forthcoming papers in this series.

 \vspace{0.1in}
 The sections are organized as follows:
 \vspace{0.1in}

I have collected some notations and conventions in Section
\ref{S:Notation}. The reader who starts browsing from the middle
of the paper and runs into some unfamiliar notation or
terminology, or is unsure of the hypotheses being made all through
the paper, may find this section helpful.

Section \ref{S:stacks} contains some general facts about stacks
over Grothendieck sites. The goal here is to set the notation, and
collect some general results for the sake of future reference in
the paper. It is by no means intended to be an introduction to
stacks over Grothendieck sites -- the reader is assumed to be
familiar with the general machinery of stacks (see \cite{Metzler},
\cite{LM}), and also with the language of 2-categories: stacks
naturally form a 2-category and we will heavily exploit this
feature.

There are five subsections to Section \ref{S:stacks}. Subsection
\ref{SS:fiber} recalls the definition of 2-fiber products and
2-cartesian diagrams. Subsection  \ref{SS:GpdvsStack} reviews the
correspondence between groupoids (in the category of presheaves
over a Grothendieck site $\C$) and stacks over $\C$. More
precisely, given such a groupoid, one can construct its {\em
quotient stack}. Conversely, any stack is of this form (up to
equivalence). Subsection \ref{SS:image} is devoted to the
definition of {\em image} of a morphism of stacks.
The main feature of a  stack $\X$ (as opposed to a sheaf) is that
any point $x$ on $\X$ comes with a group  attached to it. This
group is called the {\em inertia} group or the {\em isotropy
group} of $x$, and is denoted by $I_x$. The inertia group is
defined as the group of 2-isomorphisms from the map $x \: * \to
\X$ to itself (viewed in the 2-category of stacks). In examples
where $\X$ is a quotient stack, $I_x$ is isomorphic to the
stabilizer  group of $x$. The group $I_x$ records how the point
$x$ has been ``over-identified'' with itself. So, for instance, if
$x$ has inertia group of order 5, we want to think of $x$ as
$\frac{1}{5}$th of a point!

Related to the notion of inertia group of a point are the notions
of the inertia sheaf and the {\em residue gerbe} of a point. The
inertia groups assemble together to form a global object $\IX$,
called the {\em inertia stack}. The inertia stack comes with a map
$\IX \to \X$, which makes it into a relative group object over
$\X$. All these are discussed briefly in Subsection
\ref{SS:inertia1} of Section \ref{S:stacks}.  In Subsection
\ref{SS:maps} we give a description of maps coming out of a
quotient stack $[X/R]$
 in terms of the groupoid $[R\sst{}X]$ itself.
The description is given in terms of maps coming out of $X$ whose
restriction to $R$ satisfy certain cocycle conditions.

From Section \ref{S:topologicalsite} we narrow down to our
favorite Grothendieck site, namely the category $\mathsf{Top}$ of
compactly generated topological spaces. The Grothendieck topology
is defined using the usual notion of covering by open
subsets.\footnote{I am convinced there is no advantage in taking
other fancy topologies, unless someone proves me wrong.}
 In $\mathsf{Top}$ we can be more
concrete and talk about some specific topological features of
stacks. Nevertheless, a stack over $\mathsf{Top}$, as it is, is
still to crude to do topology on. Some general considerations,
however,
 can already be made at this level. There are three subsections
 to Section \ref{S:topologicalsite}.
 In Subsection \ref{SS:representable},
we define the notion of a {\em representable} map between stacks,
and discuss certain properties of representable maps such as local
homeomorphism, open, closed and so on. In Subsection
\ref{SS:operations}, we go over certain basic operations
(intersection, union, closure, image, and inverse image) on
substacks. A new feature of stacks over the site $\mathsf{Top}$ is
that, to any such stack $\X$ we can associate an honest
topological space $\Xm$, the {\em coarse moduli space} of $\X$.
This is discussed in Subsection \ref{SS:coarse}. There is natural
map $\X \to \Xm$ which is universal among maps from $\X$ to
topological spaces. So, in a way, $\Xm$ is the best approximation
of $\X$ by a topological space. We sometime call $\Xm$ the {\em
underlying topological space} of $\X$.

Sections \ref{S:gerbes}  and \ref{S:criterion} should really be
thought of as subsections of Section \ref{S:topologicalsite}. In
Section \ref{S:gerbes} we quickly recall what a gerbe is, and look
at certain specific features of gerbes over $\mathsf{Top}$.
Section \ref{S:criterion} is devoted to an easy but very useful
technical fact about representable maps. It is a criterion for
verifying whether a map of stacks is representable. It says that,
in order to check whether a map is representable, it is enough to
do so after base extending the map along an epimorphism. This
saves a lot of hassle later on  where we have to frequently deal
with representable maps.

In Section \ref{S:pretopological} we add flesh to the bone and
introduce pretopological stack. These are stack over
$\mathsf{Top}$ which admit a {\em chart}, that is, a representable
epimorphism $p \: X \to \X$ from a topological space $X$.
Equivalently, a pretopological stack is equivalent to the quotient
stack of a topological groupoid. The existence of a chart makes a
stack more flexible, as we can use charts to transport things back
to the world of topological spaces. However, pretopological stack
are still too crude and not good enough to do topology on. The
point is that, to introduce enough flexibility, so as to be able
to carry out topological constructions, one has to impose some
conditions on the chart $p \: X \to \X$. This will be done in
Section \ref{S:TopSt}. But before that, we establish some general
facts about pretopological stacks.

Different topological groupoids may give rise to the same quotient
stack (up to equivalence). This phenomenon is usually referred to
as {\em Morita equivalence} of groupoids. This is briefly
discussed in Section \ref{S:Morita}.

A pretopological stack roughly corresponds to a Morita equivalence
class of topological groupoids. Furthermore, by choosing suitable
representatives in the corresponding Morita equivalence classes, a
morphism of stacks can also be realized on the level of
topological groupoids. An important question that arises  is how
to compute 2-fiber products of pretopological stack using their
groupoid presentations. We answer this question in Section
\ref{S:fiber} by giving an explicit description of a topological
groupoid that represents the 2-fiber product. This description  is
very useful when we want to prove that a 2-fiber product satisfies
a certain property.

In Section \ref{S:inertia} we go back to inertia sheaves and
residue gerbes of points on a stack, now with the additional
assumption that the stack admits a chart. In Section
\ref{S:coarse2} we look at the coarse moduli space of
pretopological stacks. It turns out that the coarse moduli space
$\Xm$ of a pretopological stack $\X$ is
 reasonably well-behaved if we
assume that $\X$ admits a chart $p\: X \to \X$ with $p$  an {\em
open} map. In this case, we prove   the invariance of coarse
moduli space under base change (Corollary \ref{C:basechange}).

There is an alternative description of the quotient stack of a
topological groupoid in terms of torsors. This is discussed in
Section \ref{S:quotient}.  This explicit description of a quotient
stack comes handy sometimes. It also explain why quotient stacks
are related to classifying spaces. For example, when $G$ is a
topological group acting trivially on a point, the quotient stack,
denoted by $BG$ exactly classifies $G$-torsors.

\vspace{0.1in} This is the end of Part I. \vspace{0.1in}

In Section \ref{S:TopSt} we  introduce topological stacks. These
are our main objects of interest, and the ones on which we can
perform enough topological constructions to have a reasonable
homotopy theory. The topological stack are essentially
pretopological stacks which admit a chart $p \: X \to \X$ that is
``nicely behaved''. To make sense of what we mean by ``nicely
behaved'', in Subsection \ref{SS:fibrations} we introduce the
notion of  {\em local fibration}. This is defined to be a class of
$\mathbf{LF}$ of continuous maps satisfying certain axioms. The
axioms are modeled on the usual notion of a local fibration; that
is, maps $f \: X \to Y$ such that, after passing to appropriate
open covers of $X$ and $Y$, become a union of fibrations. Local
 Serre fibrations, local Hurewic  fibrations and locally
 cartesian maps are examples of such a class $\mathbf{LF}$.

 For any choice of $\mathbf{LF}$ we have the corresponding theory of topological
 stacks, and the properties of $\mathbf{LF}$ reflect the properties of
 the corresponding category of topological stacks.
 We can use this to adjust  our class of local fibrations to ensure
 that the corresponding topological stacks have the desired properties.

The role played by the notion of fibration becomes more apparent
in Section \ref{S:gluing} where we prove some gluing lemmas for
topological stack. Roughly speaking, the point is that, in the
category of pretopological stack gluing along {\em closed}
substacks is problematic. For instance, very simple push-out
diagrams in the category of topological spaces may no longer be
push-out diagrams when viewed in the category of pretopological
stacks. In particular, it is not possible do define homotopy
groups of pretopological stacks (because, for example, a join of
two circles does not satisfy the expected universal property when
viewed in the category of pretopological stacks). However, once we
impose the fibrancy condition on charts, the push-out diagrams
start behaving better, provided we put certain {\em local
cofibrancy} conditions on the arrows in the diagram. The bottom
line is, more restrictions on the charts will result in higher
flexibility in performing push-outs.

Parallel to the theory of Deligne-Mumford stacks in algebraic
geometry, we develop a theory of Deligne-Mumford topological
stacks.
 This is done in Section \ref{S:DM}.
A {\em weak Deligne-Mumford} topological stack is defined to be a
topological stack that admits a chart $p \: X \to \X$ that is a
local homeomorphism. These are a bit too general and not as
well-behaved as they are expected to be (as will be seen in some
examples in the text). The good notion is that of a {\em
Deligne-Mumford} topological stack; it is, by definition, a weak
Deligne-Mumford stack that is locally a quotient stack of a
properly discontinuous action (Definition \ref{D:pd}). We will not
make any finiteness assumptions on the stabilizers.

Dropping the finiteness condition in the definition of a
Deligne-Mumford stack, however, causes some problems with
existence of 2-fiber product, i.e. they may no longer be
Deligne-Mumford. In Section \ref{S:DMFiber} we see that, under
very mild conditions,
 2-fiber products of Deligne-Mumford topological stacks
 will again be Deligne-Mumford stacks (Corollary \ref{C:DMfiber}).

In  Section \ref{S:homotopy} we start setting up the basic
homotopy theory of topological stacks. We define homotopies
between maps, and homotopy group of topological stacks. All the
structures carried by homotopy groups of topological spaces (e.g.
action of $\pi_1$ on higher homotopy groups, Whitehead products
and so on) can be carried over to topological stacks. The homotopy
groups of topological stacks, however, carry certain extra
structure. More precisely, for a point $x$ on a topological stack
$\X$, we have natural group homomorphisms
  $$\pi_{n}(\Gamma_x,x) \to \pi_n(\X,x),$$
where $\Gamma_x$ is the residue gerbe at $x$. The importance of
these maps is that they relate local invariants (i.e. residue
gerbes) to  global invariants (homotopy groups). When $n=1$,
$\pi_{1}(\Gamma_x)$ is isomorphic to the inertia group $I_x$ at
$x$. Therefore, we have natural maps
  $I_x \to \pi_1(\X,x)$. Note that the left hand side depends on the
  base point, while the right hand  side does not (up to isomorphism).
We discuss this map in some detail in Section  \ref{S:covering}.

In Subsection \ref{SS:Galois} of Section \ref{S:covering}, we
develop a Galois theory of covering spaces for topological stack
and show that every locally path connected semilocally 1-connected
topological stack has a universal cover (Corollary
\ref{C:universalcover}). In fact, we prove an equivalence between
the category of $\pi_1(\X)$-sets and the category of covering
spaces of $\X$ (similar to Grothendieck's theory in SGA1). The
covering theory developed here generalizes the covering theory of
orbifolds developed by Thurston in \cite{Thurston} and the
covering theory of graphs of groups developed by Bass in
\cite{Bass}.

In Subsection  \ref{SS:hidden} of Section \ref{S:covering},  we
investigate the role played by the maps $I_x \to \pi_1(\X,x)$
mentioned above in the Galois theory of covering spaces. One main
result is that, a Deligne-Mumford stack is globally a quotient of
a (properly discontinuous) group action if and only if all the
maps $I_x \to \pi_1(\X,x)$ are injective (Theorem
\ref{T:uniformization}). This is a very useful and practical
criterion, as the maps $I_x \to \pi_1(\X,x)$ are usually very easy
to compute. This, in particular, gives a necessary and sufficient
condition for an orbifold to be a {\em good orbifold} (in the
sense of Thurston). In a separate paper \cite{Noohi2} we show how
the maps $I_x \to \pi_1(\X,x)$ can be used to compute the
fundamental group of the coarse moduli space $\Xm$. Essentially,
the idea is that $\pi_1(\X)$ is obtained by killing the images of
all the maps $I_x \to \pi_1(\X,x)$, for various points
$x$.\footnote{This, in particular, gives a formula for computing
the fundamental group of the the (coarse) quotient of a
topological group acting on a topological space (possibly with
fixed points).}

In Section \ref{S:examples} we look at some examples of
topological stacks. There are five subsections. In Subsection
\ref{SS:pathological} we collect some pathological examples. In
Subsection \ref{SS:gerbes} we talk a little bit about {\em
topological} gerbes. In Subsection \ref{SS:orbifolds} we show that
Thurston's orbifolds are Deligne-Mumford topological stacks. In
Subsection  \ref{SS:weighted} we look at weighted projective
lines. In Subsection  \ref{SS:graphs} we show that graphs of
groups are also Deligne-Mumford stacks. We point out how certain
general results of Serre and Bass are easy consequences of the
homotopy theory of topological stacks.


Finally, in Section \ref{S:algebraic}, we construct a functor from
the 2-category of algebraic stacks (of finite type over
$\mathbb{C}$) to the 2-category of topological stacks. We show
that, under this functor, Deligne-Mumford stacks go to
Deligne-Mumford topological stacks. We also prove a Riemann
existence for algebraic stack. In particular we deduce that, the
algebraic fundamental group of an algebraic stack (of finite type
over $\mathbb{C}$) is isomorphic to the profinite completion of
the (topological) fundamental group of the corresponding
topological stack.


\vspace{0.2in}\noindent{\bf Acknowledgment} I would like to thank
Kai Behrend and Angelo Vistoli, for sharing with me their ideas,
as well as their personal notes on topological stack. I also thank
Angelo Vistoli for inviting me  to visit Bologna, during which we
had interesting discussions about topological stacks. I would like
to thank Rick Jardine for always being helpful and supportive.
This paper  was written during my Imperial Oil post-doctoral
fellowship at the University of Western Ontario.

 \newpage


\section{Notation and conventions}{\label{S:Notation}}

For the convenience of the reader, in this section we collect some
of the notations and conventions used throughout the paper, and
fix some terminology. These are mostly items that are not
explicitly mentioned during the text.

\vspace{0.1in}

   Throughout the paper, all topological spaces are assumed to be
compactly generated.

  We use calligraphic symbols ($\X$,$\Y$,$\Z$,...) for stacks, categories
  fibered in groupoids and so on, and Roman
  symbols ($X$,$Y$,$Z$,...) for spaces, sheaves and so on.

    We fix a final object $*$, {\em the} point,
in the category $\mathsf{Top}$
 of topological
spaces. Except for in Section \ref{SS:GpdvsStack}, whenever we
talk about a {\em point} $x$ in a stack $\X$ what we mean is a
morphism $x \: * \to \X$ of stacks.

For an object $W$ in a Grothendieck site $\C$, we use the same
notation $W$ for the sheaf associated to the presheaf represented
by $W$, and also for the corresponding stack. For the presheaf
represented by $W$ we use $W_{pre}$.

   In the 2-category of stacks, the 2-isomorphism are
synonymously referred to as {\em transformations}, {\em
equivalences} or {\em identifications}. Given a pair $f,g \: \Y
\to \X$ of morphisms  of stacks, we use the notation
$f\twomor{\varphi}g$  or $\varphi \: f\twomor{}g$ to denote a
2-isomorphism between them. For ($*$-valued) points $x$ and $y$ in
a stack, we use the alternative notation $x \rsa y$ for an
identification to emphasis the resemblance with a path.

The 1-morphisms in a 2-category are referred to as {\em morphisms}
or {\em maps}.

Given a pair of morphisms of stacks $f \: \A \to \X$ and $p \:\Y
\to \X$, by a {\em lift} of $f$ to $\Y$ we mean a morphism
$\tilde{f} \: \A \to \Y$ together with an identification $\varphi
\: f \twomor{} p \circ \tilde{f}$ as in the following 2-cell:

$$\xymatrix@=24pt@M=4pt{                         & \Y \ar[d]^p   \\
      \A \ar[r]_{f} \ar@{} @<2pt> [r]  | (0.45){}="a" \ar @{.>} [ru]^{\tilde{f}}
                            \ar@{} @<-2pt> [ru]  | (0.3){}="b"
                      & \X       \ar @{=>}_{\varphi}  "a";"b"  }$$

When $\al$ and $\be$ are elements in some groupoid (e.g. groupoid
of sections of a stack) such that $t(\al)=s(\be)$, we use the
multiplicative notation $\al\cdot\be$ for their composition. We do
the same with composition of 2-isomorphism.

   A {\em substack} is always assumed to be saturated and full
   (Definition \ref{D:substack}).

   In some sloppy moments, we call an {\em isomorphism}
(of groupoids/stacks)  what should really be called an {\em
equivalence}.

    We use the notation $\mathcal{F} \mapsto \mathcal{F}^a$
     for the stackification and sheafification functors.

The terms {\em \'{e}tale} and  {\em local homeomorphism} are
synonymous. To be consistent with the classical terminology, we
use the term ``covering space/stack" in the topological setting
for what we would call a``finite \'{e}tale map" in the algebraic
setting.

   A groupoid object is denoted by $[R \sst{} X]$,
$[R \sst{s,t} X]$, or $[s,t \:R \sst{} X]$ and
 the symbols $s$ and $t$ are exclusively used for source and target maps.
For a subset $U \subseteq X$, the {\em orbit}
$t\big(s^{-1}(U)\big)$ of $U$ is denoted by $\OO(U)$.

When $[R \sst{} X]$ is a groupoid object in a category (say
$\mathsf{Top}$) and $Y \to X$ is  a map, we denote the
corresponding  {\em pull-back groupoid} by $[R|_Y \sst{} Y]$,
where $R|_Y:=(Y\x Y)\x_{X\x X} R$.


\newpage
{\LARGE\part{Preliminaries}}{\label{PART1}} \vspace{0.5in}

\section{Quick review of stacks over sites}{\label{S:stacks}}

In this section we collect a few facts about stacks that will be
used throughout the paper. It should not be regarded as an
introduction to stacks: the reader is assumed to be familiar with
the notions of Grothendieck topology, sheaf and presheaf, category
fibered in groupoids, stack and prestack, stackification and
2-category.

Throughout this section, we fix a category $\C$ with a
Grothendieck topology on it. All categories fibered in groupoids
will be over $\C$. We will be sloppy and not distinguish between a
presheaf and the corresponding category fibered in groupoids.

Let $\X$ be a category fibered in groupoids over $\C$. For an
object $W$ in $\C$, we have the corresponding {\bf groupoid of
sections} $\X(W)$. By Yoneda's lemma, there is a natural
equivalence of groupoids
  $$\Hom(W_{pre},\X)\cong \X(W),$$
where the left hand side is computed in the 2-category of
categories fibered in groupoids over $\C$. Here, $W_{pre}$ stands
for the category fibered in groupoids associated to the presheaf
represented by $W$.

The assignment $W \mapsto \X(W)$ gives rise to a {\em lax}
presheaf of groupoids.\footnote{This lax presheaf  can be
strictified, in the sense that it is  ``equivalent'' to the strict
presheaf $W \to \Hom(W_{pre},\X)$. But we don't really care.}
Conversely, given a (lax) presheaf of groupoids, one can construct
the associated category fibered in groupoids. For the sake of
brevity, we usually specify a category fibered in groupoids $\X$
by its groupoids of sections $\X(W)$.

\subsection{2-fiber products}{\label{SS:fiber}}

There is a notion of 2-fiber product for categories fibered in
groupoids. Consider a diagram
  $$\xymatrix@=12pt@M=10pt{
                     &  \Y \ar[d]^f  \\
            \Z  \ar[r]_(0.47)g        &     \X      }$$
where $\X$, $\Y$ and $\Z$ are categories fibered in groupoids over
$\C$. The {\bf 2-fiber product} $\Y\x_{\X}\Z$ is defined to be the
category given by

{\small
       $$\Ob(\Y\x_{\X}\Z)=\left\{\begin{array}{rcl}
              (y,z,\al) & | & y \in \Ob\Y, z \in \Ob\Z, \ \text{s.t.} \
                                      p_{\Y}(y)=p_{\Z}(z)=c \in \Ob\C;  \\
                        &   & \al \: f(y) \to g(z)  \ \text{an arrow in} \ \X(c).
                                 \end{array}\right\}$$
\vspace{0.1in}
  $$\Mor_{\Y\x_{\X}\Z}\big( (y_1,z_1,\al),(y_2,z_2,\be) \big)=\left\{\begin{array}{rcl}
              (u,v) & | &  u\: y_1 \to y_2 ,\ v\: z_1 \to z_2 \ \text{s.t.:} \\
                    &   & \bullet \ p_{\Y}(u)=p_{\Z}(v) \in \Mor\C,  \\
                    &   & \bullet \ \text{In} \Mor\X \ \text{we have}, \\
                    &   &   \hspace{0.5in}
         \xymatrix@R=8pt@C=10pt{
             f(y_1) \ar[r]^{f(u)}  \ar[d]_{\al} \ar@{}[rd]|{\circlearrowleft}
                                                   &  f(y_2) \ar[d]^{\be} \\
                            g(z_1)  \ar[r]_{g(u)}  &  g(z_2)   }
                                                          \end{array}\right\}$$}
Here $p_{\Y} \: \Y \to \C$ stands for  the structure map which
makes $\Y$ fibered  over $\C$.

It is easily checked that $\Y\x_{\X}\Z$ is again fibered in
groupoids over $\C$. The 2-fiber product should be viewed as a
homotopy fiber product. When $\X$, $\Y$ and $\Z$ are stacks, the
2-fiber product $\Y\x_{\X}\Z$ is again a 2-stack. We will not
prove this here, but the idea is that the process of
stackification consists of successive applications of  homotopy
limits combined with {\em filtered} homotopy colimits; in
particular, it commutes with finite homotopy limits.

\vspace{0.1in}
 The 2-fiber product has the following universal property:
  \begin{quote}
    For any category fibered in groupoids $\W$, the natural functor
    from the groupoid $\Hom(\W, \Y\x_{\X}\Z)$ to the groupoid of triples
    $(u,v,\al)$, where $u \: \W \to \Y$ and $v \: \W \to \Z$ are morphisms
    and $\al \: f \circ u \twomor{} g\circ v$ a 2-isomorphism
    (e.g. a transformation of functors relative to $\C$),
    is an equivalence.
  \end{quote}
 \vspace{0.1in}

 We say that a diagram
   $$\xymatrix@=12pt@M=10pt{
           \W \ar[r]^u \ar[d]_v &  \Y \ar[d]^f  \\
            \Z  \ar[r]_g        &     \X      }$$
 of categories fibered in groupoids (or stacks)
 is {\bf 2-commutative}, if there is a 2-isomorphism
 $\al \: f \circ u \twomor{} g\circ v$. A 2-commutative diagram is sometimes
 called a {\bf 2-cell} and is denoted by

 $$\xymatrix@=16pt@M=6pt{ \W \ar[r]^{u}  \ar@{} @<-4pt> [r]| (0.4){}="a"
               \ar[d]_{v} \ar@{} @<4pt> [d]| (0.4){}="b"
                                   & \Y \ar[d]^{f}     \ar @{=>}^{\al}  "a";"b" \\
                                        \X \ar[r]_{g}        & \X           }$$
 A 2-commutative diagram as above is called
{\bf 2-cartesian}, if
 the induced map $(u,v,\al) \: \W \to \Y\x_{\X}\Z$ is an equivalence.

\subsection{Groupoids vs. stacks}{\label{SS:GpdvsStack}}
Let $G$ be a (discrete) groupoid, and let $X$ be a set. Given a
map of sets $f \: X \to \Ob(G)$, we can pull back the groupoid
structure over to $X$ by setting $\Hom(x,y):=\Hom_G\big(p(x),p(y)
\big)$, for every $x,y \in X$. We obtain a groupoid $H$, with
$\Ob(H)=X$. The groupoid $H$ maps fully faithfully to $G$; this
map is an equivalence if and only if $f$ is surjective. The set
$R=\Mor(H)$ of morphisms of $H$ fits in the following 2-cartesian
diagram of groupoids:
    $$\xymatrix@=12pt@M=10pt{ R \ar[r] \ar[d]          &  \ar[d]^{\De} G \\
                 X\x X \ar[r]_{(f,f)}     &     G\x G      }$$

\noindent which can also be thought of as an alternative
definition for $R$ (hence for the pull back groupoid $H$).

We can do the same construction globally. Let $\X$ be a category
fibered in groupoids over $\C$ (or a presheaf of groupoids over
$\C$, if you wish), and let $X$ be a presheaf of sets (viewed as a
category fibered in  groupoids). Let $f \: X \to \X$ be a map of
categories fibered over $\C$. Define the presheaf (of sets) $R$ by
the following 2-cartesian diagram diagram:
    $$\xymatrix@=12pt@M=10pt{ R \ar[r] \ar[d]          &  \ar[d]^{\De} \X \\
                 X\x X \ar[r]_{(f,f)}      &    \X \x \X       }$$
\noindent Equivalently, $R\cong X\x_{\X}X$. This gives us a
groupoid object $[R \sst{} X]$ in the category of presheaves of
sets. For each $W \in \Ob\C$, the groupoid $[R(W) \sst{} X(W)]$ is
the groupoid induced from $\X(W)$ on the set $X(W)$ via the map
$X(W) \to \Ob\X(W)$, exactly as we saw in the previous paragraph.

\vspace{0.1in} \noindent{\bf Notation}. Let $[R \sst{} X]$ be a
groupoid object in the category of presheaves. We denote the
corresponding (strict) presheaf of groupoids, and  also the
corresponding category fibered in groupoids, by $\lfloor X/R
\rfloor$. \vspace{0.1in}

 \begin{defn}{\label{D:epi}}
    Let  $f \: \X \to \Y$ be a morphism of categories fibered in groupoids.
    \begin{itemize}

      \item $f$ is called an {\bf epimorphism}, if for every
    $W \in \Ob\C$
    and every $y \in \Ob\Y(W)$, there exists a covering family $\{U_i \to W\}$
    such that for every $i$, $y_{U_i}$ is (equivalent
    to an element) in the image of $\X(U_i)$.

      \item $f$ is called a {\bf monomorphism}, if for every $W \in \Ob\C$,
      the induced map $\X(W) \to \Y(W)$ is fully faithful.

    \end{itemize}

   \end{defn}

The following proposition is immediate.

 \begin{prop}{\label{P:fullyfaithful}}
 Let $\X$ be a category fibered in groupoids, $X$ a presheaf of sets
 (viewed as a category fibered in groupoids),
 and $f \: X \to \X$ a morphism. Set $R=X\x_{\X}X$, and let $[R \sst{} X]$
 be the corresponding groupoid. Then we have the following.

 \begin{itemize}

    \item[$\mathbf{i.}$] The natural
        map $\lfloor X/R \rfloor \to \X$ is a monomorphism. In particular,
        if $\X$ is a prestack, then $\lfloor X/R \rfloor$ is also
        a prestack.

    \item[$\mathbf{ii.}$] If  $f \: X \to \X$ is an epimorphism, then
        the induced map $\lfloor X/R \rfloor \to \X$ induces
        an equivalence  of stacks $\lfloor X/R \rfloor^a \to \X^a$, where
        $-^a$ stands for stackification.

 \end{itemize}

 \end{prop}

  \begin{defn}{\label{D:quotient}}
    Let $[R \sst{} X]$ be a groupoid object in the category of presheaves
    of groupoids over $\C$. We define the {\bf quotient stack} $[X/R]$
    to be the stack asscociated to $\lfloor X/R \rfloor$, that is,
    $\lfloor X/R \rfloor^a$. When $G$ is a presheaf of groups
    relative to a presheaf (of sets) $X$, we define $[X/G]$ to be the quotient
    stack of the trivial groupoid $[G\x X \sst{} X]$. This
    is sometimes denoted by $B_XG$ (or simply $\B G$, if the base $X$
    is understood), and  is called the {\bf classifying stack} of $G$.
  \end{defn}

 A special case of the Proposition \ref{P:fullyfaithful} is the
 following.

  \begin{prop}{\label{P:quotient1}}
    Let $\X$ be a stack, and let $X \to \X$ be an epimorphism
    from a presheaf of sets $X$
    (viewed as a category fibered in groupoids)
     to $\X$. Consider the groupoids object $[R \sst{} X]$
    (in the category of presheaves
    of sets over $\C$) defined by $R=X\x_{\X}X$.  Then, we have a
    natural equivalence $[X/R] \risom \X$.
  \end{prop}

 The following well-known fact is along the same lines of the previous
 proposition.

  \begin{prop}{\label{P:stackification}}
    Let $f \: \Y \to \X$ be a morphism of categories fibered in groupoids
    that is both a monomorphism and an epimorphism.
    Then the natural map $\Y^a \risom \X^a$ is an equivalence.
 \end{prop}

 The following facts also come handy sometimes.

  \begin{prop}{\label{P:sheaf}}
    Let $[R \sst{} X]$ be a groupoid object in the category of
    sheaves over $\C$. Then $\lfloor X/R \rfloor$ is a prestack.
  \end{prop}

   \begin{prop}{\label{P:prestack}}
      Let $\X$ be a prestack over $\C$. Then the natural map $\X \to \X^a$
      is a monomorphism.
   \end{prop}

 \vspace{0.2in}
\subsection{Image of a morphism of stacks} {\label{SS:image}}

 \begin{defn}{\label{D:substack}}
  Let $\X$ be a stack over $\C$. By a {\bf substack} of $\X$ we mean
  a subcategory $\Y$ (fibered in groupoids over $\C$) of $\X$ such that
  for every $W \in\Ob\C$, the subcategory $\Y(W) \subseteq \X(W)$
  is full and saturated (equivalently, $\Y$ is a full saturated
  subcategory of $\X$ as abstract categories). A subcategory
  is called {\em saturated} if whenever it contains an object then
  it contains the entire isomorphism class of that object.
 \end{defn}

 \begin{defn}{\label{D:image}}
    Let $f \: \Y \to \X$ be a morphism of stacks. The
 {\bf image} $\im(f)$  of $f$ is defined to be the smallest  substack
 of $\X$ through which $f$ factors.
\end{defn}

 \noindent An alternative definition for $\im(f)$ is that it
  satisfies the following properties:

  \begin{itemize}
    \item[$\mathbf{I1.}$] $\im(f)$ is a subcategory
     fibered in groupoids of $\X$ (not a priori full or saturated)
    through which $f$ factors;

    \item[$\mathbf{I2.}$] $f \: Y \to \im(f)$  is an epimorphism;

    \item[$\mathbf{I3.}$]
      $\im(f)$ is the largest
       subcategory fibered in groupoids of $\X$  satisfying
      ($\mathbf{I1}$) and ($\mathbf{I2}$).
  \end{itemize}

 \begin{prop}{\label{P:compose}}
  \end{prop}
 {\em
 \begin{itemize}

   \item[$\mathbf{i.}$] Let $g \: \Z \to \Y$ and
    $f \: \Y \to \X$ be morphisms of stacks.
    Then $\im(f|_{\im(g)})=\im(f\circ g)$.

   \item[$\mathbf{ii.}$] $f \: \Y \to \X$ is an epimorphism
        if and only if $\im(f)=\X$.

  \end{itemize}
 }

 \vspace{0.1in}

We also have the following way of computing the image:

 \begin{prop}{\label{P:image}}
    Let $\X$ be a stack, and let $f \: X \to \X$ be a map from a presheaf of
    sets $X$ to $\X$. Let $[R \sst{} X]$ be the corresponding groupoid
    in presheaves of sets (see Section \ref{SS:GpdvsStack}).
    Then, we have a natural equivalence $[X/R] \to \im(f)$
    (see Definition \ref{D:quotient}).
 \end{prop}

    \begin{proof}
       By Proposition \ref{P:fullyfaithful}, we have a monomorphism
       $\lfloor X/R \rfloor \to \X$. Let $\mathcal{I}$ be the image
       of this map. By Proposition
       \ref{P:stackification}, this induces an equivalence
       $[X/R] \to \mathcal{I}$.  Since $X \to [X/R]$ is an epimorphism,
       we have $\im(f)=\mathcal{I}$ by Proposition \ref{P:compose}$\mathbf{i}$.
    \end{proof}

\subsection{Inertia groups and  inertia stacks}{\label{SS:inertia1}}

Let $W$ be an object in the base category $\C$. By a $W$-{\bf
valued  point} of $\X$ we mean  an object $x$ in
 the groupoid $\X(W)$ (or the corresponding morphism $x \: W \to\X$).
 For such a point $x$,
 we can talk about the group of automorphisms of $x$, viewed as an object
 in the groupoid $\X(W)$. If we think of $x$ as a
 map $W \to \X$, then this is just the group
 of self-transformations of $x$.
 This group is called the {\bf inertia group}
 of $x$, and is denoted by $I_x$.
 Other synonyms  for inertia group are  {\em stabilizer group},
 {\em isotropy group} and {\em ramification group}. We may use these
 terms interchangeably here and there.

 The inertia groups $I_x$ of a $W$-valued point $x$ is in fact the set
 of $W$-points of a sheaf of groups $\mathbb{I}_x$ relative to $W_{pre}$,
 which we call the {\bf inertia presheaf} at $x$.
 As a presheaf over $\C$, the inertia presheaf $\mathbb{I}_x$
 is defined as follow:
     {\small $$\mathbb{I}_x(W') = \left\{ \begin{array}{rcl}
           (f,\al)  & \vert & f \: W' \to W \ \text{a morphism;} \\
              & &   \al \ \  \text{an automorphism of $f^*(x) \in \Ob\X(W')$}
                                        \end{array} \right\}$$}
 \noindent The forgetful map makes  $\mathbb{I}_x$ a
 sheaf of groups relative to $W_{pre}$. When $W_{pre}$ is a sheaf,
 then so is $\mathbb{I}_x$. This is always the case when the topology
 on $\C$ is subcanonical.

  \begin{defn}{\label{D:pi0}}
    Let $\X$ be a category fibered in groupoids over $\C$.
    We define $\pi_0\X$ to be the sheaf associated to the
    presheaf of sets on $\C$ defined by
    $$ W \mapsto \{\text{isomorphism classes in} \ \X(W)\}.$$
  \end{defn}

 There is a natural map $\X \to \pi_0\X$. This map is universal
 among  maps $\X \to Y$, where $Y$ is a sheaf over $\C$.
 For any $W$-point $x$ of $\X$,
 we can compose it with $\X \to \pi_0\X$ to obtain a $W$-point for $\pi_0\X$,
 which we denote by $\bar{x}$.

 \vspace{0.1in}
 The inertia presheaves $\mathbb{I}_x$ measure the difference between
 $\X$ and $\pi_0\X$ in the following sense.

 \begin{prop}{\label{P:inertiasheaf}}
    Let $x$ be a $W$-point of $\X$.
    Then we have a natural 2-cartesian diagram
        $$\xymatrix@=12pt@M=10pt{ \B_W\mathbb{I}_x  \ar[r] \ar[d] & \X \ar[d]  \\
                     W \ar[r]_(0.4){\bar{x}}        &     \pi_0\X      }$$
    \noindent where
    $B_W\mathbb{I}_x=[W/\mathbb{I}_x]$ is the classifying stack of the
    inertia presheaf $\mathbb{I}_x$ (see Definition \ref{D:quotient}).
 \end{prop}

  \begin{proof}
      Let  $X_0$ denote the presheaf
         $$ W \mapsto \{\text{isomorphism classes in} \ \X(W)\}.$$
      It follows formally from the definitions   that
      the following diagram is  2-cartesian
        $$\xymatrix@=12pt@M=10pt{
           \lfloor W/\mathbb{I}_x\rfloor   \ar[r]\ar[d]   & \X \ar[d]  \\
            W_{pre} \ar[r]        &      X_0   }$$
     The desired 2-cartesian diagram is the stackification of this diagram.
  \end{proof}

  \begin{defn}{\label{D:residue}}
    Let $\X$ be a stack over a site $\C$, and let $x \:W \to \X$ be a $W$-point.
    We define the {\bf residue stack} $\Gamma_x$ of $\X$ at $x$ to be the image
    of $x$, that is $\Gamma_x=\im(x)$.
 \end{defn}

 Using the definition of the image (Section \ref{SS:image})
 the following proposition is immediate

  \begin{lem}{\label{L:cartesian}}
      The following diagram is 2-cartesian:
           $$\xymatrix@=12pt@M=10pt{
                       \im(x) \ar@{^{(}->}[r] \ar[d]  &  \X \ar[d]  \\
                  \im(\bar{x})  \ar@{^{(}->}[r]  &  \pi_0\X       }$$
 \end{lem}

  \begin{cor}{\label{C:cartesian}}
    Notation being as above, assume $\bar{x} \: W \to \pi_0\X$ is a
    monomorphism. Then we have a    2-cartesian diagram
         $$\xymatrix@=12pt@M=10pt{ \Gamma_x \ar[r] \ar[d] &  \X \ar[d] \\
                   W  \ar[r]_(0.45){\bar{x}}        &     \pi_0\X      }$$
  \end{cor}

 Combining Lemma \ref{L:cartesian} with Proposition \ref{P:inertiasheaf},
 we get the following

  \begin{cor}{\label{C:residue}}
    We have a natural 2-cartesian diagram
        $$\xymatrix@=12pt@M=10pt{
                   B_W\mathbb{I}_x \ar[r] \ar[d] &  \Gamma_x\ar[d]  \\
                    W \ar[r]_(0.46){\bar{x}}        &     \pi_0\X      }$$
    In particular, when $\bar{x} \: W \to \pi_0\X$ is a monomorphism,
    we have a natural equivalence $B_W\mathbb{I}_x \risom \Gamma_x$.
  \end{cor}

\vspace{0.2in}

The inertia sheaves can be assembled in a more global object. For
a  stack  $\X$, we define its
 {\bf inertia stack}
 $\IX \to \X$, by the following 2-fiber product:

   $$\xymatrix@=12pt@M=10pt{  \IX \ar[r] \ar[d] &  \X \ar[d]^{\De}  \\
                 \X \ar[r]_(0.37){\De}        &    \X\x \X       }$$

For an object  $W$ in $\C$, the groupoid of sections $\IX(W)$ is
naturally equivalent to the following groupoid:

  $$\Ob = \left\{ \begin{array}{rcl}
           (x,\al)  & \vert & x \in \Ob\X(W);  \\
              & &   \al \ \  \text{an automorphism of $x$}
                                        \end{array} \right\}$$
 \vspace{0.1in}
  $$\Mor\big((x,\al),(x',\al')\big)=\left\{
         \gamma \in \Mor_{\X(W)}(x,x') \ \vert
                            \ \al\cdot\gamma=\gamma\cdot\al'     \right\}$$

\vspace{0.1in}
 The map $\IX \to \X$ corresponds the forgetful map.
 \vspace{0.1in}

 The inertia stack and the inertia presheaf of a $W$-point $x$
 are related by the following 2-cartesian diagram:
   $$\xymatrix@=12pt@M=10pt{ \mathbb{I}_x \ar[r] \ar[d] &  \IX \ar[d]  \\
                W \ar[r]_(0.45)x        &     \X      }$$

In particular,  the groupoid of possible liftings of $W$ to $\IX$
is equivalent to a set whose element are in natural bijection with
the group of self-transformations of $x$ (which is the same as the
inertia group $I_x$). This is in fact true if $W$ is any stack (or
category fibered in groupoids) over $\C$. In particular, when
$W=\X$, we find that the set of sections (up to identification) of
the map $\IX \to \X$ is in natural bijection with the group of
self-transformations of the identity functor $\id \: \X \to \X$.

\vspace{0.2in}
  We state the following proposition for future use. It is easy to prove.

  \begin{prop}{\label{P:mono}}
   Let $f \: \Y  \to  \X$ be a monomorphism.
   Then we have a natural 2-cartesian diagram
     $$\xymatrix@=12pt@M=10pt{ \IY \ar[r] \ar[d] & \IX \ar[d]  \\
                  \Y \ar[r]        &   \X       }$$
  \end{prop}

\subsection{Maps coming out of a quotient stack}{\label{SS:maps}}

Let $\X=[X/R]$ be the quotient stack of a groupoid object (in
presheaves of sets over $\C$)
 $[s,t \:R \sst{} X]$, and let $\Y$ an arbitrary stack.
The groupoid $\Hom_{\mathbf{St}}(\X,\Y)$ of maps from $\X$ to $\Y$
has a description in terms of cocycles. More precisely, define the
groupoid $C(\X,\Y)$ as follows:

{\small
       $$\Ob C(\X,\Y)=\left\{\begin{array}{rcl}
          (f,\varphi)
               & | & f\: X \to \Y, \ f\circ s \twomor{\varphi} f\circ t \\
               &   & \text{such that} \hspace{0.1in}
                   $$\xymatrix@M=4pt@R=8pt@C=-10pt{
                  f \circ d_0 \ar@{=>}[rr]^{\varphi\circ\mu}
                   \ar@{=>}[rd]_{\varphi\circ pr_1}   &   &  f \circ d_2  \\
                       &  f \circ d_1 \ar@{=>}[ur]_{\varphi\circ pr_2} &    }$$

                                            \end{array}\right\}$$
\vspace{0.1in}
  $$\Mor_{C(\X,\Y)}\big( (f,\varphi), (f',\varphi') \big)=
           \left\{\begin{array}{rcl}
              \al & | &  f \twomor{\al} f' \ \text{an identification}\\
                        &   &  \text{such that}  \hspace{0.2in}
         \xymatrix@R=10pt@C=10pt{
             f\circ s \ar @{=>}[r]^{\varphi}
                    \ar@{=>}[d]_{\al\circ s} \ar@{}[rd]|{\circlearrowleft}
                                         &  f\circ t \ar@{=>}[d]^{\al\circ t} \\
                    f'\circ s  \ar@{=>}[r]_{\varphi'}  &  f'\circ t   }
                                                          \end{array}\right\}$$}

 \noindent Here,  $d_0,d_1,d_2 \: R\x_X R \to X$ are the obvious maps,
         $pr_1,pr_2\: R\x_X R   \to R$ are the projection maps, and
         $\mu \: R\x_X R \to R$ is the multiplication.

 \begin{prop}
   The natural functor
   $$\Sigma \:\Hom_{\mathbf{St}}(\X,\Y) \to C(\X,\Y)$$
   is an equivalence of groupoids.
 \end{prop}

 \begin{proof}[Sketch of proof.] We show how to construct an
   inverse $\Theta$ for $\Sigma$.

   First, we claim that there is a presheaf of sets $Y$ over $\C$
   with an epimorphism $q \: Y \to \Y$ such that, for every presheaf
   of sets $X$ and every map $f \: X \to \Y$, there exists a
   lift $\tilde{f} \: X \to Y$. The presheaf $Y$ is defined as follows:
      $$W \mapsto \Ob\Hom(W_{pre},\Y).$$
   It is left to the reader to verify that $Y$ has the claimed property.

   Set $T=Y\x_{\Y}Y$. By Proposition \ref{P:quotient1},
   we have an equivalence $[Y/T]\risom \Y$.
   Now, given an object  $(f,\varphi)$ in $C(\X,\Y)$,
   pick a lift $\tilde{f}: X \to Y$ of $f$. The transformation
   $f\circ s \twomor{\varphi} f\circ t$ exactly translates as
   a map $R \to T$ (of presheaves of sets) for which the following
   diagram commutes:

     $$\xymatrix@=12pt@M=10pt{ R
           \ar[r] \ar[d] &  \ar[d] T \\
             X\x X \ar[r]        &    Y\x Y       }$$
  The cocycle condition satisfied by $\varphi$ also translates as saying that
  the map $R \to T$ makes  $[R\sst{}X] \to [T\sst{} Y]$
    a morphism of groupoids.

  The induced map on
  the quotient stacks gives us the desired map
  $\Theta(f,\varphi) \: \X \to \Y$. The effect of
  $\Theta$ on morphisms of $ C(\X,\Y)$ is defined in a similar way.
 \end{proof}

\section{Stacks over the topological site}{\label{S:topologicalsite}}

From now on, all stacks are over the site $\mathsf{Top}$ defined
below, unless otherwise specified.

 Let $\mathsf{Top}$  be the category of compactly generated
  topological with continuous  maps. We endow $\mathsf{Top}$ with a
 Grothendieck  topology  defined
 by taking the open coverings to be the usual open
 coverings of topological spaces.  This is a subcanonical topology.
 That means, the presheaf $W_{pre}$ represented by any $W \in \mathsf{Top}$
 is indeed a sheaf (so we denote it again by $W$).

We denote the 2-category of stacks over $\mathsf{Top}$ by
$\mathsf{St}_{\mathsf{Top}}$; it contains $\mathsf{Top}$ as a full
sub-category (Yoneda).

  Throughout the paper we fix a final object $*$ in $\mathsf{Top}$.

 \begin{defn}{\label{D:point}}
    By a {\bf point} $x$ of a stack $\X$ over $\mathsf{Top}$ we mean a
    $*$-valued point; that is,
    a map of stacks $x \: * \to \X$. We sometimes abuse notation and denote this
    by $x \in \X$.
 \end{defn}

\subsection{Representable maps}{\label{SS:representable}}

In $\mathsf{St}_{\mathsf{Top}}$, there are certain maps
 which behave much like maps
of topological spaces. In particular,
 many properties of continuous maps of topological
spaces makes sense for them.

  \begin{defn}{\label{D:representable1}}
   Let $f \: \Y \to \X$ be a morphism of stacks. We say $f$ is
   {\bf representable}, if for any map $X \to \X$ from a topological space
   $X$, the fiber product $Y:=X\x_{\X}\Y$ is (equivalent to) a topological space.
 \end{defn}

  \begin{prop}{\label{P:rep}}
       Let $f \: \Y \to \X$ be a representable morphism of stacks. Then,

        \begin{itemize}

          \item[$\mathbf{i.}$] for every stack $\Z$, the map
     $$\Hom_{\mathsf{St}}(\Z,\Y) \to \Hom_{\mathsf{St}}(\Z,\X)$$
     is  faithful. In particular, for every topological space $W$,
     $\Y(W) \to \X(W)$ is  faithful.

          \item[$\mathbf{ii.}$] for every topological space $W$,
           and every $W$ point
            $y \: W \to \Y$, the natural group homomorphism $I_y \to I_{f(y)}$
            is injective.
        \end{itemize}
    \end{prop}

\begin{proof}
Easy.
\end{proof}

  \begin{defn}{\label{D:local}}
    Let $\mathbf{P}$ be a property of morphisms of topological spaces.
      \begin{itemize}

    \item We say $\mathbf{P}$ is {\bf  invariant under base change}, if for
        every map $f \: Y \to X$  which is $\mathbf{P}$,
        the base extension $f' \: Y' \to X'$ of $f$ along any $q \: X' \to X$
        is again
        $\mathbf{P}$. We say $\mathbf{P}$ is {\bf invariant under restriction},
        if the above condition is satisfied when $q \: X' \to X$ is
        an embedding (i.e.
        inclusion of a subspace).

    \item We say $\mathbf{P}$ is {\bf local on the target},
        if for every map $f \: Y \to X$, the base extension
        $f' \: Y' \to X'$ of $f$ along a {\em surjective local homeomorphism}
         $q \: X' \to X$
        being $\mathbf{P}$ implies that $f$ is $\mathbf{P}$.
        ({\em Remark}. If we replace `surjective local homeomorphism'
        by `open covering' we get the same notion.
        If  $\mathbf{P}$ is invariant under restriction,
        we could also replace `surjective local homeomorphism'
        by  `epimorphism', arriving at the same notion.)
 \end{itemize}

    \end{defn}

 Remark that a map of topological spaces    is an epimorphism if and
 only if it admits local sections.

    \begin{defn}{\label{D:representable2}}
      Let $\mathbf{P}$ be a property of morphisms of topological spaces
      that is invariant under restriction and local on the target.
      We say a representable
      map $f \: \Y \to \X$ of stacks is $\mathbf{P}$, if
      there exists an epimorphism $q \: X \to \X$ such that
      the the base extension $Y \to X$ of $f$ along $q$
       is $\mathbf{P}$.
 \end{defn}

   If $\mathbf{P}$ in the above definition is also invariant under base
   extension, then it follows that for {\em every}
   $q \: X \to \X$ the base extension $Y \to X$ is $\mathbf{P}$.

    \begin{ex}{\label{E:P}}
    \end{ex}

        \begin{itemize}

    \item[$\mathbf{1.}$] Any of the following is invariant under base change
          and  local on the target:
    \vspace{0.1in}
    \begin{quote}
          {\em open, epimorphism, surjective, embedding }
          (see below for definition), {\em closed embedding, open embedding,
          local homeomorphism, covering map, finite fibers, discrete fibers.}
    \end{quote}

 \vspace{1mm}
    \item[$\mathbf{2.}$] Any of the following is invariant under restriction
          and  local on the target:
    \vspace{0.1in}
     \begin{quote}
        {\em closed,  closed onto image}.
      \end{quote}

 \end{itemize}

 \vspace{0.1in}

   \begin{defn}{\label{D:embedding}}
     We say a map
  $f \: Y \to X$ of topological spaces is an {\bf embedding},
  if $f$ induces a homeomorphism from $Y$ onto the subspace $f(Y)$ of $X$.
  When the image of $f$ is open, we say $f$ is an {\bf open embedding}.
  When the image of $f$ is closed, we say $f$ is a {\bf closed embedding}.
  A composite of open and closed embeddings is called a {\bf locally-closed
  embedding}.
  All these notions are invariant under base change and local on the
  target, so we can define
  them for representable  maps of stacks (Definition \ref{D:representable2}).
 \end{defn}

  \begin{lem}
     Every embedding is a monomorphism.
   \end{lem}

    \begin{proof}
       Easy.
    \end{proof}

  \begin{defn}{\label{D:embedded}}
    Let $\X'$ be a substack of $\X$. We say that $\X'$ is an {\bf embedded
    substack} if the inclusion map $\X' \hra \X$ is an embedding.
    {\bf Open substack} and {\bf closed substack} are defined
    in a similar fashion.
  \end{defn}

 Note that the inclusion of a substack $\X' \hra \X$ is not necessarily
 representable, so not every substack is embedded.

\subsection{Operations on substacks}{\label{SS:operations}}
  We will briefly go over the notions of
  {\em intersection}, {\em union}, {\em inverse image} and {\em image}
  for substacks of a stacks over $\mathsf{Top}$. What we are really
  interested in  the case where substacks are {\em embedded}.

\vspace{0.1in}
  Let $\X$ be a stack over $\mathsf{Top}$, and let
  $\{\X_i\}_{i\in I}$  be a family of substacks.
  The {\bf intersection} $\bigcap_{i\in I}\X_i$
  is defined to be the intersection of the categories $\X_i$
  in the set theoretic
  sets (note that substacks are assumed to be saturated and full).
  It is easy to see that the intersection of substacks is again a substack.

  Let $\X'$ be a substack of $\X$ and let $\Y \to \X$ be a morphism
  of stacks. We define the {\bf inverse image} $\Y'$ of $\X'$ to
  be the image (Definition \ref{D:image})
  of the natural monomorphism $\Y\x_{\X}\X' \to \Y$.
  The inverse image is again a substack.

   \begin{lem}{\label{L:operations}}
    \end{lem}
  {\em
     \begin{itemize}
          \item[$\mathbf{i.}$] Taking inverse image commutes with
            intersection.

          \item[$\mathbf{ii.}$] Intersection of {\em embedded} substacks is again
            embedded. Intersection of closed substacks is again a closed
            substack.

          \item[$\mathbf{iii.}$] Inverse image of an embedded (resp. closed,
              open) substack under an
            arbitrary morphism of stacks is again an
             embedded (resp., closed, open) substack.

        \end{itemize}
    }
   \begin{proof}
   Straightforward.
    \end{proof}

  Part  ($\mathbf{ii}$) of the lemma implies that, for any substack
  $\X'$ of a stack $\X$, there is a canonical (unique)  smallest embedded
  substack $\X''$ of $\X$ that contains $\X'$. In general, given
  a collection of substacks $\{\X_i\}_{i \in I}$ of $\X$,
  there is a canonical (unique)
  smallest embedded substacks of $\X$ containing all the $\X_i$.

  Similarly, given a closed substack $\X'$ of a stack $\X$, we define
  its {\bf closure} $\bar{\X'}$ to be the intersection
  of all closed substacks containing $\X'$. The closure $\bar{\X'}$
  is a closed substack of $\X$.

   The {\bf union} $\bigcup_{i\in I}\X_i$
   of a family $\{\X_i\}_{i \in I}$
   of embedded substacks of $\X$ is defined
   to be the smallest  embedded substack
   of $\X$ containing all the $\X_i$. An alternative description
   for the union of substacks is given by the following lemma.

    \begin{lem}{\label{L:union}}
       An object $x \in \X(W)$ is in $\bigcup_{i\in I}\X_i(W)$
       if and only if the corresponding map $W \to \X$ has the property
       that $\bigcup_{i\in I}W_i=W$, where $W_i$ is the inverse image
       of $\X_i$ in $W$ (which is a subspace of $W$).
    \end{lem}

   \begin{lem}{\label{L:operations2}}
    \end{lem}
  {\em
     \begin{itemize}
          \item[$\mathbf{i.}$] Taking inverse image commutes with
            union.

          \item[$\mathbf{ii.}$] Union of open substacks is again an open
            substack.

        \end{itemize}
    }
   \begin{proof}
   Follows from Lemma \ref{L:union}.
    \end{proof}

   Finally, given a map $f \: \Y \to \X$ of stacks, we define the {\bf embedded
   image}
   $f(\Y)$ of $f$ to be the smallest embedded substack of $\X$ containing
   $\im(f)$. Embedded image corresponds to the usual notion of
   image of maps of topological spaces. It is typically bigger than $\im(f)$
   (for example, take $f$ to be
   a surjective non epimorphic map of topological spaces).

   When all the stacks involved are topological spaces, the
   above notions coincides with the original  notions
   in classical topology.

\subsection{The coarse moduli space}{\label{SS:coarse}}

 To an stack $\X$ we associate a topological
 space $\Xm$, called
 the {\em coarse moduli space} of $\X$ (or, loosely speaking,
 the {\em underlying space} of $\X$). There is a natural
 map $\pim \: \X \to \Xm$ which is universal among maps from $\X$ to topological
 spaces.

 The coarse moduli space $\Xm$ is defined as follows. As a set it is equal to
 $\pi_0(X(*))$, that is, the set of maps from $*$ to $\X$, up to
 identifications.
 For any open substack $\U \subseteq \X$, we have a natural inclusion
 $\U_{mod} \subseteq \Xm$. These are defined to be the open sets of $\Xm$.

  \begin{ex}{\label{E:quotient}}
   Let $X=[R \sst{} X]$ be a groupoid in topological spaces, and
   let $\X=[X/R]$ be the quotient stack.
   Then we have a natural homeomorphism $\Xm \cong X/R$, where the
   latter is the (naive) quotient space of the equivalence
   relation on $X$ induced by  $R$.
  \end{ex}

  \begin{rem}{\label{R:pi0}}
    The coarse moduli space $\Xm$ of a stack $\X$ should not be confused with
    $\pi_0\X$ of Section \ref{SS:inertia1}  (Definition \ref{D:pi0}).
    All we can say is that there is a map $\pi_0\X \to \Xm$
    which induces a bijection on the set of points (i.e. $(\pi_0\X)(*) \to \Xm(*)$
    is a bijection). This map is neither a monomorphism nor an epimorphism
    in general (also see Corollary \ref{C:gerbe1}).
  \end{rem}

  \begin{prop}{\label{P:functoriality}}
     Let $f \: \X \to \Y$ be a map of stacks.
   \begin{itemize}

    \item[$\mathbf{i.}$]The   induced map $f_{mod} \: \Xm \to \Ym$
     is continuous. If $g \: \X \to \Y$ is another map that is
     2-isomorphic to $f$, then $f_{mod}=g_{mod}$.

    \item[$\mathbf{ii.}$]  {\bf Functoriality}.
    There is a natural map $\pim \: \X \to \Xm$.
     This map is functorial, in the sense that the following diagram
     commutes:
         $$\xymatrix@=12pt@M=10pt{  \X   \ar[r]^f \ar[d]_{\pim}
                                               &  \Y \ar[d]^{\pim}  \\
                      \Xm  \ar[r]_{f_{mod}}    &  \Ym         }$$

    \item[$\mathbf{iii.}$]{\bf Universal property}.
     For any map $f \: \X \to Y$ to a topological space,
     there is a unique continuous map $f_{mod} \: \Xm \to Y$ which makes
     the following diagram commute:
          $$\xymatrix@=12pt@M=10pt{  \X \ar[d]_{\pim} \ar[r]^f   &   Y     \\
                        \Xm  \ar[ru]_{f_{mod}}       &          }$$

   \end{itemize}
 \end{prop}

 \begin{proof}[Proof of part {\em ($\mathbf{i}$)}]
  This is true since the inverse image
  of an open substack is again an open substack.

  \vspace{0.1in}
  \noindent{\em Proof of part} ($\mathbf{ii}$). Let $W \in \mathsf{Top}$
  and let $a \in \Ob\X(W)$. We abuse notation and denote the corresponding
  map   $W \to \X$ also by $a$. We define
  $\pim(a) \: W \to \Xm$ to be the map that sends $w \in W$ to $a(w) \in \Xm$.
  In other words, $\pim(a)$ is obtained from $a \: W \to \X$ by passing
  to the coarse moduli spaces, which is continuous by Part ($\mathbf{i}$).
  Again, by Part ($\mathbf{i}$),  2-isomorphic element of $\X(W)$ gives
  rise to the same elements in $\Xm(W)$. The commutativity of the diagram
  is obvious.

  \vspace{0.1in}
  \noindent{\em Proof of part} ($\mathbf{iii}$). Follows from Part ($\mathbf{i}$).
 \end{proof}

   \begin{lem}{\label{L:embedded}}
      Let $\X' \subseteq \X$ be an embedded substack. Assume
      the induced map $\Xm' \to \Xm$ is a bijection.
      Then $\X'=\X$.
    \end{lem}

    \begin{proof}
       We have to show that for every topological space $W$
       we have $\X'(W)=\X(W)$. If this is not the case, we can find a
       $W$ and a map $f \: W \to \X$ such that $W':=f^{-1}(\X') \subset W$
       is not equal to $W$. Pick a point $w \in W$ that is not in $W'$.
       Then $f(w)$ is a point of $\X$ that is not in $\X'$. This
       contradicts the fact that $\Xm' \to \Xm$ is a bijection.
    \end{proof}

    \begin{prop}{\label{P:embedded}}
       Taking inverse image gives a bijection between the subspaces
       of $\Xm$ and embedded substacks of $\X$. The inverse is given by
       taking embedded image (see Section \ref{SS:operations}).
    \end{prop}

   \begin{proof}
      Follows from Lemma \ref{L:embedded}.
    \end{proof}

     \begin{cor}{\label{C:embedded}}
       Let $f \: \Y \to \X$ be a map of stacks. Then
       $f(\Y)=\pim^{-1}\big(f_{mod}(\Ym)\big)$, where $f(\Y)$
       stands for the embedded image of $f$ (Section \ref{SS:operations}).
    \end{cor}

\section{Gerbes}{\label{S:gerbes}}

  In this short section, we quickly review gerbes over $\mathsf{Top}$.
  The following proposition also serves as a reminder of the definition
  of a gerbe.

  \begin{prop}{\label{P:relativegerbe}}
    Let $f \: \X \to \Y$ be a morphism of stacks. The following conditions
    are equivalent:

        \begin{itemize}

          \item[$\mathbf{i.}$] For any topological space $W$
          and  any $y \in \Y(W)$,
          there is an open covering $\{U_i\}$ of $W$ such that each $y|_{U_i}$
          is equivalent (in the groupoid $\Y(U_i)$) to
          $f(x_i)$, for some $x_i \in \X(U_i)$. For
          $x, x' \in \X(W)$, and for any  isomorphism
          $\be \: f(x) \twomor{} f(x')$ in $\Y(W)$, there is an open  covering
          $\{U_i\}$ of $W$ such that each restriction $\be|_{U_i}$
          can be lifted to an isomorphism (in the groupoid $\X(U_i)$)
          $\al_i \: x|_{U_i} \twomor{} x'|_{U_i}$.

          \item[$\mathbf{ii.}$] The maps $f \: \X \to \Y$ and
          $\De \: \X \to \X\x_{\Y}\X$ are epimorphisms.

          \item[$\mathbf{iii.}$] For any map $g \: \Z \to \X$ of stacks,
          $g$ is an epimorphism if and only if $f\circ g$ is so.

          \item[$\mathbf{iv.}$]  For any map $g \: Z \to \X$ from a
          topological space $Z$ to $\X$,
          $g$ is an epimorphism if and only if $f\circ g$ is so.
        \end{itemize}
  \end{prop}

   \begin{proof}
     Standard.
   \end{proof}

  \begin{defn}{\label{D:relativegerbe}}
    Let $f \: \X \to \Y$ be a morphism of stacks. We say that
    $\X$ is a {\bf relative gerbe}  over $\Y$ (via $f$) if the equivalent
    conditions of Proposition \ref{P:relativegerbe} are satisfied.
    When $\Y$ is a topological space, we simply say that $\X$ is a {\bf gerbe}.
    This is justified by the next proposition.
  \end{defn}

  \begin{prop}{\label{P:gerbe}}
  \end{prop}
  {\em
        \begin{itemize}

          \item[$\mathbf{i.}$]
          Let $\X$ be a stack. Then $\X$ is a gerbe relative to $\pi_0\X$.
          Conversely,    suppose $\X$ is a gerbe relative to a sheaf $Y$.
          Then there is a natural isomorphism of sheaves $\pi_0\X \risom Y$.

          \item[$\mathbf{ii.}$] Let $\X$ be a gerbe relative to a
          topological space $Y$. Then there is a unique homeomorphism
            $\Xm \cong Y$ (respecting the map coming from $\X$).

        \end{itemize}
   }

  \begin{proof}[Proof of part {\em ($\mathbf{i}$)}]

     This is more or less the definition (as in Proposition
     \ref{P:relativegerbe}.$\mathbf{i}$).

\vspace{0.1in}
   \noindent{\em Proof of part} ($\mathbf{ii}$). By Part ($\mathbf{i}$), we have
   $Y=\pi_0\X$, so $Y$ satisfies
   the universal property of $\pi_0\X$
   (the paragraph after Definition \ref{D:pi0}). In particular, $Y$ satisfies
   the universal property of $\Xm$
   (Proposition \ref{P:functoriality}.$\mathbf{iii}$),
   so it is uniquely homeomorphic to $\Xm$.
  \end{proof}

  \begin{cor}{\label{C:gerbe1}}
    Let $\X$ be a stack. Then $\X$ is a gerbe if and only if the natural map
    $\pi_0\X \to \Xm$ is an isomorphism.
  \end{cor}

   \begin{cor}{\label{C:gerbe2}}
   Let $\X$ be a gerbe
   and let $x\: W \to \X$ be a $W$-point in $\X$.
   Then we have a natural 2-cartesian diagram

     $$\xymatrix@=12pt@M=10pt{
         B_W\mathbb{I}_x    \ar[r] \ar[d] &  \X \ar[d]  \\
           W      \ar[r]        &  \Xm         }$$
   \end{cor}

   \begin{proof}
     This follows from Corollary \ref{C:gerbe1} and
     Proposition \ref{P:inertiasheaf}.
    \end{proof}

  \begin{cor}{\label{C:trivial}}
    Let $\X$ be a gerbe (over  $\Xm$), and assume $\pim \: \X \to \Xm$
    has a section $s \: \Xm \to \X$. Then we have a natural equivalence
    $\X\cong B_{\Xm}\mathbb{I}_s$.
  \end{cor}

 A gerbe $\X$ is called {\bf trivial} if $\pim \: \X \to \Xm$
 has a section. By Corollary \ref{C:trivial}, every trivial gerbe
 is of the form $B_WG$, where $W$ is a topological space and $G$ is a
 sheaf of groups over $W$. In general, for any gerbe $\X$, the map
 $\X \to \Xm$ has a section after replacing $\Xm$ with an open cover.
 Therefore, any gerbe $\X$ can be covered by open substacks of the form
 $B_WG$.

 For more examples of gerbes see \ref{SS:gerbes}.


\section{A useful criterion for representability}{\label{S:criterion}}

 All stacks are over the site $\mathsf{Top}$.

 This section concerns a useful technical lemma for proving
 representability of maps of stacks (Lemma \ref{L:representable}).
 Before proving the lemma, we need some preliminary results.

  \begin{lem}{\label{L:epimorphism}}
   Let $f\: \Y \to \X$ be a morphism of stacks over a base category $\C$.
   Suppose there exists an epimorphism  $\X' \to \X$
   such that the base extension $\X'\x_{\X}\Y \to \X'$
   is an equivalence. Then  $f$ is also an equivalence.
 \end{lem}

  \begin{proof}
    To show that $f$ is an equivalence we have to show that, for
    every $X \in \C$, the induced map $\Y(X) \to \X(X)$ is an equivalence
    of groupoids. By formal nonsense, it is enough to show that,
    for every map $p \: X \to \X$, the base extension of $f$ along $p$,
    which itself satisfies the condition of the lemma,
    is an equivalence.  So we are reduced to the case $\X=X$, for some
    $X \in \C$. Since for every stack
    $\X'$ there exists an epimorphism
    $X' \to \X'$ from a sheaf of sets $X'$, we may also assume that
    $\X'=X'$, for some sheaf of sets $X'$. So, we end up with a
    cartesian diagram
        $$\xymatrix@=12pt@M=10pt{ X' \ar[r] \ar[d]_= & \Y \ar[d]^f  \\
                     X'\ar[r]        &   X         }$$
    where $X' \to X$ is an epimorphism of sheaves   of sets.
    We have to show that $f$ is an isomorphism, that is, for every $W \in \C$,
    the groupoid $\Y(W)$ is equivalent to the set $X(W)$ via $f$.

    First, we claim that, for every  pair of $W$-points $a$ and $b$ in $\Y(W)$
    such that $f(a)=f(b)$, there is a {\em unique} isomorphism in $\Y(W)$
    between $a$ and $b$. If $a$ and $b$
    are both in the image of $X'$, this follows from the fact the above
    diagram is cartesian.  In the general case, we can replace $W$ by an open
    covering over which $a$ and $b$ are both in the image of $X'$,
    so on each of the opens in the covering, we find a unique isomorphism
    between the restrictions of $a$ and $b$. Since $\Y$ is a stack, these
    isomorphisms  glue to an isomorphism between $a$ and $b$ over $W$.
    This proves the claim.

    This claim implies two things at the same time. The first
    one is that $\Y(W)$ is equivalent to a set, so we might as well
    assume $\Y$ is a sheaf of sets, and switch the notation from $\Y$ to $Y$.
    The second thing is that the   map $f \: Y \to X$ of   sheaves
    of sets is injective.

    All that is left to check is that $f$ is an
    epimorphism (note that we are dealing with {\em sheaves}).
    But this is obvious from the commutativity of the above
    diagram,  because if we pre-compose $f$ with $X' \to X$ we get $X' \to X$,
    which is an epimorphism.
  \end{proof}

  \begin{cor}{\label{C:locallyspace}}
    Let $\X$ be a stack, and let $\{\U_i\}_{i \in I}$ be a covering of
    $\X$ by open substacks (Definition \ref{D:embedding}).
    (Here, by covering we mean that the map $\coprod \U_i \to \X$ is an
    epimorphism, or, equivalently, the induced map $\coprod \U_{i,mod} \to \Xm$
    of topological spaces, which is just a union of open embeddings,
    is surjective.)
    Assume each $\U_i$ is equivalent to a topological space. Then, so is $\X$.
  \end{cor}

  \begin{proof}
    Consider the map of $\pi_{mod} \: \X \to \X_m$. By assumption,
    the base extension  of this map via the epimorphism
    $\coprod \U_i^{mod} \to \Xm$ is an equivalence of stacks. The result
    follows from Lemma \ref{L:epimorphism}.
  \end{proof}

 Now, we come to the main lemma.

 \begin{lem}{\label{L:representable}}
    Consider a 2-cartesian diagram of stacks
    $$\xymatrix@=12pt@M=10pt{ \X' \ar[r]^a \ar[d]_{f'} &  \X \ar[d]^f  \\
                 \Y' \ar[r]_b        &       \Y   }$$
    \noindent in which the horizontal arrows are epimorphisms.
    If $f'$ is representable, then so is $f$.
 \end{lem}

  \begin{proof}
    We may assume $\Y$, $\X'$ and $\Y'$ are  representable (hence denoting them
    by $Y$, $X'$ and $Y'$, respectively). We have to show that $\X$ is
    representable. First we assume that $b$ is a disjoint union of open
     embeddings, say $Y'=\coprod V_i \to Y$.  Let $X'=\coprod U_i$ be the
     corresponding decomposition for $X'$. Then $\coprod U_i \to \X$
     is a covering of
     $\X$ by open substacks, each of which equivalent to a topological space.
     So, by Corollary \ref{C:locallyspace}, $\X$ itself is equivalent
     to a topological
     space.


     For the general case, by what we just proved,
     we may replace $Y$ with an open covering, so
     we may assume $b \: Y' \to Y$ has a section. But in this case the result
     is obvious. Proof is complete.
  \end{proof}

\section{Pretopological stacks}{\label{S:pretopological}}

 In the previous section we looked at stacks over the topological
 site $\mathsf{Top}$. These are simply categories fibered in groupoids
 over $\mathsf{Top}$ (or, loosely speaking, presheaves of groupoids
 over $\mathsf{Top}$) which satisfy the descent condition. We saw that
 we can pinpoint  a certain class of morphisms of such stacks
 ({\em representable} morphisms)   to which we
 could attribute usual  properties of maps of
 topological spaces.

   However, stacks over $\mathsf{Top}$,  are still too crude to do topology
  on. We need to be able to approximate such a stack by an honest
  topological space so as to, possibly, be able to talk about its topological
  properties. The next definition is our first approximation
  of a reasonable notion of a {\em topological} stack. The full-fledged
  definition (i.e. that of a {\em topological stack}) will be given
  in Section \ref{S:TopSt}.

   \begin{defn}{\label{D:stackwithchart}}
    Let $\X$ be a stack over $\mathsf{Top}$. A {\bf chart} for
    $\X$ is  a representable epimorphism $p \: X \to \X$ from
    a topological space $X$ to $\X$. If such a chart exists,
    we say $\X$ is a {\bf pretopological stack}.
   \end{defn}

  First we  notice that the condition on
 representability of the diagonal which appears in standard texts on
 stacks is implied by the definition.

  \begin{prop}{\label{P:diagonalrep}}
   Let $\X$ be a pretopological stack.  Then,
   the diagonal $\De \: \X \to \X\x\X$ is representable.
  \end{prop}

      \begin{proof}

     Let $p \: X \to \X$ be a   chart for $\X$. Base extend $\De$ by
     $(p,p) \: X \x X \to \X\x\X$. We obtain the following cartesian diagram:
         $$\xymatrix@=12pt@M=10pt{
                   X\x_{\X} X \ar[r]  \ar[d] &  \X \ar[d]^{\De}  \\
                      X\x X  \ar[r]_{(p,p)}       &   \X\x\X       }$$
     \noindent  Since $X\x_{\X} X$ is representable and $(p,p)$ is an
     epimorphism,  Lemma \ref{L:representable} applies.
   \end{proof}

  \begin{cor}{\label{C:representable}}
   Let $f \: X \to \X$ be a map from a topological space to a
   pretopological stack. Then $f$ is representable.
  \end{cor}

 \begin{cor}{\label{C:fiber}}
    The 2-category of pretopological stacks is closed under 2-fiber
    product. (Also see Section \ref{S:fiber}.)
 \end{cor}

 \begin{proof}
   Let $\Y \to \X$ and $\Z \to \X$ be morphisms of pretopological
   stacks. Let $Y \to \Y$ and  $Z \to \Z$ be charts. Then it is
   easy to check that $Y\x_{\X}Z$, which is a topological space
   by Corollary \ref{C:representable}, is a chart for
   $\Y\x_{\X}\Z$.
 \end{proof}

    \begin{cor}{\label{C:image}}
       Let $f \: \Y \to \X$ be a morphism of pretopological stacks.
       Then $\im(f)$ is pretopological (but not necessarily
       embedded).
    \end{cor}

      \begin{proof}
       Denote $\im(f)$ by $\I$.
         Let $p \: X \to \X$ be a chart for $X$. Then, the induced
         map $X \to \im(f)$ is an epimorphism; so, $\im(f\circ p)=\I$.
         Let $R:=X\x_{\I}X=X\x_{\Y}X$. Since $\Y$ is pretopological,
         Corollary \ref{C:representable} implies that $R$ is an honest
         topological space. Therefore, $X \to \I$ is representable,
         because its base extension along the epimorphism $X \to \I$
         is the map $R \to X$, which is representable
         (Lemma \ref{L:representable}).
      \end{proof}

 \begin{cor}{\label{C:inertia}}
   Let $\X$ be a pretopological stack. Then the inertia stack $\IX$
   is representable over $\X$ (i.e. the map $\IX \to \X$ is representable).
 \end{cor}

  \begin{proof}
    This follows from Proposition \ref{P:diagonalrep},
    because $\IX \to \X$ is the base extension
    of the representable map
    $\X \llra{\De} \X\x\X$ along  $\X \llra{\De} \X\x\X$.
 \end{proof}

 Let $\X$ be a pretopological stack, and let $p \: X \to \X$ be a chart for it.
 Set $R=X\x_{\X}X$. Since $p$ is representable, $R$ is again a topological
 space. Thus, we obtain a groupoid in topological spaces $[R \sst{} X]$.
 By Proposition \ref{P:quotient1}, we have an equivalence $[X/R] \to \X$.
 Conversely, given a groupoid in topological space $[R \sst{} X]$,
 we can form  a quotient stack $\X=[X/R]$. The natural map $P \: X \to \X$
 is then a representable (Lemma \ref{L:representable}) epimorphism.
 Therefore, $\X$ is indeed a pretopological stack.

 The above correspondence can be made precise by saying that, the
 2-category of pretopological stacks equipped with a chart, is equivalent
 to the 2-category of  groupoids in topological spaces.
 This  correspondence gives rise to a correspondence between
 pretopological stacks (this time without  a fixed chart) and
 {\em Morita equivalence classes} of  topological groupoids.
 We will say a few words on this in Section \ref{S:Morita}.

 Under the above correspondence,
 any property of the source (and target) map of a groupoid $[R \sst{} X]$
 is reflected as a property of the chart $p \: X \to [X/R]$
  (assuming this property is invariant under base change).
 For instance, the source and target maps of $[R \sst{} X]$ are open, \'{e}tale,
 fibration  etc. if and only if the chart
 $p \: X \to [X/R]$ is so.

 When $G$ is a topological group acting on a topological space
 (or, more generally, a presheaf of groups acting on a presheaf
 of sets), we can form the groupoid $[G\x X \sst{} X]$ where
 the source and target maps are the projection map and the action map,
 respectively. This groupoid is called the {\bf  action groupoid} of this
 action. The following simple lemma tells us when a topological
 groupoid is the action groupoid of a discrete group action.

  \begin{lem}{\label{L:actiongpd}}
   Let $[R \sst{} X]$ be a topological groupoid.  Assume $X$ is connected
   and
   $R$ is a disjoint union of components, each of which
   mapping homeomorphically to $X$ via source and target maps.
   Then $[R \sst{} X]$ is the action groupoid of a discrete group
   acting on $X$.
  \end{lem}

   \begin{proof}
    Easy.
    \end{proof}

  The {\bf stabilizer group} of a groupoid $[R \sst{s,t} X]$ is defined to be
  the relative (topological) group $\mathbb{I}_X \to X$, where $\mathbb{I}_X$
  is defined by the following cartesian diagram:

    $$\xymatrix@=12pt@M=10pt{  \mathbb{I}_X \ar[r] \ar[d] &  R \ar[d]^{(s,t)}   \\
                   X \ar[r]_(0.37){\De}        &    X\x X       }$$

 If $\X=[X/R]$, for some topological groupoid $[R \sst{} X]$,
 we have the following cartesian diagram

   $$\xymatrix@=12pt@M=10pt{ \mathbb{I}_X \ar[r] \ar[d] &  \IX \ar[d]  \\
                X \ar[r]        &  \X        }$$
where $\IX$ is the inertia stack of $\X$ (Section
\ref{SS:inertia1}).

 If we denote the chart $X \to [X/R]$ be $p$, then (the sheaf
 represented by $\mathbb{I}_X$) can be identified with what
 we would call $\mathbb{I}_p$ by the terminology of Section \ref{SS:inertia1}.

 The fiber of $\mathbb{I}_X \to X$ over a point $x$ in $X$ is isomorphic,
 as a topological group, to the
 inertia group $\mathbb{I}_{p(x)}$ at $p(x)$, so $\IX$
 can be thought of as organizing
 all the inertia groups in one global family. We will
 see in  Section \ref{S:homotopy} that an element in a group $\mathbb{I}_x$
 can be viewed
 as a ``ghost loop'' at $x$ which is not visible if we only look at it in
 the coarse moduli space $\Xm$ of $\X$ (Section \ref{SS:coarse}).
 So $\IX$ is parameterizes all these ghost loops.

 The inertia stack of a topological space is trivial (i.e $\IX\cong\X$),
 but the converse is
 not
 true (Example \ref{E:vistoli} below).
 We call  a pretopological stack  with trivial inertia stack
 a {\bf quasitopological space} or
 {\bf quasirepresentable}. A pretopological stack $\X$ is a quasitopological
 space if and only if for every topological space $W$, the groupoid of sections
 $\X(W)$ is equivalent to a set (so $\X$ is really a sheaf of sets).
 Similarly one can define the notion
 of a {\bf quasirepresentable map} of stacks.

 One can produce examples of pretopological spaces that are not a topological
 space by considering group actions that are not
 properly discontinuous. I learned the following example from
 Angelo Vistoli.

  \begin{ex}{\label{E:vistoli}}
  \end{ex}
   \begin{itemize}

        \item[$\mathbf{1.}$]
          Let  $\mathbb{Q}$, viewed as a discrete group,
          act on $\mathbb{R}$ by translations.
          Then $[\mathbb{R}/\mathbb{Q}]$ is a quasitopological
          space (since the action
          has no fixed points), but it
          is not a topological space.
          An easy way to see this is to compare the fiber products
          $\mathbb{R}\x_{{\mathbb{R}}/{\mathbb{Q}}}\mathbb{R}$ and
          $\mathbb{R}\x_{[{\mathbb{R}}/{\mathbb{Q}}]}\mathbb{R}$. The former
          is homeomorphic to  $\mathbb{R}\x\mathbb{Q}$, where $\mathbb{Q}$
          is endowed with its subspace topology from $\mathbb{R}$. But the
          latter is homeomorphic to  $\mathbb{R}\x\mathbb{Q}$, where $\mathbb{Q}$
          is endowed with the discrete topology.

        \vspace{0.1in}
        \item[$\mathbf{2.}$] Consider the previous example, but now assume
          that  $\bbQ$ is endowed with the subspace topology from $\bbR$.
          Again, the quotient stack $[\mathbb{R}/\mathbb{Q}]$ is
          a quasitopological space, but not a topological space (see the next
          proposition).

        \end{itemize}

 \vspace{0.1in}
    The following proposition tells us exactly when a quotient stack is a
    topological space.

    \begin{prop}
          Let $[R \sst{} X]$ be a topological groupoid. Then
          the quotient stack $[X/R]$ is a quasitopological
          space if and only if the map $R \to X\x X$ is injective.
          It is a topological space if and only if the following
          two conditions hold:
         \begin{itemize}

          \item $R \to X\x X$ is an embedding;

          \item $X \to X/R$ is an epimorphism.

        \end{itemize}
    In this case, we have $[X/R]\cong X/R$.
    \end{prop}

\section{A few words on Morita equivalence}{\label{S:Morita}}

 We saw in the previous section that a  topological groupoid
 gives rise to a pretopological stack by the   quotient stack
 construction. Several topological
 groupoids may lead to equivalent quotient stacks. Such groupoids
 are called {\bf Morita equivalent}.

 Since the language of groupoids and Morita equivalences is frequently
 used in the literature, we briefly mention how it fits in the stacky point of
 view.

 \vspace{0.1in}
 All groupoid are assumed to be topological groupoids.
 \vspace{0.1in}

 Let $[R' \sst{} X']$ and $[R \sst{}  X]$ be groupoids, and let
 $F \: [R' \sst{} X']\to[R \sst{}  X]$ be a map of groupoid.
 By this we mean a pair of maps $F_0 \: X' \to  X$ and $F_1 \: R' \to R$
 satisfying the obvious conditions. This induces a map of quotient
 stacks $f \: [X'/R'] \to [ X/R]$.
 We say $f$ is an {\em elementary Morita
 equivalence} if $F_0 \: X' \to  X$ is an epimorphism and the
 induced map $R' \to R|_{X'}=(X'\x X')\x_{ X\x  X} R$ is a
 homeomorphism.\footnote{We can alternatively define
 elementary Morita equivalence
 by requiring $F_0$ to be an open covering. Everything in this
 section remains valid, except
 we have to modify the statement of Proposition \ref{P:Morita}.}

 It is easy to show that, when $F$ is an elementary Morita equivalence,
 the induced map $f \: [X'/R'] \to [X/R]$ of quotient stacks is
 an equivalence. The converse is not in general true. That is,
 if $F \: [R' \sst{} X']\to[R \sst{}  X]$ is a map of groupoids
 such that the induces map $f \: [X'/R'] \to [ X/R]$ is an equivalence,
 it does not necessarily follow that $F$ is an elementary Morita equivalence.
 But we have the following result.

  \begin{prop}{\label{P:Morita}}
    Let $F \: [R' \sst{} X']\to[R \sst{}  X]$ be a map of groupoids.
     Assume in addition that $[R \sst{}  X]$ is an \'{e}tale groupoid
     (that is, the source and target maps are local homeomorphisms).
     Then the induced map $f \: [X'/R'] \to [ X/R]$ of quotient stacks is
    an open embedding if an only if $U=f(X')$ is an open subset of $ X$
    and the induced map $[R' \sst{} X']\to[R|_U \sst{} U]$ is an elementary
    Morita equivalence. Furthermore, $f$ is an equivalence if and only if
    $\OO(U)= X$.
  \end{prop}

  We will not use this result in this paper and leave it as an exercise to the
  reader.

 \vspace{0.1in}

 Given two groupoids $[R' \sst{} X']$ and $[R \sst{}  X]$ whose
 quotient stack are equivalent, it is not necessarily true that
 there is map $F \: [R' \sst{} X']\to[R \sst{}  X]$ inducing the
 equivalence. All we can say is that, there is another
 groupoid $[R'' \sst{}  X'']$, and a pair of elementary
 Morita equivalences as in the diagram

  {\footnotesize $$\xymatrix@=12pt@M=10pt@L=6pt{
          & \ar[ld]_(0.6){\text{Elem. Morita}} [R'' \sst{}  X'']
                             \ar[rd]^(0.63){\text{Elem. Morita}} & \\
                                  [R' \sst{} X'] &  &  [R \sst{}  X]   }$$}

 More generally, given a morphism $f \: [X'/R'] \to [ X/R]$,
 we can find a groupoid $[R'' \sst{}  X'']$, an elementary
 Morita equivalence $[R'' \sst{}  X''] \to [R' \sst{} X']$, and a
 map of groupoids $F \: [R'' \sst{}  X'']\to[R \sst{}  X]$
 that induces $f$ after passing to quotient stacks:
   {\footnotesize $$\xymatrix@=12pt@M=10pt@L=6pt{
            & \ar[ld]_(0.6){\text{Elem. Morita}} [R'' \sst{}  X'']
                                   \ar[rd]^F & \\
                                   [R' \sst{} X'] &  &  [R \sst{}  X]   }$$}

 \noindent In fact, given a finite collection  of maps
 $f_1,f_2,\cdots,f_n \: [X'/R'] \to [ X/R]$, we can choose  $[R'' \sst{}  X'']$
 so that it works simultaneously  for all $f_i$.

 We can also interpret 2-morphisms of the 2-category of pretopological
 stacks in this way. We will describe briefly how it works, but
 will not be using it later.

 Let $F,G \: [R' \sst{} X']\to[R \sst{}  X]$ be maps of
 groupoids and $f,g \: [X'/R'] \to [ X/R]$ the induced maps
 on the quotient stacks. To describe the 2-morphisms between
 $f$ and $g$ we make use of the ``morphism groupoid'' $[P \sst{} R]$
 of $[R \sst{}  X]$, where $P$ is defined by the following fiber product:

  $$\xymatrix@=12pt@M=10pt{
          P  \ar[r] \ar[d] &  R \ar[d]^{(s,t)}  \\
            R\x R  \ar[r]_(0.45){(s,s)}        &    X\x  X        }$$
 The source map of the  groupoid $[P \sst{} R]$ is the top horizontal
 map of the above square. The target map is given
 by the composition
   $$\xymatrix@C=2pt@M=10pt{P=(R\x R)\x_{ X\x  X}R
            \ar[rrr]^(0.55){(\iota,id,id)} & & &
              R\x_X R\x_X R \ar[rr]^(0.6){\text{mult.}} && R} $$
  \noindent where $\iota$ is the involution and the last arrow is
  the multiplication of three composable arrows. We leave it to the
  reader to figure out the rest of the structure on $[P\sst{}R]$.
 ({\em Remark}. This is a special case of the construction of Section
  \ref{S:fiber}:  take all three groupoids to be $[R\sst{} X]$.)

 There are natural maps of groupoids $S, T \: [P\sst{}R] \to [R\sst{} X]$.
 It can be shown that, the set
 of 2-morphisms $f \twomor{} g$ between $f,g \: [X'/R'] \to [X/R]$ is in
 natural bijection with  the set
 of groupoid morphisms $H \: [R' \sst{} X'] \to [P\sst{}R]$ such that
 $S\circ H=F$ and $T\circ H=G$.

 \vspace{0.1in}

 One can formulate the above discussion as an equivalence between the 2-category
 of pretopological stacks and the 2-category of groupoids with elementary
 Morita equivalences ``inverted''. We will not need it here,
 so we will leave it as a challenge for the reader!

\section{Fiber products}{\label{S:fiber}}

 Let $\C$ be a Grothendieck site. As we saw in Section
 \ref{SS:fiber}, the 2-category of categories fibered in
 groupoids over $\C$ has 2-fiber products (sometimes called
 {\em homotopy fiber products}). The 2-subcategory  $\mathsf{St}_{\C}$
 of stacks is closed  under this 2-fiber product. This is due to the fact that
 stackification  commutes with 2-fiber products.

 Now, let $\C=\mathsf{Top}$ be the site topological spaces. We saw in Section
 \ref{S:pretopological} that
 the 2-subgategory  $\mathsf{PretopSt}_{\mathsf{Top}}$ of pretopological
 stacks is also closed
 under 2-fiber products (Corollary \ref{C:fiber}).
 Here we present another proof that can also be used later
 on to show that the 2-subcategories of topological stacks, and Deligne-Mumford
 topological stack (as defined in the subsequent sections) are also closed
 under 2-fiber products.

 Assume we are given a diagram of pretopological stacks
    $$\xymatrix@=12pt@M=10pt{    &  \ar[d]^f \Y  \\
              \Z  \ar[r]_g        &     \X      }$$
 As we saw in Section \ref{S:Morita}, we can choose groupoid presentations
 for $\X$, $\Y$ and $\Z$, and also for the maps between them.
 So we may assume that the above diagram comes from the following
 diagram of topological groupoids:
     $$\xymatrix@=12pt@M=10pt{      &  [S \sst{} Y] \ar[d]^{F}  \\
                      [T \sst{} Z] \ar[r]_{G}   & [R \sst{} X]   }$$
 We define the {\em fiber product}\footnote{Despite its messy
 appearance, the construction is quite simple. It is simply the construction
 of homotopy fiber product of (discrete) groupoids translated in the
 categorical language.}
  of this diagram to
 be the topological
 groupoid whose base (i.e. the``object space") is
         {\small   $$P_0:=(Y\x Z) \x_{X\x X} R \ \ \,(\text{using}
                \xymatrix@=12pt@M=8pt@C=28pt{
                      Y\x Z \ar[r]^(0.47){(F_0, G_0)}  & X\x X.)}$$ }
 and whose``arrow space" is
        {\small  $$P_1:=(S\x T)\x_{X\x X} R \ \ \,(\text{using}
            \xymatrix@=12pt@M=8pt@C=44pt{
                S\x T \ar[r]^(0.47){(F_0\circ s, G_0\circ s)}  & X\x X).}$$ }
        The source map of the groupoid $[P_1 \sst{} P_0]$
        is just the base extension
        of   the map
        $$\xymatrix@=12pt@M=10pt@C=20pt{ S\x T \ar[r]^(0.47){(s, s)}  & Y\x Z,}$$   so it
        will be nicely behaved (e.g. smooth, \'{e}tale, LF etc.) provided
        those of $[S \sst{} Y]$ and $[T \sst{} Z]$ are.
        The involution $\iota_P \: P_1 \to P_1$ is defined by
{\small $$\xymatrix @M=10pt@C=36pt{ (S\x T)\x_{X\x X} R
\ar[r]^(0.49){(\iota_S,\iota_T,\varphi)}
                                               &(S\x T)\x_{X\x X} R }$$}
        where $\iota$ stands for the involution, and
        the map $\varphi$ is defined by the composition
{\small $$\xymatrix@M=10pt@C=4pt{ (S\x T)\x_{X\x X} R
\ar[rrr]^(0.51){(F_1\circ \iota_S, G_1,\id)} &&&
                         R\x_{X} R\x_{X}R \ar[rr]^(0.63){\text{mult.}}   & &  R  }$$ }
        \noindent Having defined
        $\iota_P$, we define the target map $t \: P_1 \to P_0$
        to be the composition $s\circ \iota_P$.
        We leave it to the reader to figure out the multiplication
        map.

        It is immediate from the construction that the following diagram
        of presheaves of groupoids is 2-cartesian
        (see Section \ref{SS:GpdvsStack} for notation):
             $$\xymatrix@=12pt@M=10pt{  \lfloor P_0/P_1 \rfloor \ar[r] \ar[d]
                                     &  \lfloor Y/S\rfloor  \ar[d]  \\
                           \lfloor Z/T\rfloor \ar[r]
                                         &  \lfloor X/R\rfloor         }$$
        Since stackification commutes with 2-fiber products, the above
        diagram remains 2-cartesian after stackification. Therefore,
        $\PP:=[P_0/P_1]$ is naturally equivalent to the 2-fiber
        $\Z\x_{\X}\Y$.

  \begin{rem}
   Of course, there is nothing specific about $\mathsf{Top}$;
    the same construction is valid for algebraic or analytic stacks.
 \end{rem}

\section{Inertia groups and residue gerbes of pretopological stacks}
{\label{S:inertia}}

 When $\X$ is a pretopological stack, for every point $x \in \X$
 the inertia group  is naturally a topological group. When viewed
 as a topological  group we denote it by $\mathbb{I}_x$,
 and when viewed as discrete group we denote it  by $I_x$.

  The residue stack $\Gamma_x$ at a
  point $x$ also has a nice description, as seen in the Proposition
  \ref{P:residue1} below. In particular, it is a gerbe.
  Whenever we deal with $*$-valued points, we switch the terminology
  from residue stacks to {\bf residue gerbe}s.
  The residue gerbe $\Gamma_x$ of a point $x$
  is of  very local nature, in the sense that, for any substack
  (not necessarily open or closed) $\Y$ of  $\X$   which
  contains $x$ (i.e., $x \: * \to \X$ factors through $\Y$),
  the residue gerbe of $x$ as a point  in $\Y$ is the same
  as the residue gerbe of $x$ viewed as a point in $\X$. In fact, the
  same is  true for the residue stack of any $W$-valued point.

 \begin{prop}{\label{P:residue1}}
    Let $\X$ be a pretopological stack and let $x$ be a  point in $\X$.
      Then we have
    a natural equivalence $\Gamma_x\cong \B \mathbb{I}_x$.
 \end{prop}

  \begin{proof}
        The result  follows from Corollary \ref{C:residue} with $W=*$.
  \end{proof}

 Consider the moduli map $\pim \: \X \to \Xm$. The residue gerbe $\Gamma_x$
 maps to the point $\bar{x}:=\pim(x) \in \Xm$.
 It may appear that $\Gamma_x$ is the fiber of $\pim$ over $\bar{x}$.
 It turns out, however, that this is not always the case. The
 next proposition tells us exactly when this is the case.

 \begin{prop}{\label{P:residue2}}
     Let $\X$ be a  pretopological stack, and let $p \: X \to \X$
     be a chart for it.  Denote the corresponding groupoid
     by  $[R \sst{} X]$. Let  $x$ be a   point in $\X$. Pick a lift
       $x' \in X$   of $x$ (so $p(x')$ is equivalent to $x$),
       and let $X'=\OO(x')$ be
      its orbit (with the subspace topology). Denote the restriction of
      $[R \sst{} X]$ to $X'$ by $[R' \sst{} X']$. Finally, let
       $\Gamma'_x=\im(X')$ (via the natural map $X' \to \X$).
   \begin{itemize}

     \item[$\mathbf{i.}$] There are 2-cartesian diagrams:

        $$\xymatrix@=12pt@M=10pt{ \Gamma_x \ar[r] \ar[d] &  \X \ar[d]
                                    && \Gamma'_x \ar[r] \ar[d] &  \X \ar[d]\\
                    {*} \ar[r]_(0.38){\bar{x}}        &     \pi_0\X
                           &&  {*}\ar[r]_(0.37){\bar{x}}        &     \Xm  }$$
     \item[$\mathbf{ii.}$]  There is a  natural equivalence
         $[X'/R']\risom \Gamma'_x $.

     \item[$\mathbf{iii.}$]  There is  a natural equivalence
         $\B \mathbb{I}_x \risom \Gamma_x$.

     \item[$\mathbf{iv.}$] $\Gamma_x$ is a substack of $\Gamma'_x$. Indeed,
        $\Gamma'_x$  is the smallest embedded substack of $\X$ containing
         $\Gamma_x$ (i.e.  $\Gamma'_x$ is the embedded image of
         $x \: * \to \X$).

     \item[$\mathbf{v.}$]  $\Gamma_x=\Gamma'_x$
     if and only if $t \: s^{-1}(x') \to X'$ is an epimorphism.
     In this case, we have a 2-cartesian diagram
       $$\xymatrix@=12pt@M=10pt{
           \B \mathbb{I}_x \ar[r] \ar[d] &  \ar[d] \X \\
             {*} \ar[r]        &  \Xm         }$$
   \end{itemize}
 \end{prop}

  \begin{proof}[Proof of part {\em ($\mathbf{i}$)}]
      The left square is 2-cartesian
     by Corollary \ref{C:cartesian}. To prove that the right square
     is 2-cartesian, let $\Gamma$ be the substack   of $\X$
     that is equivalent to the 2-fiber product $*\x_{\Xm}\X$; we
     have to show that $\Gamma=\Gamma_x'$.
     There is  a natural map
     $X' \to \Gamma$ which makes $\Gamma'_x$ a
     substack of $\Gamma$. It is enough to show that this map is
     an epimorphism. This follows from the fact that the following diagram is
     2-cartesian:
            $$\xymatrix@=12pt@M=10pt{
           X' \ar[r] \ar[d] &  \ar[d] X \\
             \Gamma \ar[r]        &    \X       }$$

\vspace{0.1in}
      \noindent{\em Proof of part} ($\mathbf{ii}$). This follows from Proposition
       \ref{P:image}.

\vspace{0.1in}
     \noindent{\em Proof of part} ($\mathbf{iii}$). This is Proposition
     \ref{P:residue1}.

\vspace{0.1in}
     \noindent{\em Proof of part} ($\mathbf{iv}$).
       Since $x \: * \to \X$ factors through $\Gamma'_x$, the inclusion
       $\Gamma_x \subseteq \Gamma'_x$ is obvious. The second statement follows
       from Corollary \ref{C:embedded}.

\vspace{0.1in}
     \noindent{\em Proof of part} ($\mathbf{v}$).
       The equality $\Gamma_x=\Gamma'_x$
       holds if and only if $x \: * \to \Gamma'_x$ is an
       epimorphism. Using the cartesian diagram
           $$\xymatrix@=12pt@M=10pt{
               s^{-1}(x') \ar[r]^(0.56)t \ar[d] &  \ar[d] X' \ar[r] & X \ar[d] \\
                {*} \ar[r]        &     \Gamma'_x   \ar[r] & \X  }$$
       (in which the vertical maps are epimorphisms), the bottom left
       map being an epimorphism is equivalent to the top left map being
       an epimorphism.

       The last statement follows from ($\mathbf{i}$) and ($\mathbf{iii}$).
 \end{proof}

 The point of the above proposition is that, given a {\em transitive} groupoid
 $[R \sst{} X]$, the quotient stack $[X/R]$ is not necessarily equivalent
 to $\B \mathbb{I}$, where $\mathbb{I}$
 is the stabilizer group (viewed as a topological group)
 of a point on $X$.
 The former is typically larger than the latter. The equality holds,
 if and only if the assumption ($\mathbf{v}$) of the proposition is satisfied.
 When ($\mathbf{v}$) is satisfied at every point on a pretopological stack
 $\X$, we can think of $\X$ as  a ``family of classifying stacks
 $\B \mathbb{I}_x$, parameterized by $\Xm$".

   \begin{ex}{\label{E:discreteorbits}}
      Let $[R \sst{} X]$ be a groupoid and assume $x \in X$ is a point
      for which the orbit $\OO(x)$
      is discrete. Then the condition ($\mathbf{v}$) of the proposition
      of the proposition is satisfied at $x$. So we have a 2-cartesian diagram
       $$\xymatrix@=12pt@M=10pt{
            \B \mathbb{I}_x \ar[r] \ar[d] &  \ar[d] \X \\
             {*} \ar[r]        &  \Xm         }$$
   \end{ex}

 There are two major classes of stacks for which the condition
  ($\mathbf{v}$) of Proposition \ref{P:residue2} is automatically
  satisfied at every point: Deligne-Mumford topological stacks
  (Section \ref{S:DM})
  and pretopological gerbes (use
  Corollary \ref{C:gerbe1} and Proposition \ref{P:residue2}.$\mathbf{i}$).

    \begin{ex}{\label{E:Q}}
      Let $X=\mathbb{Q}$ with the subspace topology induced from
      $\mathbb{R}$. Let $\mathbb{Q}$ (viewed as a discrete group)
      act on $X$ by translations, and let $\X$
      be the quotient stack. Since this action
      is transitive, we have $\Xm=*$. So, $\X$ is a quasitopological space
      that is not a topological space.
    \end{ex}

    \begin{ex}{\label{E:Z2}}
      Another counter-example is obtained by taking $X$ to be the trivial
      topology on a set with 2 elements and letting $\mathbf{Z}_2$ act
      on it by swapping the elements. The coarse moduli space of
      the quotient stack is obviously $*$.
      We have a natural representable map
      $* \to  [X/\mathbf{Z}_2]$ that identifies $*$ with a substack of
      $[X/\mathbf{Z}_2]$, but this substack is not embedded.
     \end{ex}


\section{The coarse moduli space in the presence of
                open charts}{\label{S:coarse2}}

  In this section we assume all out pretopological stacks admit open
  charts.

 In the presence of an {\em open} chart for a pretopological stack $\X$,
 the coarse moduli space is better-behaved. We will address this issue
 in this section. Recall that, by an open chart $p \: X \to \X$ for
 $\X$ we mean one for which $p$ is an open map (Definition
 \ref{D:representable2} and Example \ref{E:P}).

 We can extend the notion of an open map to
  morphisms of pretopological stacks that are not necessarily representable.

  \begin{defn}{\label{D:open}}
  We say that a (not necessarily representable) map $f \: \Y \to \X$ of
  (not necessarily pretopological)
  stacks is {\bf  open} if there exists diagram of stacks
      $$\xymatrix@=12pt@M=10pt{
             \Z \ar[rd]^g \ar[d]_h &     \\
             \Y \ar[r]_(0.44)f        &  \X        }$$
  where $h$ is a  surjective map and $g$ is a representable
  open map. (For a not necessarily representable map
     $h$ surjective means $h_{mod}$ is surjective.)
 \end{defn}

   It is  easy to show that base extension of an open
 map is open, and so is the composition of two open maps.
 It is also easy to show that when $f$ is representable this
 definition coincides with the previous definition.

 \begin{lem}{\label{L:open}}
    Consider a diagram of stacks (not necessarily pretopological)
     $$\xymatrix@=12pt@M=10pt{
             \Z \ar[rd]^g \ar[d]_h &     \\
             \Y \ar[r]_(0.44)f        &  \X        }$$
     in which $h$ is surjective
     and $g$ is open. Then $f$ is open.
    \end{lem}

   \begin{proof}
     Since $g$ is open,  there is a stack $\Z'$ and  maps
     $h' \: \Z' \to \Z$ and $g' \: \Z' \to \X$ such that $h'$ is surjective
     and $g'$ is representable and open. The following diagram
     does the job:
     $$\xymatrix@=12pt@M=10pt{
             \Z' \ar[rd]^{g'} \ar[d]_{h\circ h'} &     \\
             \Y \ar[r]_(0.44)f        &  \X        }$$
    \end{proof}

 \vspace{0.1in}
  \noindent{\bf Exercise}.
  Show that if $f\: \Y \to \X$ is open, then so is $f_{mod} \: \Ym \to \Xm$,
  but the converse is not necessarily true.

  \begin{lem}{\label{L:naivequotient}}
  Let $\X$ be a pretopological stack, and
     $[R \sst{} X]$   a groupoid presentation for it.
  \begin{itemize}

   \item[$\mathbf{i.}$]
      There is a  natural homeomorphism $X/R \risom \Xm$,
      where $X/R$ is the (naive) quotient space of the equivalence
      relation induced on $X$ from $R$.

   \item[$\mathbf{ii.}$] If the source  map of the groupoid
     is open, then so is $X \to \Xm=X/R$.

  \end{itemize}
   In particular, if $\X$ is a pretopological stack that admits an
   open chart, then the moduli map $\pim \: \X \to \Xm$ is open.
  \end{lem}

  \begin{proof}
   Easy.
  \end{proof}

 If for a stack $\X$ the moduli map $\pim \: \X \to \Xm$ is open, then
 it is uniquely characterized by the following three properties:  open,
  continuous, bijection. By the above proposition, this is the case when
 $\X$ is the quotient stack of a topological groupoid whose source (hence also
 target) map is open.

  \begin{cor}[{\bf Invariance under base change}]{\label{C:basechange}}
  Let $\X$ be a stack such that the moduli map $\pim \: \X \to \Xm$ is open
  (e.g.  a pretopological stack that admits an open chart),
  and let $f \: Y \to \Xm$ be a
  continuous map. Set $\Y=Y\x_{\Xm}\X$. Then, the projection
  map $\pi \: \Y \to Y$ makes $Y$ into the coarse moduli space of $\Y$.
 \end{cor}

 \begin{proof}
   It is easily checked that $\pi \: \Y \to Y$ induces a natural continuous
   bijection $\Ym \to Y$. The map  $\pi \: \Y \to Y$ is open,
   being  base extension of an   open map. This implies that
   $\Ym \to Y$ is also open (Lemma \ref{L:open}), hence a homeomorphism.
 \end{proof}

\section{Quotient stacks as classifying spaces for torsors}{\label{S:quotient}}

  Let $[R \sst{} X]$ be a groupoid in $\mathsf{Top}$. There is an alternative
  description of the quotient stack $[X/R]$ in terms of torsors which
  is both technically and conceptually very important. Roughly speaking,
  the quotient stack $[X/R]$ {\em parameterizes}  torsors for the
  groupoid $[R \sst{} X]$ (see Definition \ref{D:torsor}). For example,
  if $[G \sst{} *]$ is the action groupoid of a topological group $G$
  acting (trivially) on a point, then maps from a topological space $Y$ into
  $[*/G]$ correspond to $G$-torsors over $Y$. In this sense,
  $\B G=[*/G]$ can be thought of as the {\em classifying space} of $G$.

     \begin{defn}{\label{D:cartesian}}
        A map  $[R \sst{} X] \to [R' \sst{} X']$ of
     groupoids  is called {\bf cartesian} if
     the following square is cartesian
       $$\xymatrix@=12pt@M=10pt{ R  \ar[r] \ar[d]_s & R' \ar[d]^s  \\
                  X  \ar[r]        &     X'      }$$
    Note that the same will be true for $t$.
    \end{defn}

     \begin{rem}{\label{R:cartesian}}
        The above definition is equivalent to saying that
        the following diagram of presheaves of groupoids
        is 2-cartesian:
            $$\xymatrix@=12pt@M=10pt{ X \ar[r] \ar[d] &  X'\ar[d]  \\
              \ar[r]   \lfloor X/R \rfloor     &    \lfloor X'/R' \rfloor  }$$
        Since stackification commutes with fiber products, this diagram
        remains 2-cartesian after stackification.
     \end{rem}

  As far as I know, the following definition is due to Kai Behrend.

   \begin{defn}{\label{D:torsor}}
     Let $[R \sst{} X]$ be a groupoid in $\mathsf{Top}$, and let $W$ be a
     topological
     space. By an $R$-{\bf torsor} (more precisely, a $[R \sst{} X]$-torsor)
     over $W$ we mean an epimorphism
     $T \to W$ of topological spaces,
     together with a cartesian map of groupoids
     $$[T\x_W T \sst{} T] \to [R \sst{} X].$$
   \end{defn}

   Given $R$-torsors $T\to W$ and $T'\to W'$, we define a {\em morphism}
   of $R$-torsors from $T$ to $T'$ to be a cartesian square
       $$\xymatrix@=12pt@M=10pt{T  \ar[r] \ar[d] & T' \ar[d]  \\
                                 W  \ar[r]        &   W'        }$$
    such that the induced diagram of groupoids commutes:

    $$\xymatrix@R=10pt@M=10pt@C=-18pt{ [T\x_W T \sst{} T]  \ar[rr]  \ar[rd] &  &
                                             [T'\x_{W'} T' \sst{} T'] \ar[dl] \\
                                       & [R \sst{} X]  &           }$$
    The category of all $R$-torsors, over various $W$,  is fibered in groupoids
    over $\mathsf{Top}$ via the forgetful functor
    (that out of everything only remembers $W$).

    %
  %
%

     %


    \begin{ex}{\label{E:torsor}}
    \end{ex}
  \begin{itemize}

          \item[$\mathbf{1.}$] Take $[R \sst{} X]$ be the action
          groupoid $[G\x X \sst{} X]$ of a topological group $G$ acting on a
          topological
          space $X$. An $R$-torsor over $W$ consists of a pair $(T,\al)$,
          where $T \to W$ is a $G$-torsor (that is, a locally trivial principal
          $G$-space over $W$), and $\al \: T \to X$ is a $G$-equivariant map.
          For a fixed $W$, a morphism between  $R$-torsors $(T,\al)$ and
          $(T',\al')$ over $W$ is a $G$-equivariant map $a \: T \to T'$
          relative to $W$ such that
          the following triangle commutes:
           $$\xymatrix@M=8pt@R=10pt@C=-4pt{
              T  \ar[rr]^(0.47)a  \ar[rd]_(0.44){\al'} &
                                                      &  T' \ar[dl]^{\al} \\
                                       & X  &           }$$
          \item[$\mathbf{2.}$] In the previous example, assume the action of
          $G$ on $X$ is trivial. Then, an $R$-torsor over $W$ is
          of a pair $(T,\be)$ consisting of a $G$-torsor $T$ over $W$
          and continuous map $\be \: W \to X$. When $X=*$ an $R$-torsor
          is simply a $G$-torsor.

          \item[$\mathbf{3.}$] We could also consider the previous examples
          with everything being relative to a base $B$. So, $X$ is now a space
          over $B$, and $G$ is a topological group over $B$, acting
          on $X$ (relative to $B$). A particular case of interest is when $X=B$.
        \end{itemize}
\vspace{0.1in}

        We say an $R$-torsor $T \to W$ is {\bf trivial} if it admits
        a section $\sigma \: W \to T$. Having fixed such a section, the set
        of $R$-torsor morphisms $a \: T \to T'$ to an $R$-torsor $T' \to W$
        is in natural bijection with the set of sections
        $\sigma' \: W \to T'$ for which the following square commutes:
             $$\xymatrix@=12pt@M=10pt{ T' \ar[r]  &   X \\
             W \ar[r]_{\sigma}    \ar@{.>}[u]^{\sigma'}    &       T \ar[u] }$$

        The bijection is given by
        $a \mapsto a(\sigma)$.

          \begin{rem}{\label{R:trivial}}
             An arbitrary $R$-torsor $T$ over $W$ is locally trivial.
             That means, there is an open covering $\{U_i\}$ of $W$
             such that the restriction of $T$ to each $U_i$ is trivial.
             This is because $T  \to W$ is an epimorphism.
          \end{rem}

         \vspace{0.1in}

        The following proposition explains why  $R$-torsors are
        related to stacks.

         \begin{thm}{\label{T:torsor}}
           Let $[R \sst{} X]$ be a groupoid in $\mathsf{Top}$.
           We have a natural equivalence (of categories fibered in groupoids
           over $\mathbf{Top}$)
             $$ [X/R]\risom \Tors_R.$$
           In particular, $\Tors_R$ is a pretopological stack over
           $\mathsf{Top}$.
         \end{thm}

         \begin{proof}
             Denote $[X/R]$ by $\X$.
             We define a functor $\Theta \:  \X \to \Tors_R$
             as follows. Let $W \in \mathsf{Top}$ and pick on object
             $f \in \X(W)$.  We denote the corresponding map $W \to \X$
             also by $f$. We define $\Theta(f)$ to be the $R$-torsor
             $T \to W$, where $T=W\x_{\X} X$ (the map from $W$ to $\X$ being $f$).
             The fact that
             $T \to W$ has an $R$-torsor structure follows from formal
             properties of 2-fiber products. Similarly, given a morphism
             $\al \: f \twomor{} f'$ in $\X(W)$, we obtain an induced map
             of fiber products  $a \: T \to T'$ (where $T'=W\x_{\X} X$, the map
             from
             $W$ to $\X$ being $f'$). Again, it follows from formal properties
             of 2-fiber products that this is indeed a map of $R$-torsors.
             This defines the functor $\Theta \:  \X \to \Tors_R$.

             We now show that $\Theta$ is an equivalence. We have to show
             that $\Theta$ is fully faithful and essentially surjective.

\vspace{0.1in}
             \noindent {\em Proof of faithfulness.}
                Let $f,f' \in \X(W)$ be $W$-points (we use the same
                notation for the corresponding maps $W \to \X$) and
                $\al,\be \: f \twomor{} f'$ morphisms between them in $\X(W)$
                (which we view as 2-isomorphisms between the
                corresponding
                maps $W \to \X$).
                Let $T,T' \to W$ be the corresponding $R$-torsors,
                with $a, b \: T \to T'$ the $R$-torsors morphisms
                associated to $\al$ and $\be$. We have to show that $a=b$
                implies $\al=\be$. Since $\X$ is a stack, it is enough to check
                the equality $\al=\be$ locally, so we may replace $W$ by an open
                cover and assume that $f$   lifts to $X$. This is
                equivalent to saying that the torsors $T$ is trivial.
                We fix a choice of a lift for $f$   (by this
                we mean a pair $(\tilde{f},\gamma)$,
                where $\tilde{f} \: W \to X$
                is a map of topological spaces and
                $\gamma \: f  \twomor{} p\circ\tilde{f}$ is an identification).
                By definition of 2-fiber product, there is a natural bijection
                between such lifts and sections of $T \to W$. So,
                $(\tilde{f},\gamma)$ corresponds to a section
                $\sigma \: W \to T$.

                  The sections $a(\sigma), b(\sigma) \: W \to T'$
                  correspond to the lifts $(\tilde{f}, \al^{-1}\gamma)$
                  and $(\tilde{f}, \be^{-1}\gamma)$ of $f'$.
                  If $a=b$, then these two lifts are equal, so
                  $\al^{-1}\gamma=\be^{-1}\gamma$. This implies $\al=\be$.

\vspace{0.1in}
                In fact, the above argument proves the fullness as well,
                provided we assume the existence of lifts.

\vspace{0.1in}
             \noindent {\em Proof of  fullness.}
                Let $W$, $f$,$f'$, $T$ and $T'$ be as in the previous part.
                Let $a \: T \to T'$ be an $R$-torsor morphism.
                We have to show that it comes from an $\al \: f \twomor{} f'$.
                We saw in the previous part that the statement is true
                 if $T$ is trivial (which
                implies that $T'$ is also trivial).
                In the general case,  choose a covering $\{U_i\}$ of $W$
                over which $T$ becomes trivial (Remark \ref{R:trivial}).
                Then, we obtain a unique family
                $\al_i \: f|_{U_i} \twomor{} f'_{U_i}$ inducing
                $a|_{U_i} \: T|_{U_i} \to T'|_{U_i}$. Since the $a|_{U_i}$ are
                compatible over the double intersections, so will be the $\al_i$
                (we are using faithfulness). So $\al_i$ glue to a global
                $\al \: f \twomor{} f'$
                over $W$. Let $b \: T \to T'$ be the $R$-torsor map
                induced
                by $\al$. Then  $a$ and $b$ are equal over every open $U_i$.
                It is easy to check that this implies that $a=b$. In other
                words, $a$ is induced from $\al$.

\vspace{0.1in}
             \noindent {\em Proof of  essential surjectivity.}
                Let $T \to W$ be an $R$-torsor over $W$. By definition,
                we have a cartesian diagram
                  $$\xymatrix@=12pt@M=10pt{
                     T\x_W T \ar[r] \ar[d]_{pr_1} &  \ar[d] R  \\
                     T      \ar[r]        &    X       }$$
                By Remark \ref{R:cartesian}, and using the natural
                equivalence     $[T/T\x_W T ] \risom W$, we obtain
                  a 2-cartesian diagram
                   $$\xymatrix@=12pt@M=10pt{
                              T \ar[r] \ar[d] &  \ar[d] X  \\
                              W \ar[r]        &        [X/R]   }$$
                Denoting the bottom map by $f \: W \to [X/R]$,
                we see that the  $R$-torsor $T \to W$ is isomorphic
                to $\Theta(f)$. This proves the essential surjectivity.
         \end{proof}

         The upshot is that the quotient
         stack $[X/R]$ can be thought of as ``the moduli stack
         of $R$-torsors'', with $X \to [X/R]$ being the ``universal'' $R$-torsor
         over it. In other words, an $R$-torsor $T$ over $W$ gives rise to a
         map $f \: W \to [X/R]$ and $T$ is just the pull back of
         of the universal $R$-torsor along $f$. Given two
         $R$-torsors $T$ and $T'$ over $R$, a morphism between them
         corresponds to a 2-isomorphism between the corresponding
         classifying maps $f, f' \: W \to [X/R]$.

         The reader is advised to work out the meaning of the above proposition
         for the special cases considered in Example \ref{E:torsor}.


\newpage
{\LARGE\part{Topological stacks}}{\label{PART2}} \vspace{0.5in}
\section{Topological stacks}{\label{S:TopSt}}

 In this section, we define our main objects of interest, the
 {\em topological stacks}. A topological stack is a pretopological
 stack which admits  a chart that is a {\em local fibration}
 (LF, for short),
 where a local fibration is defined to be a member
 of  a class $\mathbf{LF}$ of continuous maps satisfying
 certain axioms (see Subsection \ref{SS:fibrations} below).

 Therefore, our notion of topological stack is not an absolute
 one and depends on the choice of the class $\mathbf{LF}$. The freedom
 in choosing our own notion of local fibration gives us flexibility
 in crafting the corresponding notion of topological stack meeting
 our particular needs. For example, if we need to do homotopy theory
 with our stacks, {\em local Serre fibrations} (see Example \ref{E:lf})
 is a good choice. Or, if we want the coarse moduli spaces to be well-behaved
 we require that all the maps in $\mathbf{LF}$ to be open.

 We will see that the desired properties of a category of
 topological stacks are often reflected as some axioms imposed on
  $\mathbf{LF}$.

 Topological stacks are the topological
 counterparts of algebraic stacks. In fact, we will see in
 Section \ref{S:algebraic}
 how to associate a topological stack to an algebraic stack
 (of finite type) over $\mathbb{C}$, much the same
 way that one associates a topological space to a scheme (of finite type)
 over $\mathbb{C}$.  In  Section \ref{S:DM}, we consider (weak)
 Deligne-Mumford topological stack, which are the topological counterparts
 of  algebraic Deligne-Mumford stacks.

This section has two subsections. In the first one we introduce
local fibrations and supply some examples. In the second
subsection, we use the notion of local fibration to define
topological stacks.

\subsection{Local fibrations}{\label{SS:fibrations}}

For a pretopological stack to behave nicely, we need to require
some kind of `fibrancy' condition on the chart $p \: X \to \X$.
The stronger the  fibrancy condition is, the more manageable our
stack becomes. What this exactly means is discussed in more
details in Section \ref{S:gluing}. Let us just mention that, the
fibrancy conditions on charts will result in having more freedom
in performing push-outs in  the corresponding category of stacks.

 By a class $\mathbf{LF}$ of {\bf local fibrations} we mean
 a collection of continuous maps satisfying the following
 axioms:

  \begin{itemize}

          \item[$\mathbf{LF1.}$] Every open embedding is a local fibration;

          \item[$\mathbf{LF2.}$] Local fibrations
          are closed under composition;

          \item[$\mathbf{LF3.}$] Being a local fibration is stable under base
             change and local on the target (Definition \ref{D:local});

             \item[$\mathbf{LF4.}$] If $f_i \: X_i \to Y$ is a family
               of local fibrations, then $\coprod f_i \: \coprod X_i \to Y$
               is also a local fibration.

        \end{itemize}

 \noindent Having fixed $\mathbf{LF}$, we refer to a map in $\mathbf{LF}$
 by an {\em LF map}.

 Note that, by virtue of $\mathbf{LF3}$, we can talk about
 LF maps of stacks. More precisely, we say a representable map $f \: \X \to \Y$
 of stacks is LF, if there is a chart $Y \to \Y$ such the base extension
 of $f$ over $Y$ is LF.

  \begin{ex}{\label{E:lf}}
     Let $\mathbf{Fib}$ be any of the following classes
     of maps:
      \begin{itemize}
       \item[$\mathbf{1.}$] {\em Serre fibrations}. These are maps that
       have the
          homotopy lifting property for finite
          CW complexes.


          \item[$\mathbf{2.}$] {\em Hurewic fibrations}. These are maps
          that have the homotopy lifting property for all topological spaces.


          \item[$\mathbf{3.}$] {\em Cartesian maps}. These are projection maps
           of products, i.e $pr_2 \: X\x B \to B$.

          \item[$\mathbf{3'.}$]
          There are variations of the notion of cartesian map in
          which the fibers are assumed to belong to a certain family of spaces.
          For instance, we can require the fiber to be  a Euclidean
          space (i.e. homeomorphic to some $\bbR^n$).

          \item[$\mathbf{4.}$] {\em Right invertible maps}. These are maps
          that admit sections.

          \item[$\mathbf{5.}$] {\em Open maps}.

          \item[$\mathbf{6.}$] {\em Homeomorphisms}.

        \end{itemize}

  Define the corresponding class $\mathbf{LF}$
  of local fibrations by declaring
   a continuous map $f \: X \to Y$ of topological spaces
   to be LF,
   if for every point $x \in X$ there are
   open sets $U\subseteq X$ containing $x$
   and  $V\subseteq Y$ containing $f(x)$ such that
   $f|_V \: U \to V$ is in $\mathbf{Fib}$.  It is easy to check that
   the resulting class of maps satisfies the axioms for local fibrations.
\end{ex}

     The local fibrations for Example ($\mathbf{3}$) are sometimes called
    {\em topological submersion}. Every submersion
    of differentiable manifolds is a local fibration for Example ($\mathbf{3'}$)
    (Implicit Function Theorem).
    Local fibrations of Example ($\mathbf{4}$)  are the epimorphisms.
    Local fibrations of Example ($\mathbf{5}$) are the open maps.
     Local fibrations of Example ($\mathbf{6}$) are the local homeomorphisms.

   \begin{ex}{\label{E:morelf}}
       Intersection of any two classes of local fibrations is again a class
         of local fibrations. So we could from new classes such as
         {\em open epimorphic local Serre fibrations}  and so on.
   \end{ex}

 Every notion of fibration comes with the corresponding notion
 of {\bf cofibration} (rather, {\em trivial} cofibration),
 which is defined using certain lifting properties. So, with
 {\em local} fibrations we have the corresponding notion
 of {\em local} trivial cofibrations.

  \begin{defn}{\label{D:LTC}}

   Let $f \: X \to Y$ be a map of topological spaces. We say
   a map $i \: A \to B$ of topological spaces
   has {\em local left lifting property}
   with respect to $f$ if in every commutative diagram
      $$\xymatrix@=12pt@M=10pt{
            A \ar[r]\ar[d]_i &    \ar[d]^f X \\
                B  \ar[r] \ar@{..>}[ru]       &    Y     }$$
  and for every point $a \in A$, there exist small enough neighborhoods
  $U \subseteq B$ of $i(a)$ such that, after replacing $i$ by
  $i|_U \: i^{-1}(U) \to U$, the dotted arrow can be filled.
   For a fixed  class $\mathbf{LF}$ of local fibrations, we say
   a map $i \: A \to B$
    is  a
   {\bf local (trivial) cofibration}
   (LTC for short),  if it satisfies the local left lifting property
   with respect to every $f \in \mathbf{LF}$.
 \end{defn}

   \begin{lem}{\label{L:LTC}}
         Fix a class $\mathbf{LF}$ of local fibrations.
    \begin{itemize}

          \item[$\mathbf{i.}$] Every open embedding is LTC.

          \item[$\mathbf{ii.}$] LTC maps
                                are closed under composition.

           \item[$\mathbf{iii.}$] LTC maps are stable under push out along
            embeddings. That is, if in the following  push out diagram
            $i$ is LTC and $j$ is an embedding, then $i'$ is also LTC:

          $$\xymatrix@=12pt@M=10pt{
                 A \ar@{^(->}[r]^j \ar[d]_i &  A' \ar[d]^{i'}  \\
                  B   \ar@{^(->}[r]        &   B'        }$$

            \item[$\mathbf{iv.}$] Let $i \: A \to B$ be LTC, and let
            $U \subseteq B$ be an open set. Then $i|_U \: i^{-1}(U) \to U$
            is LTC.

          \item[$\mathbf{v.}$] Being  LTC is local on the target, that is, if
             $i \: A \to B$ is a map of topological spaces such that, for
             some collection  $\{U_i\}$ of open subsets of $B$ whose
             union contains $i(A)$, the maps
              $i^{-1}(U_i) \to U_i$ are   LTC, then so is $i$.

          \item[$\mathbf{vi.}$]  Let $i \: A \to B$ be a map
            of topological spaces. Let $p \: B' \to B$ be an epimorphism,
             and let $A'=B'\x_{B}A$. If $i': A' \to B'$ is LTC then so is $i$.
          \item[$\mathbf{vii.}$] Assume $i \: A\coprod A' \to B$ is LTC.
              Then $\overline{i(A)}\cap i(A')=\emptyset$.
        \end{itemize}
    \end{lem}

   \begin{proof} Parts ($\mathbf{i}$)-($\mathbf{v}$)
   are straightforward.

   \vspace{0.1in}
   \noindent {\em Proof of part} ($\mathbf{vi}$).
       Since being LTC is local on the target, we may assume that $p$
       has a section $s \: B \to B'$. Consider a diagram
         $$\xymatrix@=12pt@M=10pt{
            A \ar[r]\ar[d]_i &    \ar[d]^f X \\
                B  \ar[r] \ar@{..>}[ru]^q      &    Y     }$$
                where $X \to Y$ is LF.
       Pre-composing the horizontal arrows with $p$, we obtain
       the following diagram
         $$\xymatrix@=12pt@M=10pt{
            A' \ar[r]\ar[d]_{i'} &    \ar[d]^f X \\
                B'  \ar[r] \ar@{..>}[ru]^{q'}       &    Y     }$$
       Since $i'$ is LTC, there is a    family $\{U'_i\}$ of open
       subsets of $B'$ whose union contains $i(A')$ and over which the
       lift $q'$ exists. Set $U_i:=s^{-1}(U'_i)$. Then $\{U_i\}$
       is a family   of open
       subsets of $B$ whose union contains $i(A)$ and over which the
       lift $q$ exists: take $q=q'\circ s$.

    \vspace{0.1in}
  \noindent  {\em Proof of part} ($\mathbf{vii}$).
      Consider the local fibration $f=\id\coprod\id \: B\coprod B \to B$.
      Let $g=i|_A\coprod i|_{A'}\: A\coprod A' \to B\coprod B$.
      The local left lifting property for the following diagram
      implies the claim:
         $$\xymatrix@=12pt@M=10pt{
             A\coprod A' \ar[r]^{g}
                      \ar[d]_j &    \ar[d]^f  B\coprod B \\
                B  \ar[r]_{\id} \ar@{..>}[ru]       &    B     }$$
    \end{proof}

  \begin{rem}
    By part ($\mathbf{vi}$) of the above proposition, we can define LTC
    for every representable map of stacks.
  \end{rem}

    \begin{rem}{\label{R:LTC}}
      For a given class $\mathbf{C}$
      of maps in $\mathbf{Top}$ satisfying property ($\mathbf{vii}$)
      of Lemma \ref{L:LTC}, let $\mathbf{LF}$ be the class  of all maps
      $f \: X \to Y$ with the property that every  $i \in \mathbf{C}$
      satisfies the local left lifting property with respect to
      $f$. It can be checked that $\mathbf{LF}$ satisfies the axioms
      $\mathbf{LF1}$-$\mathbf{LF4}$. This gives an alternative
      way to produce classes of local fibrations.
    \end{rem}

        Let us look at some examples of LTC maps, corresponding
        to classes of local fibrations considered in Example \ref{E:lf}.

    \begin{ex}{\label{E:LTC}}
          Numbering corresponds to the numbering in  Example \ref{E:lf}.
    \begin{itemize}

          \item[$\mathbf{1.}$] Every CW inclusion of finite CW complexes
             is LTC. To prove this, use the fact that every such inclusion is
             locally homotopy equivalent to the $t=0$ inclusion $X \hra I\x X$,
             where $X$ is a finite CW complex.

          \item[$\mathbf{2.}$] Every $i \: A \to B$ that is a
            Hurewic cofibration (or only locally on $B$ so) is LTC.
            This follows from \cite{Strom}, Lemma 4.

          \item[$\mathbf{3.}$] Every embedding that is a local retract
            is LTC.

          \item[$\mathbf{3'.}$] {\em Locally cartesian maps with Euclidean
            fibers}. Every embedding that
            locally satisfies the Tietze extension property is LTC. For
            instance,
            every inclusion of a locally-closed subspace in a normal topological
            space is LTC.

          \item[$\mathbf{4.}$] The only LTC maps are open embeddings.


          \item[$\mathbf{6.}$] Every embedding is LTC.
        \end{itemize}
    \end{ex}

\subsection{Topological stacks}

Fix your favorite class $\mathbf{LF}$ of local fibrations (see the
previous subsection).

  \begin{defn}{\label{D:topological}}
   Let $\X$ be a  pretopological stack. We say $\X$ is
   {\bf topological} if there is a chart $p \: X \to \X$
   that is LF.
  \end{defn}

For a choice $\mathbf{LF}$ of a class of local fibrations, we have
the corresponding theory of topological stacks. The smaller the
class of local fibrations, the smaller the corresponding category
of topological stacks, the more flexibility  in performing
push-outs in the corresponding 2-category of topological stacks
(as we will see in Section \ref{S:gluing}).

Our favorite choices for $\mathbf{LF}$ are the ones listed in
Example \ref{E:lf}, or their combinations (Example
\ref{E:morelf}). After Example $4$, which gives us the good old
pretopological stacks, the topological stacks of Example $1$ give
us the most general theory. The ones for Example $3'$ appear as
underlying topological stacks of differentiable and also algebraic
stacks.  Topological stack of Example 5 are the ones studied in
Section \ref{S:coarse2}. Topological stacks of Example 6 are what
we call {\em weak Deligne-Mumford stacks}. These are studied in
the next Section.



 Topological stacks form a full 2-subcategory of $\mathsf{St}_{\mathsf{Top}}$,
 which we denote by $\mathsf{TopSt}$. The Yoneda functor identifies
 $\mathsf{Top}$ with a full 2-subcategory of $\mathsf{TopSt}$.

   \begin{prop}{\label{P:basic}}
    Let $\X$ be a pretopological stack.
   \begin{itemize}

          \item[$\mathbf{i.}$] If $\X$ is topological,
          then so is every embedded substack of $\X$.

          \item[$\mathbf{ii.}$] Let $\{\U_i\}_{i\in I}$ be a
          covering of $\X$ by open substacks. If every $\U_i$ is topological,
          then so is $\X$.

          \item[$\mathbf{iii.}$] The 2-category $\mathsf{TopSt}$ is
          closed under  2-fiber products.
        \end{itemize}

   \end{prop}

  \begin{proof} Easy.
  \end{proof}

  There is an inherent difficulty in performing push-outs in the
  2-category of topological stacks. In fact, the notion of
  {\em local fibration} is precisely designed to bring this issue under control.

  An unfortunate fact is that, even simple push-outs of topological {\em spaces}
  may no longer be push-outs  when viewed in $\mathsf{TopSt}$.
  It is  also unreasonable to expect to be able to quotient
  out a topological stack modulo a, say, closed substack. For this
  we {\em need} to impose certain conditions on the closed substack,
  or else the quotient can be shown not to exists in any reasonable
  sense.\footnote{For example, the {\em suspension} construction
  usually fails for topological stacks.}

  Similar problems occur when we try to glue topological stacks along
  substacks. In Section \ref{S:gluing}
   we will discuss the issue of  push-outs in $\mathsf{TopSt}$
   in more detail.

\section{Deligne-Mumford topological stacks}{\label{S:DM}}

 In this section we introduce a class of topological stacks which
 are counterparts of Deligne-Mumford stacks in algebraic geometry.
 We call them {\em Deligne-Mumford topological stacks}. In fact, there is a
 weaker notion, called {\em weak Deligne-Mumford topological stack}, which
 could also be regarded as topological counterpart of algebraic
 Deligne-Mumford stacks. The weak Deligne-Mumford topological stacks, however,
  could behave pathologically in certain situations, since they are a bit
  too general.

 We begin with a definition which is modeled on the notion of a ``slice''
 of a group action. Group actions for which slices exist will satisfy the
 following property.

  \begin{defn}{\label{D:mild}}
   Let $G$ be a topological group acting continuously on a topological space
   $X$ and let $x$ be a fixed point of this action. We say this action
   is {\bf mild} at $x$ if every open neighborhood of $x$ contains
   an invariant open neighborhood (that is, the invariant opens form
   a basis at $x$).
 \end{defn}

For example, any finite group action is mild at every fixed point.
Any continuous action of a compact Lie group on a topological
space also has this property, but in this paper we are only
concerned with the discrete groups, so we will not get into that.

  \begin{defn}{\label{D:pd}}
      Let $G$ be a discrete group acting on a topological space $X$.
      We say that
      this action is {\bf properly discontinuous} if, for every $x \in X$,
      there is an  $I_x$-invariant open neighborhood $U$ of $x$ such that:
       \begin{itemize}

          \item[$\diamond$]   $I_x$ acts mildly on $U$;

          \item[$\diamond$] for every  $g \in G\backslash I_x$,
               $U\cap g(U)=\emptyset$.

        \end{itemize}
    Here $I_x$ stands for the stabilizer group at $x$.
    Note that the mildness condition is automatic if $I_x$ is finite.
 \end{defn}

\begin{defn}{\label{D:DM}}
 A pretopological stack $\X$ is called a {\bf weak Deligne-Mumford
 topological stack} if  there is a chart $p \: X \to \X$ that is
 a local homeomorphism. A weak
 Deligne-Mumford topological stack $\X$ is called a
 {\bf  Deligne-Mumford topological stack},
 if for every point $x$ in $\X$ and every open substack
 $\U \subseteq \X$
 containing $x$ there is an open substack $\V \subseteq \U$ containing
 $x$ such that
 $\V \cong [V/I_x]$, for some topological space $V$ with  an action of $I_x$
 that is mild at the (unique) fixed point of $V$ lying above $x$.
 \end{defn}

 Weak Deligne-Mumford topological
 stacks are exactly topological stack for $\mathbf{LF}$=local homeomorphisms.
  Weak Deligne-Mumford topological  stacks form a full 2-subcategory of
 $\mathsf{TopSt}$,
 which we denote by $\mathsf{WeakDM}$. Deligne-Mumford topological stacks
 form a full 2-subcategory of $\mathsf{WeakDM}$, which we denote by
 $\mathsf{DM}$.
 The Yoneda functor identifies
 $\mathsf{Top}$ with a full 2-subcategory of $\mathsf{DM}$.

   \begin{prop}{\label{P:basicDM}}
    Let $\X$ be a pretopological stack.
   \begin{itemize}

          \item[$\mathbf{i.}$] If $\X$ is (weak) Deligne-Mumford,
          then so is every embedded substack of $\X$.

          \item[$\mathbf{ii.}$] Let $\{\U_i\}_{i\in I}$ be a
          covering of $\X$ by open substacks. If every $\U_i$ is
          (weak) Deligne-Mumford,
          then so is $\X$.

          \item[$\mathbf{iii.}$] The 2-category $\mathsf{WeakDM}$ is
          closed under  2-fiber products. (See Section \ref{S:DMFiber}
          for the case of $\mathsf{DM}$.)
        \end{itemize}
   \end{prop}

 We can characterize \'etale groupoids whose quotient
 is Deligne-Mumford  as follows.

  \begin{prop}{\label{P:DMchart}}
    Let $[R \sst{} X]$ be an \'etale
    groupoid and set $\X=[X/R]$. The necessary and sufficient
    condition for $\X$ to be a
    Deligne-Mumford topological stack is that, for every $x \in X$
    and every open $U'$
    containing $x$, there exist an open $x \in U \subseteq U'$ such that
    the $[R|_U \sst{} U]$ is isomorphic to the action groupoid of
    an $I_x$-action on $U$.
  \end{prop}

   \begin{proof}
      One implication is obvious. Assume now that $\X$ is Deligne-Mumford.
      By definition, there is a topological space $V$ with a mild $I_x$ action
      such that $[V/I_x]$ is equivalent to an open substack of $\X$ containing
      $p(x)$. This gives a local homeomorphism $q\: V\to \X$. Since
      $p \: X \to \X$ is a local homeomorphism, we may assume,
      after possibly
       shrinking $V$ to a smaller
      $I_x$-invariant open neighborhood of $x$
      (that we can, since $I_x$ acts mildly), that there is a
      lift $\tilde{q} \: V \to X$ of $q$
      to $X$ (as a pointed map). Since $q=p\circ \tilde{q}$ and $p$ are local
      homeomorphism, $\tilde{q}$ is also a local homeomorphism. By shrinking
      $V$  again around $x$, we may assume that $\tilde{q}$ is an open
      embedding. The image $U=\tilde{q}(V) \subseteq X$ is
      the open neighborhood
      of $x$ we were after.
      We can make $U$ as small as we want by shrinking $V$.
   \end{proof}

    \begin{cor}{\label{C:pd1}}
      Let $G$ be a discrete group acting on a topological space $X$.
      Then $[X/G]$ is  Deligne-Mumford if an only if the action
      is properly discontinuous (Definition \ref{D:pd}).
    \end{cor}

     \begin{proof}
        If the action is properly discontinuous, $[X/G]$ is obviously
        Deligne-Mumford. Assume now that $[X/G]$ is Deligne-Mumford.
        Take an arbitrary point $x \in X$.
        Consider the action groupoid $[G\x X \sst{} X]$.  The open neighborhood
        $U$ of $x$ constructed in the Proposition \ref{P:DMchart} satisfies
        the property required in the definition \ref{D:pd}.
     \end{proof}

      \begin{cor}{\label{C:pd2}}
        A pretopological stack is Deligne-Mumford if and only if
        it can be covered by open
        substacks of the form $[X/G]$, where $G$ is a discrete group
        acting properly discontinuously   on a topological space $X$.
      \end{cor}

     \begin{cor}{\label{C:finite}}
          Let $G$ be a finite group acting on a topological
    space $X$. Then $[X/G]$ is a Deligne-Mumford topological stack.
     \end{cor}

  \begin{prop}{\label{P:strong}}
   Let $[R \sst{} X]$ be an \'{e}tale groupoid such that
    the diagonal map $R \to X\x X$ is a closed map onto its image (where
    the image is endowed with the subspace topology from $X\x X$) and that,
     for every $x \in X$, the stabilizer group $I_x$ is finite.
     Assume further that $X$ is locally connected. Then
     $[X/R]$ is a  Deligne-Mumford topological stack.
     \end{prop}
  \begin{proof}
       Let $x$ be an arbitrary point on $X$, and let $H=I_x$ be its stabilizer
       group (which is finite). First we consider a special case.

\vspace{0.1in}
       \noindent {\em Special case.} Assume our groupoid has
       the following special form: $R=\coprod_{h \in H} R_h$
       such that $R_1$ is the identity section, and that
       the restriction of $t$ and $s$ to each $R_h$ is an open embedding from
       $R_h$ to $X$. Note that the image under $s$ and $t$ of every $R_h$
       contains $x$. Let $V$ be an open neighborhood of $x$ that is contained
       in  every $s(R_h)$. Call the pre-image under $s$ of $V$ in $R_h$ by
       $V_h$. Set $U'=\cap t(V_h)$. Then, $U'$ is an invariant neighborhood of
       $x$, and so is the connected component $U$ of $x$ in $U'$.
       The restriction of $R$ over $U$ is a groupoid of the form
       $[H\x U \sst{} U]$, such that the restriction of the source
       and target maps to each
       layer $\{h\}\x U$ of it are isomorphisms onto $U$. Lemma
       \ref{L:actiongpd} implies that
       this groupoid is in fact the groupoid associated to an action
       of $H$ on $U$.

\vspace{0.1in}
       \noindent {\em General case.} We reduce  to the special
       case considered above. Recall that $H=I_x$ is a subset of $R$, and that
       the source map $s \: R \to X$ is a local homeomorphism.
       Since $H$ is finite, we can find a neighborhood $W$ of $x$, and
       neighborhoods $W_h$ for each $h \in H$, such that $s$ induces
       a homeomorphism from $W_h$ to $W$. Let $A=R-\cup W_h$. Since
       $\De \: R \to X\x X$ is closed onto its image, $\De(A)$ is a closed
       subset of $Y:=\De(R)$. Note that, by the construction of $A$, $(x,x)$
       is not in $\De(A)$, so $Y-\De(A)$ is an  open neighborhood of
       $x$ in $Y$. Hence, there is
       an open neighborhood $U$ of $x$ in $X$ such that $U\x U$ is contained
       in $Y-\De(A)$, or equivalently, $U\x U$ does not intersect $\De(A)$.
       It is easy to check that the restriction $R_U$ of $R$ to $U$
       (which is defined to be $R_U=\De^{-1}(U\x U)$, with the same source and
       target maps) is a groupoid which satisfies the property required in the
       special case above. The reason for this is that, indeed, $R_U$ can
       be identified with an open subset of $\coprod W_h$, and the source
       and target maps for $R_U$ are the ones induced from $\coprod W_h \to X$.
 \end{proof}

   Proposition \ref{P:strong} can be rephrased as follows.

  \begin{prop}{\label{P:strong2}}
    Let $\X$ be a locally connected (Definition \ref{D:lcDM})
    weak Deligne-Mumford topological stack. Assume every
    point in $\X$ has  a finite inertia group, and that the diagonal
    $\X \to \X\x \X$ is a  closed map onto its image (see Definitions
    \ref{D:local} and \ref{D:representable2}
    and the ensuing example).
    Then $\X$ is a Deligne-Mumford stack.
 \end{prop}

   For an example of a weak Deligne-Mumford topological stack that is not
    Deligne-Mumford see Example \ref{E:vistoli}. For more examples of
    Deligne-Mumford topological stacks  see
    Section \ref{S:examples}.

\section{Fiber products of Deligne-Mumford topological stacks}
 {\label{S:DMFiber}}

  Unfortunately, the 2-category of Deligne-Mumford topological stacks
  fails to be closed under fiber products. But the failure is not so
  dramatic. We will see that,
   under some mild extra hypotheses, fiber products
  of Deligne-Mumford topological stacks will again be
  Deligne-Mumford
  (Corollary \ref{C:DMfiber}). The counterexamples could safely be regarded
  as pathological.

    \begin{defn}{\label{D:locallyquotient}}
      We say that a pretopological stack is  {\bf locally discrete-quotient}
      if it can be covered by open substacks each of which  equivalent
      to the quotient stack of a discrete group acting on a
      topological space.
    \end{defn}

  A locally discrete-quotient stack is  weak Deligne-Mumford,
   but it
  fails to be  Deligne-Mumford in general (e.g. pick a non-properly
  discontinuous group action and take its quotient).

    \begin{defn}{\label{D:controlable}}
      Let $f \: \Y \to \X$ be a map of locally discrete-quotient stacks.
      We say $f$ is {\bf controlable}, if for every point $y \in \Y$,
      there is an open substack $\V$ around $y$ and an open substack $\U$
      around $f(y)$ such that  $f$ maps $\V$ into $\U$ and
       \begin{itemize}

        \item there are presentation $\V=[V/H]$ and $\U=[U/G]$ as
           quotient stacks by discrete group actions;

        \item for a suitable choice of such presentations,
          the restriction $f|_{\V} \: \V \to \U$ is induced by a
          group homomorphism $\varphi \: H \to G$ and a
           $\varphi$-equivariant map $V \to U$.

       \end{itemize}
    \end{defn}

   Remark that, it is not true in general that a map of quotient stacks
   $[V/H] \to [U/G]$ comes from a map of action groupoids
   $[H\x V \sst{} V] \to [G\x U \sst{} U]$. Also it is {\em not} true in general
   that a map of action groupoids   $[H\x V \sst{} V] \to [G\x U \sst{} U]$
   is induced by a group homomorphism $\varphi \: H \to G$ and a
   $\varphi$-equivariant map $V \to U$.

   If $f \: \Y \to \X$ is controlable and $\Y$ is a
    Deligne-Mumford topological stack,
   then it can be shown 
   that
   the second condition of the above definition holds for the presentation
   $\V=[V/I_x]$ too.

   \begin{prop}{\label{P:locallyquotient}}
     Consider a diagram
         $$\xymatrix@=12pt@M=10pt{    &  \Y\ar[d]^f  \\
                \Z  \ar[r]_(0.45)g        &        \X   }$$
     of pretopological stacks, and assume $\Y$ and $\Z$ are
     locally discrete-quotient stacks.
     Then so is $\Z\x_{\X}\Y$. In particular,
     the 2-category   of locally discrete-quotient stacks
     is closed under 2-fiber products.
   \end{prop}

   \begin{proof}
      We may assume that $\Y=[Y/H]$ and $\Z=[Z/K]$ are quotient stacks,
      where $H$ and $K$ are discrete groups, and that
      the above diagram is induced from a diagram of groupoids
          $$\xymatrix@=12pt@M=10pt{            &  [H\x Y \sst{} Y] \ar[d]^F  \\
                   [K\x Z \sst{} Z]  \ar[r]_(0.53)G  &  [R \sst{} X]   }$$
      Before proceeding with the proof, we fix some notations.
      Let $F_0 \: Y \to X$ and $F_1 \: H\x Y \to R$
      denote the structure maps of $F$. For each $h \in H$, let $F_h \: Y \to R$
      be the composition
        $$Y \risom \{h\}\x Y \hra H\x Y \llra{F_1} R.$$
      Define $G_0 \: Z \to X$, $G_1 \: K\x Z \to R$, and $G_k \: Z \to R$
      similarly.

      The construction of Section \ref{S:fiber} shows that the the fiber product
      $\Z\x_{\X}\Y$ is the quotient stack
      of the groupoid $[P_1 \sst{} P_0]$, where
      $P_0=Y\x_{\X}Z=(Y\x Z)\x_{X\x X}R$ and
      $P_1=(H\x K)\x P_0$. The source map is simply the
      projection onto the second factor. This already suggests that
      $[P_1 \sst{} P_0]$ is the action groupoid of an action of $H\x K$ on
      $P_0$. We give the action explicitly, using the description
      of the target map of the groupoid $[P_1 \sst{} P_0]$ given in   Section
      \ref{S:fiber}. Let $(h,k) \in H\x K$ be an arbitrary element.
      The action of $(h,k)$ on $P_0=(Y\x Z)\x_{X\x X}R$ sends a triple
      $(y,z,r)$ to
        $$\big(\,h(y),\,k(z),\,[F_{h^{-1}}\big(h(y) \big)]\cdot[r]\cdot[G_k(z)]\,\big)$$
    where the elements in the square brackets belong to $R$, and $\cdot$
    stands for composition of composable arrows in  $R$. This
    is easily seen to be a group action.
   \end{proof}

   The above proposition implies that the fiber product of Deligne-Mumford
   topological stacks is always a locally discrete-quotient stack.
   The next proposition
   tells us when this locally discrete-quotient stack is   Deligne-Mumford.

    \begin{thm}{\label{T:DMprod}}
      Consider a diagram
         $$\xymatrix@=12pt@M=10pt{    &  \Y\ar[d]^f  \\
                \Z  \ar[r]_(0.45)g        &        \X   }$$
     of locally discrete-quotient stacks, and assume $\Y$ and $\Z$
     are Deligne-Mumford.
     Also, assume  that $f$ and $g$ are controlable. Then  the 2-fiber product
     $\Z\x_{\X}\Y$
     is a Deligne-Mumford topological stack.
    \end{thm}

    \begin{proof}
       We may assume that $\Y=[Y/H]$
     (resp. $\Z=[Z/K]$), for some topological space $Y$ (resp. $Z$)
     with a properly discontinuous $H$ (resp. $K$) action
     (Corollary \ref{C:pd1}), and that
     $\X=[X/D]$, where $D$ is a discrete group acting on $X$.
     We may also assume that our diagram of stacks is induced from the following
     diagram of groupoids
         $$\xymatrix@=12pt@M=10pt@C=40pt{                    &  [H\x Y \sst{} Y]
                          \ar[d]^{F=(\varphi, F_0)}  \\
                   [K\x Z \sst{} Z]  \ar[r]_{ G=(\psi, G_0)}
                                                       &  [D\x X \sst{} X]   }$$
     where the vertical map is induced from a group homomorphism
     $\varphi \: H \to D$, and a $\varphi$-equivariant map $F_0 \: Y \to X$,
     and    the horizontal map is induced from a group homomorphism
     $\psi \: K \to D$, and a $\psi$-equivariant map $G_0 \: Z \to X$.
     Set $P:=Y\x_{\X}Z=(Y\x Z)\x_{X\x X}(D\x X)$ (this is what we called
     $P_0$ is Section \ref{S:fiber}).   In the course of proof of
     Proposition \ref{P:locallyquotient} we showed that there is an
     action of $H\x K$ on $P$ such that the quotient stack
     $[P/H\x K]$ is equivalent to $\Z\x_{\X}\Y$. It is enough to show that
     this action is properly discontinuous.
     Pick   an arbitrary point $p_0 \in P=(Y\x Z)\x_{X\x X}(D\x X)$.
     This point can be represented as
     $p_0=(y_0,z_0,x_{0s},x_{0t},d_0,x_0)$, where
     $y_0\in Y$, $z_0\in Z$, $x_{0s},x_{0t},x \in X$, and $d_0\in D$. Note that,
     $F_0(y_0)=x_{0s}$, $G_0(z_0)=x_{0t}$, $x_0=x_{0s}$, and
     $d_0(x_{0s})=x_{0t}$, where
     $d_0 \: X\to X$ is the automorphism coming from the
     action of $D$ on $X$. Denote by $I_{y_0}$, $I_{z_0}$, $I_{x_{0s}}$ and
     $I_{x_{0t}}$
     the stabilizer groups of the corresponding group actions at the
     respective points. Observe that, since $d(x_{0s})=x_{0s}$, conjugation
     by $d_0$ maps $I_{x_{0s}}$ isomorphically to $I_{x_0t}$. So we {\em identify}
     $I_{x_{0s}}$ and $I_{x_{0t}}$ via this isomorphism and call them both $I$.

     Having made these remarks, we can now shorten the notation for
     a point $p=(y,z,,x_s,x_t,d,x) \in  P$ by simply denoting it by
     $p=(y,z,d)$.

     We have natural group homomorphisms $I_{y_0} \to I$ and $I_{z_0} \to I$,
     induced by $\varphi$ and $\psi$, respectively.
     We claim the following:

       \begin{itemize}

          \item[$\mathbf{a.}$] The stabilizer group at $p_0$
          of the action of
             $H\x K$ on $P$  is equal to the fiber product
             $I_{y_0}\x_{I}I_{z_0} \subseteq H\x K$. (Note that
             $$I_{y_0}\x_{I}I_{z_0}
                    =\{(h,k) \in I_{y_0}\x I_{z_0}\vert \ \
                                   d_0\;\psi(k)\;d_0^{-1}=\varphi(h)\}.)$$

          \item[$\mathbf{b.}$] There exists an
          $I_{y_0}\x_{I}I_{z_0}$-invariant
          open neighborhood $U$ of $p_0$ in $P$ such that the action of
          $I_{y_0}\x_{I}I_{z_0}$ on $P$ is mild at $p_0$,
          and for any $(h,k) \in H\x K$
          not in  $I_{y_0}\x_{I}I_{z_0}$, $(h,k)(U)\cap U=\emptyset$.
       \end{itemize}
       Such open neighborhood $U$ of $p_0$  satisfies the property required
       in Definition \ref{D:pd}, hence, the proposition
       follows.

     \vspace{0.1in}
       \noindent {\em Proof of} ($\mathbf{a}$).
          If we decipher  the description of the
        action  of $H\x K$ given in Proposition \ref{P:locallyquotient}
        in the case where the groupoid $[R \sst{} X]$ is the action groupoid
        $[D\x X \sst{} X]$, we see that
        an arbitrary $(h,k) \in H\x K$ sends
        $p=(y,z,d)$ to
        $$\big(\,h(y),\,k(z),\, \varphi(h^{-1})\;d\;\psi(k)\,\big)$$
        This immediately implies ($\mathbf{i}$).

    \vspace{0.1in}
    \noindent {\em Proof of} ($\mathbf{b}$). Set
         $$U=\{(y,z,d) \in P\vert\ d=d_0\}.$$
         This an open neighborhood of $p_0=(y_0,z_0,d_0)$.
         By    definition of the fiber product, $I_{y_0}\x_{I}I_{z_0}$ is
         the set of pairs $(h,k) \in H\x K$ such that
         $\varphi(h^{-1})d_0\psi(k)=d_0$. So if $(h,k)$ {\em is not} in
         $I_{y_0}\x_{I}I_{z_0}$, the third component of any element in
         $(h,k)(U)$ is different from $d_0$, so  $(h,k)(U)\cap U=\emptyset$.
         Similarly, if $(h,k)$ {\em is} in
         $I_{y_0}\x_{I}I_{z_0}$, we have  $(h,k)(U)=U$.
         All that remains to be proved is that the action of
         $I_{y_0}\x_{I}I_{z_0}$
         on $U$ is mild. To prove this, note that we have an
         $I_{z_0}\x_{I}I_{y_0}$-equivariant
         isomorphism
         $U \risom Y\x Z$ (forgetting the factor $d$).
         Since the actions of $I_{y_0}$ on $Y$
         and $I_{z_0}$ on $Z$ are mild, so is the action of $I_{z_0}\x I_{y_0}$
         on $Y\x Z$. In particular, the action of the subgroup
         $I_{z_0}\x_{I}I_{y_0}$ is also mild. Therefore, the action of
         $I_{z_0}\x_{I}I_{y_0}$ on $U$ is mild.
         The proof of the proposition is now complete.
   \end{proof}


   The above proposition is useless unless we know how to determine if a
   map of locally discrete-quotient stacks is controlable. The
   Proposition \ref{P:controlable} tell us that it is the case in quite
   general situations.

    \begin{defn}{\label{D:lcDM}}
      Let $\X$ be a pretopological stack. We say $\X$ is
      {\bf locally connected}, if there exists a
      chart $p \: X \to \X$  such that $X$
      is a locally connected topological space (also see
      Section \ref{SS:connectivity}).
   \end{defn}

      \begin{lem}{\label{L:lcDM}}
        Let $\X$ be a locally connected weak Deligne-Mumford stack.
        Then for any \'{e}tale chart $p \: X \to \X$, $X$ is locally
        connected.
      \end{lem}

    \begin{proof}
      See Lemma \ref{L:lpcDM}.
    \end{proof}

    \begin{prop}{\label{P:controlable}}
       Let $f \: \Y \to \X$ be a map of stacks, where $\X$ is a locally
       discrete-quotient stack and $\Y$ is Deligne-Mumford.
       Then $f$ is controlable, if at least one  of the following holds:
        \begin{itemize}
          \item $\X$ is a topological space;

          \item  $\Y$ is locally connected
                    (see Definition \ref{D:lcDM});

          \item  $\Y$ has a finitely generated inertia group at every point.

        \end{itemize}
    \end{prop}

    \begin{proof}
       The case where $\X$ is a topological space is obvious. We do the other
       two cases.

       Pick a point $y$ on $\Y$, and let $H=I_y$. We may assume $\Y=[Y/H]$,
       where $Y$ has an $H$-action with a unique
       fixed point $y_0$, lying above $y$, at which $H$ acts mildly.
       Assume that either $Y$ is connected
       and locally connected (Lemma \ref{L:lcDM}),
       or that $H$ is finitely generated.
       We may also assume that $\X=[X/G]$, with $G$ discrete.
       Finally, by shrinking $Y$ around $y_0$, we may assume that
       the composite map $Y \to \X$ lifts to $X$; that is, $f$
       is induced by a map of action groupoids
       $F \: [H\x Y \sst{} Y] \to  [G\x X \sst{} X]$. The notation
       $F_0 \: Y \to X$ and $F_1 \: H\x Y \to G\x X$ should be clear.
       We consider the two cases
       separately.

\vspace{0.1in}
       \noindent{\em $Y$ is connected.} Since $Y$ is connected,
       under the map $F_1 \: H\x Y \to G\x X$, every layer $\{h\}\x Y$, $h \in H$,
       should factor through  a layer $\{g\}\x X$, for some $g \in G$.
       This gives a map $\varphi \: H \to G$ that induces the original map
       $F\: [H\x Y \sst{} Y] \to  [G\x X \sst{} X]$ of groupoids. It is easily
       seen  that $\varphi$ is a group homomorphism, and the map $Y \to X$
       is $\varphi$-equivariant.

\vspace{0.1in}
       \noindent{\em $H$ is finitely generated.} Let $x_0=F_0(y_0)$ be the image
       of $y_0$ in $X$.
       Let $h_1,h_2,\cdots, h_n$ be generators for $H$.
       For each $i=1,2,\cdots, n$, the map $F_1 \: H\x Y \to G\x X$
       maps the point $(h_i,y_0)$   to some
       $(g_i,x_0)$,  $g_i \in G$. Let $V_i \subseteq Y$ be a neighborhood
       of $y_0$ such that $F_1(\{h_i\}\x V_i) \subseteq \{g_i\}\x X$.
       Let $V$ be an $H$-invariant neighborhood of $y_0$ contained
       in $\bigcap U_i$. Then $F_1(\{h_i\}\x V) \subseteq \{g_i\}\x X$,
       for every $i$. But since $F_1$ comes from a map of groupoids and
       $h_i$ generate $H$,
       it follows that, for every $h \in H$, we have
        $F_1(\{h \}\x V_i) \subseteq \{g\}\x X$, for some $g \in G$.
        That is, we have a map $\varphi \: H \to G$ that induces
         the original map
       $F\: [H\x Y \sst{} Y] \to  [G\x X \sst{} X]$ of groupoids. It is easily
       seen that $\varphi$ is a group homomorphism, and the map $Y \to X$
       is $\varphi$-equivariant.
     \end{proof}

    \begin{cor}{\label{C:DMfiber}}
     Consider a diagram
        $$\xymatrix@=12pt@M=10pt{    &  \Y\ar[d]^f  \\
             \Z  \ar[r]_g        &        \X   }$$
     of locally discrete-quotient stacks, and assume $\Y$ and $\Z$ are
     Deligne-Mumford. Assume at least one of the following holds:

 \begin{itemize}

          \item $\X$ is a topological space;

          \item each of $\Y$ and $\Z$
              is either locally connected (see Definition \ref{D:lcDM})
              or has finitely generated inertia groups.

        \end{itemize}
      Then   $\Z\x_{\X}\Y$ is a Deligne-Mumford topological stack.
    \end{cor}

   One can construct examples where a 2-fiber product of Deligne-Mumford
   topological stacks is not Deligne-Mumford. I will give a flavor
   of how such an example is constructed, but refrain from giving the messy
   details. Let $\mathbb{F}$ be the free group on countably
   many generators. Take $\X=\B\mathbb{F}=[*/\mathbb{F}]$ and $\Z=*$.
   The interesting one is $\Y$, whose construction is a bit more involved.
   Let $A=\{0\}\cup\{\frac{1}{n}\ |\ n \in \bbN\}$.
   Set $Y=\mathbb{F}\x A/\mathbb{F}\x \{0\}$  with the topology
   that is discrete on the complement
   of the distinguished point (which we  call $O$), and has a local
   basis at $O$    consisting of ``balls of radius $\frac{1}{n}$''
   (these are all point, in various copies of $A$ appearing in $Y$, whose
   distance to $0$ is less than  $\frac{1}{n}$). Notice that this topology
   is weaker than the quotient topology on $\mathbb{F}\x A/\mathbb{F}\x \{0\}$.
   There is an obvious {\em mild} ``rotation''
   action of $\mathbb{F}$  on $Y$ fixing  $O$
    and acting freely on the rest.
   Define $\Y=[Y/\mathbb{F}]$. Outside the distinguished point, $\Y$
   is equivalent to the discrete set $\{\frac{1}{n}\ |\ n \in \bbN\}$.
   At the distinguished
   point $\Y$ looks like $\B\mathbb{F}$
   (i.e. $\Gamma_O\cong \B\mathbb{F}$). What is interesting is the way
   the set $\{\frac{1}{n}\ |\ n \in \bbN\}$ is ``approaching'' $\B\mathbb{F}$.

   The fun part is  to construct an {\em uncontrolable} map $f \: \Y \to \X$
   which satisfies certain nice properties. I will not get into that, but the
   outcome is that, the fiber product $\Z\x_{\X}\Y$ (which can
   be thought of as the fiber of $f$ over a point)  will turn out to be a
   quasitopological space that is not a topological space. In particular,
   it can not be a Deligne-Mumford topological stack.

   This example also shows that   a covering
   spaces  of a non locally connected Deligne-Mumford topological
   stack is not necessarily a Deligne-Mumford stack (see Proposition
   \ref{P:strongcovering}),
   because in this case $\Z\x_{\X}\Y$ is a covering space of $\Z$.

\section{Gluing topological stack along  substacks}{\label{S:gluing}}

 Throughout this section we fix a class $\mathbf{LF}$ of local fibrations
 (Section \ref{SS:fibrations}). By a substack we mean an embedded substack.

 In point set topology, given an equivalence relation on a topological space,
 one can always construct a quotient space,
 satisfying the expected universal property.
 It turns out that in the world of pretopological stacks this is far
 from being true; even sticking to the special case of  Deligne-Mumford does
 not help much. More depressing is
 the fact that even the
 good old
 quotients of topological {\em spaces} may loose their universal property when
 viewed in the 2-category of topological stack.
 So we have got to be very careful with our intuition.
 Fortunately,
 by choosing the right class of local fibrations,
 we can fine-tune our category of topological stacks so that
 our desired gluing and quotient
 constructions can be performed
 within the realm of topological stacks.

 It turns out that, for the purpose of homotopy theory, the required conditions
 on $\mathbf{LF}$ are quite mild (Section \ref{S:homotopy}).

 The following lemma is essentially due to Angelo Vistoli.

  \begin{prop}{\label{P:gluing1}}
   Let $L$, $Y$ and $Z$ be  topological spaces. Suppose we are given
   a push-out diagram in the category of topological spaces
     $$\xymatrix@=12pt@M=10pt@M=10pt@=12pt{
               L \ar@{^(->} [r]^i \ar@{^(->}[d]_j &  \ar[d] Y  \\
                             Z \ar[r]        &    Y\vee_L Z       }$$
   with $i$ and $j$ embeddings. Suppose further that $i \: L \hra Y$ is  LTC
   (see Definition \ref{D:LTC}).
   Then, this diagram remains
   a push-out in the 2-category of topological stacks. In other words,
   for every topological stack $\X$,
   given a pair of maps $f  \: Y \to \X$ and $g  \: Z \to \X$   that are
   compatible along   $L$ (i.e. we are given an identification
   $\varphi \: f|_L \twomor{} g|_L$),  $f$ and $g$ can be glued in
   a unique (up to a unique identification) way to a
   map $f \vee_L g \: Y\vee_L Z \to \X$.
\end{prop}

 Before giving the proof, let us make the statement precise.
 Let $\Y$, $\Z$ and $\LL$ be pretopological stacks, and $i \: \LL \hra \Y$
 and $j \: \LL \hra \Z$ be embeddings.
 Define the  groupoid $\Glue(\Y,\Z,\LL; \X)$ of
 {\em gluing data}
 as follows: \label{pushout}

 {\small
       $$\Ob\Glue(\Y,\Z,\LL; \X)=\left\{\begin{array}{rcl}
          (f,g,\varphi)
               & | & f\: \Y \to \X, \ g\: \Z \to \X, \ f|_{\LL}
                                  \twomor{\varphi} g|_{\LL}
                                            \end{array}\right\}$$
\vspace{0.1in}
  $$\Mor_{\Glue(\Y,\Z,\LL; \X)}\big( (f,g,\varphi), (f',g',\varphi') \big)=
           \left\{\begin{array}{rcl}
              (\al,\be) & | &  f \twomor{\al} f', \ g \twomor{\be} g' \
                                    \text{identifications}\\
                        &   &  \text{such that}  \hspace{0.2in}
         \xymatrix@R=10pt@C=12pt{
             f|_L \ar @{=>}[r]^{\varphi}
                    \ar@{=>}[d]_{\al|_{\LL}} \ar@{}[rd]|{\circlearrowleft}
                                             &  g|_L \ar@{=>}[d]^{\be|_{\LL}} \\
                 f'|_{\LL}  \ar@{=>}[r]_{\varphi'}  &  g'|_{\LL}   }
                                                          \end{array}\right\}$$}

     Proposition \ref{P:gluing1} follows from the following theorem.

   \begin{thm}{\label{T:gluing}}
      Notation being as in Proposition \ref{P:gluing1}, consider the natural
      restriction functor
       $$\varrho \: \Hom_{\mathsf{St}}(Y\vee_L Z, \X) \to \Glue(Y,Z,L; \X)$$
           $$ H \mapsto (H|_{\Y}, H|_{\Z},\id). $$
       This functor is an equivalence of groupoids. Furthermore,
       if we do not assume
       that $i \: L \to Y$ is LTC, then $\varrho$ is still fully faithful.
   \end{thm}

  \begin{proof}
     Fix an $\mathbf{LF}$ chart
     $p \: X \to \X$.  Let $[R \sst{} X]$
     be the corresponding groupoid.

    \vspace{0.1in}
        \noindent {\em Proof of faithfulness}. Let $H,H' \: Y\vee_L Z \to \X$ be
          morphisms, and $\Psi_1,\Psi_2 \: H \twomor{ }H'$
           identifications between them.
          Assume $\Psi_1|_Y=\Psi_2|_Y$ and $\Psi_1|_Z=\Psi_2|_Z$. We have
          to show that  $\Psi_1=\Psi_2$. Since $\X$ is a stack, it is
          enough to do this locally. So, after replacing $Y\vee_L Z$ by an
          open covering, we may assume that $H$ and $H'$
          lift to $X$.  Denote the lifts
          by $h,h' \: Y\vee_L Z \to X$. (Recall that, a {\em lift} of $H$
          is a  pair $(h,\gamma)$, with $h  \: Y\vee_L Z \to X$
          and $\gamma \: p\circ h \twomor{} H$; we supress $\gamma$
          in the notation.)
          This gives
          a map $(h,h')\: Y\vee_L Z \to X\x X$.
          The 2-isomorphisms $\Psi_1$ and $\Psi_2$
          precisely  correspond to lifts of this map to $R$:
              $$\xymatrix@=22pt@M=8pt{   & R  \ar[d] \\
                    Y\vee_L Z \ar@/^/ [ur]^{\psi_1}  \ar@/_/ [ur]_{\psi_2}
                                \ar[r]_{(h,h')}   &     X\x X     }$$
          Since $\Psi_1|_Y=\Psi_2|_Y$, we have $\psi_1|_Y=\psi_2|_Y$.
          Similarly, we have $\psi_1|_Z=\psi_2|_Z$. Therefore, since $Y\vee_L Z$
          is a push out in the category of topological spaces, we have
          $\psi_1=\psi_2$. This implies that $\Psi_1=\Psi_2$.

    \vspace{0.1in}
        \noindent {\em Proof of fullness}.
          Let $H,H' \: Y\vee_L Z \to \X$  be
          morphisms. An arrow  from $\varrho(H)$
          to $\varrho(H')$
          in the groupoid $\Glue(Y,Z,L;\X)$ is a pair of identifications
          $\Psi_Y \:H|_Y \twomor{}  H'|_Y$ and
          $\Psi_Z \:H|_Z \twomor{}  H'|_Z$ such that $\Psi_Y|_L=\Psi_Z|_L$.
          We have to show that there exists an identification
          $\Psi \: H \twomor{} H'$ such that $\Psi|_Y=\Psi_Y$
          and $\Psi|_Z=\Psi_Z$.

          Having proved the uniqueness of $\Psi$ in the previous part,
          it is now enough to prove
          the existence locally (because $\X$ is a stack).    So, again, we may
          assume that we have lifts $h,h' \: Y\vee_L Z \to X$
           of $H$ and $H'$. We want to construct
           a lift $\psi$ as in the following diagram:
            $$\xymatrix@=22pt@M=8pt{   & R  \ar [d] \\
                    Y\vee_L Z \  \ar@{.>} [ur]^{\psi}  \ar[r]_{(h,h')}
                                                            &     X\x X     }$$
          We have a lift of $\psi_Y : Y \to R$ of $(h,h')|_Y$ obtained from
           the 2-isomorphism $\Psi_Y$, and a similar one $\psi_Z : Z \to R$
           obtained from $\Psi_Z$. Since the restrictions of $\Psi_Y$ and
           $\Psi_Z$ to $L$ are equal, we have $\psi_Y|_L=\psi_Z|_L$.
           Therefore, $\psi_Y$ and $\psi_Z$ glue to a map
           $\psi \: Y\vee_L Z \to R$. The 2-isomorphism
           $\Psi \: H \twomor{} H'$
           corresponding to $\psi$  is what we were looking for.

        \vspace{0.1in}
          This much of the proof works for any pretopological stack $\X$,
          without  LTC assumption on $i \: L \hra Y$.

        \vspace{0.1in}
        \noindent {\em Proof of essential surjectivity}.
        Assume we are given a gluing datum
           $$(f,g,\varphi) \in \Ob\Glue(Y,Z,L; \X).$$
        We have to show that there exists $H \: Y\vee_L Z \to \X$
        such that $(H|_Y,H|_Z,\id)$ is isomorphic to
        $(f,g,\varphi)$ in $\Glue(Y,Z,L; \X)$.

        In view of the
        previous parts, it is enough to prove the
        existence locally.  In other words, we have to show that
        for every point $t \in L$, there is an open
        neighborhood $U$ of $t$ in $Y\vee_L Z$
        over which the gluing is
        defined. (If $t \in (Y\vee_L Z)\backslash L$  the existence is obvious.)
        By shrinking  $Y$ and $Z$ around $t$ we may assume that
        that $f$ and  $g$ have lifts  $\tilde{f} \: Y \to X$ and
        $\tilde{g} \: Z \to X$.
        The identification $f|_L \twomor{\varphi} g|_L$
        corresponds to a lift $\phi$ as in the following diagram:
         $$\xymatrix@=22pt@M=8pt{   & R  \ar [d] \\
                    L \ar    [ur]^{\phi}  \ar[r]_(0.4){(\tilde{f}|_L,\tilde{g}|_L)}
                                                            &     X\x X     }$$
        Consider the following commutative diagram:
           $$\xymatrix@=12pt@M=10pt{
             L  \ar[r]^{\phi} \ar[d]_i &  R\ar[d]^s  \\
             Y \ar[r]_{\tilde{f}}  \ar@{.>}[ur]^{F}      &   X        }$$
        The map $i \: L \hra Y$ is LTC and $s$ is LF,
        so by shrinking $Y$ around $t$  we can fill the dotted line.
        Shrink $Z$ accordingly, so that $Y$ and $Z$ are
          intersect $L$ in the same open set. (This can be done, since
          $L$ is a subspace of $Z$.)
        The maps $t\circ F \: Y \to X$ and
        $\tilde{g} \: Z \to X$
        now agree (on the nose)  along $L$, so they glue to a map
        $h \: Y\vee_L Z \to X$. The composition $H=p\circ h$ is the desired
        gluing of $f$ and $g$.
  \end{proof}


  It is an interesting question what kinds of colimits
  in $\mathsf{Top}$ remain a ``colimit'' in $\mathsf{TopSt}$.
  It is easy to construct examples where this fails.

   \begin{ex}
    Consider the  push-out   diagram
     $$\xymatrix@=18pt@M=10pt{
      \mathbb{R}^2 \ar[r]^{x \mapsto -x} \ar[d]_{\id} & \mathbb{R}^2 \ar[d]  \\
             \mathbb{R}^2  \ar[r]     &    \mathbb{R}^2/\mathbb{Z}_2    }$$
     This diagram is no longer a push-out diagram in  $\mathsf{TopSt}$.
     For example, take $\X=[\mathbb{R}^2/\mathbb{Z}_2]$ and let
     $f=g \: \mathbb{R}^2 \to \X$  be the obvious map. Then,
     there is no induced map on the colimit $\mathbb{R}^2/\mathbb{Z}_2 \to \X$.
     (Hint: if such map existed, it would have a lift to $\bbR^2$
     around a small enough neighborhood of
     $0 \in \mathbb{R}^2/\mathbb{Z}_2=\mathbb{R}^2$.)
    \end{ex}

  The next thing we want to consider is the question of
  existence of push outs in the 2-category of topological stacks.
 We formulate the question in  the following form.

 \vspace{0.1in}
 \noindent{\bf Gluing Problem.}
  Let $\LL$, $\Y$ and $\Z$ be  topological stacks. Suppose we are given
   a  diagram
     $$\xymatrix@=12pt@M=10pt@M=10pt@=12pt{
                 \LL \ar@{^(->} [r]^i \ar@{^(->}[d]_j &    \Y  \\
                        \Z          &        }$$
   with $i$ and $j$   embeddings.
   Can we can glue $\Y$ and $\Z$ along $\LL$
   to obtain a topological stack $\Y\vee_{\LL}\Z$?
 \vspace{0.1in}

 A few remarks on this gluing problem are in order. First of all,
 the assumptions on the maps $i$ and $j$ being embeddings is very
 crucial. Not all push-outs  exist in $\mathsf{TopSt}$.
 For instance, if $\X$ is
 a non-trivial $\bbZ_2$-gerbe over $S^1$, it can be shown that
 $I\x\X\big/\{1\}\x\X$ can not exist in any reasonable sense. In particular,
 {\em suspensions} of stacks do not always exist.

 The second thing that should be made clear is
 what  it means to glue  pretopological stacks! We require
 the glued pretopological
 stack $\Y\vee_{\LL}\Z$ to satisfy the following properties:

  \begin{itemize}

          \item[$\mathbf{G1}$] There are substacks $\Y'$ and $\Z'$
            and $\LL' := \Y' \cap \Z'$ of $\Y\vee_{\LL}\Z$,
             together with equivalences  $\Y \risom \Y'$,
             $\Z' \risom \Z$ and $\LL \risom \LL'$ that fit into
             a commutative 2-cell
           $$\xymatrix@=6pt@M=6pt{  & \Y  \ar[rd] \ar@{=>} [d] & \\
          \LL \ar[rr] \ar[ru]^i \ar[dr]_j & \ar@{=>} [d] &  (\Y\vee_{\LL}\Z)\\
                    & \Z \ar[ur] &           }$$


          \item[$\mathbf{G2}$] The equivalences in  $\mathbf{G1}$
          induce an equivalence
            $$(\Y\backslash\LL)\coprod(\Z\backslash\LL) \risom
                                (\Y\vee_{\LL}\Z)\backslash\LL'.$$

        \end{itemize}

   When $\Y$, $\Z$ and $\LL$ are topological spaces, it is not difficult
   (but not so trivial either) to show that these conditions characterize
   $\Y\vee_{\LL}\Z$ up to homeomorphism. In the case of pretopological stacks,
   however, it is presumably not true in general. So, in this generality,
   there might be several ways of gluing  $\Y$ and $\Z$ along $\LL$.
   We will see that, once we impose certain conditions on the stacks and
   on the
   embeddings $i$ and $j$, then we can state some uniqueness results.

  The reason for formulating the gluing problem using the rather strange looking
  conditions  $\mathbf{G1}$ and $\mathbf{G2}$,
  instead of a universal property,
  is that, first of all, these two conditions
  are obviously necessary, and second of all,
  a gluing $\Y\vee_{\LL}\Z$ always exists in this sense (but it
  may not satisfy the universal property).

   The question whether $\Y\vee_{\LL}\Z$ does indeed satisfy
  the universal property of a push out
  can then be discussed as a separate
  problem (see Theorem \ref{T:gluing3}).

   \begin{thm}{\label{T:gluing2}}
     Let $\LL$, $\Y$ and $\Z$ be  pretopological stacks. Suppose we are given
   a   diagram
     $$\xymatrix@=12pt@M=10pt@M=10pt@=12pt{
                 \LL \ar@{^(->} [r]^i \ar@{^(->}[d]_j &    \Y  \\
                       \Z          &        }$$
      where $i$ and $j$ are embeddings. Then,
              $\Y$ and $\Z$ can be glued along $\LL$
               to give  a pretopological stack  $\Y\vee_{\LL}\Z$
              (see axioms $\mathbf{G1}$ and $\mathbf{G2}$ above).
   Furthermore, if   $\Y$ and $\Z$ are topological and
             admit neat LF charts with respect to $\LL$ (see Definition
             \ref{D:neat} below), then
             $\Y\vee_{\LL}\Z$ admits a chart $X \to \Y\vee_{\LL}\Z$ whose
             restrictions
             over both $\Y$ and $\Z$ are LF.
             (See Remark \ref{R:gluing} for a discussion about when
             $\Y\vee_{\LL}\Z$ becomes topological.)
 \end{thm}
  %
%

  To prove this theorem we need some preparation.

    \begin{defn}{\label{D:neat}}
     Recall that we have fixed a class  $\mathbf{LF}$ of local fibrations.
     Let $Y$ be a topological space and $L \subseteq Y$ a subspace.
     An epimorphic
      local fibration $f \: L' \to L$ is said to be {\bf extendable},
     if there is a topological space $Y'$ containing $L'$ as a subspace,
     and an epimorphic
      local fibration $\tilde{f} \: Y' \to Y$ extending $f$, such that
      $\tilde{f}^{-1}(L)=L'$. If such an extension
     is only possible
     after replacing $L'$ by an open covering,  we say
     $f$ is {\bf locally extendable}. We say $(Y,L)$ is a {\bf neat} pair
     (or $L$ is a neat subspace of $Y$), if every local fibration
     $f \: L' \to L$ is locally extendable.

     For a pair $(\Y,\LL)$ of topological stacks, we say
     an LF chart $Y \to \Y$ is  neat, if the corresponding
     subspace pair $(Y,L)$ is  neat.
    \end{defn}

    Being neat is not a very strong condition, as we see in the following
    example.

    \begin{ex}{\label{E:neat}}
    \end{ex}
    \begin{itemize}
           \item[$\mathbf{1.}$] The notions
           of  locally extendable and neat
            are both local
           on the target, in the sense that, given a family of open
           subspaces $\{U_i\}$ of $Y$ whose union contains  $L$,
           we can check
           local extendability,  or neatness, by doing so for every
           $L\cap U_i \hra U_i$.

          \item[$\mathbf{2.}$] Every subspace $L \subseteq Y$ that is a local
           retract is  neat.

             To see this, we use part
             ($\mathbf{1}$) to reduce to the case where $L \subseteq Y$
             is a retract. Let $r \: Y \to L$ be the retraction map.
             Then, for any local
             fibration $f \: L' \to L$, the projection map
             $Y\x_{L}L' \to Y$ is a local fibration that extends  $f$.
             Note that $L' \hra Y\x_{L}L'$ is again a retract.

          \item[$\mathbf{3.}$] Assume $\mathbf{LF}$ is (contained in) the class
          of locally cartesian maps (Example \ref{E:lf}.$\mathbf{3}$)
          or local homeomorphisms (Example \ref{E:lf}.$\mathbf{6}$),
          then any subspace $L \subseteq Y$ is  neat.

          \item[$\mathbf{4.}$] Let $\mathbf{LF}$ be the class of all
           epimorphisms. Then every subspace $L \subseteq Y$ is neat.

             To see this, let $f \: L' \to L$ be an epimorphism.
             Since local extendability is local on the target
             by part ($\mathbf{1}$), we may assume $f$ admits a section
             $s \: L \to L'$. We can now take  $Y' = L'\vee_{L}Y$,
             and extend $f$ over to $Y'$ in the obvious way.
        \end{itemize}

  \vspace{0.1in}

 \noindent{\bf Exercise.} Show that if $\Y$ is a weak Deligne-Mumford
 topological  stack and $\LL \hra \Y$ a neat substack, then the condition
 of Definition \ref{D:neat} is satisfied for {\em every} \'{e}tale chart
 $Y \to \Y$. (Hint: adopt the proof of Lemma \ref{L:lpcDM}.)
  \vspace{0.1in}

   The following simple lemma is used in the proof of
   Theorem \ref{T:gluing2}.

  \begin{lem}{\label{L:gluing}}\end{lem}{\em

            \begin{itemize}

              \item [$\mathbf{i.}$] Let $Y$ be a topological space and
              $L \subseteq Y$ a subspace.
              Let $\{U_i\}$ be an open covering of $L$. Then
              $\coprod U_i \to L$ is extendable. If $(Y,L)$
              is a neat pair, then the extension can also be chosen to
              be neat.

              \item [$\mathbf{ii.}$]
              Let $\LL \hra \Y$ be an embedding of topological stacks,
              and let $Y \to \Y$ be a neat chart with $L \subseteq Y$
              the (neat) subspace corresponding to $\LL$.
              Let $L' \to L$ be an epimorphic local fibration,
              and think of it as a chart for $\LL$. Then $L' \to \Y$
              is locally extendable to an LF chart  $Y' \to \Y$ for $\Y$.

            \end{itemize}
          }

\vspace{0.1in}
          \begin{proof}[Proof of part {\em ($\mathbf{i}$)}]
             For each $U_i$
             in the covering find an open $V_i$ in $Y$ such that
             $V_i\cap L=U_i$, and set $Y'=(\coprod V_i)\coprod (Y\backslash L)$.

\vspace{0.1in}
             \noindent {\em Proof of part} ($\mathbf{ii}$).
                 Since $L \subseteq Y$ is neat, we can assume,
                 after possibly replacing $L'$ by an open covering,
                 that $L' \to Y$ extends to an LF epimorphism
                  $Y' \to Y$,
                 for some topological space $Y'$ containing $L'$ as
                 a subspace. Composing this map with the chart $Y \to \Y$,
                 give the desired extension of $L' \to \LL$.
          \end{proof}

   \begin{proof}[Proof of Theorem \ref{T:gluing2}]
       We assume that $\Y$, $\Z$ and $\LL$ are
        topological stacks, and that $\Y$ and $\Z$ admit nice charts
        with respect to $\LL$. The pretopological case follows by taking
        $\mathbf{LF}$ to be the class of all epimorphisms
        (see also Example \ref{E:neat}.$\mathbf{4}$).

         Pick LF charts $p \: Y \to \Y$ and $q\: Z \to \Z$
         such that the (invariant) subspaces
         $L_1 \subseteq Y$ and $L_2 \subseteq Z$ corresponding to $\LL$
         are neat.
          We would like to glue $Y$ to $Z$ along these subspaces to build
         a chart for the would-be $\Y\vee_{\LL}\Z$.
         But the problem is that
         $L_1$ may not be homeomorphic to $L_2$. Thanks to Lemma \ref{L:gluing},
         we can overcome this problem
          as follows.

     For any two LF charts $L_1$ and $L_2$ for $\LL$, there exists
     a third chart $L$ which factors through $L_1$ and $L_2$ via
     local fibrations (e.g. take $L=L_1\x_{\LL}L_2$).
     By  Lemma \ref{L:gluing}.$\mathbf{ii}$, we may replace $L$
     by some
     open cover and assume that it extends to a chart for $\Y$.
     By refining the open cover, we may assume that $L$ extends to
     a chart for $\Z$ as well (here we used Lemma \ref{L:gluing}.$\mathbf{i}$).

     From now on, we will assume that $L_1$ and $L_2$ are homeomorphic.
     Let $[S \sst{} Y]$ and $[T\sst{}Z]$
     be the groupoid associated to the chart $Y \to \Y$and $Z\to\Z$.
     Denote the restriction
     of $S$ to $L_1$ by $[R_1 \sst{} L_1]$ and the restriction of
     $T$ to $L_2$ by $[R_2 \sst{} L_2]$. We have $R_1 \cong L_1\x_{\LL} L_1$
     and $R_2 \cong L_2\x_{\LL} L_2$, so the homeomorphism between
     $L_1$ and $L_2$ induces an isomorphism between the groupoids
     $[R_1 \sst{} L_1] \risom [R_2 \sst{} L_2]$. We can now glue
     the groupoids $[S \sst{} Y]$ and $[T \sst{} Z]$ along this isomorphism.
     The quotient stack of this new groupoid is the desired gluing
     $\Y\vee_{\LL}\Z$.
    \end{proof}

  There are two issues with the gluing of topological stacks: 1) It is not
  guaranteed that when we glue topological stacks along neat substacks
  we obtain topological stacks;
  2) The pretopological stack $\Y\vee_{\LL}\Z$ may not satisfy the universal
   property of push-out.

\vspace{0.1in}
  We discuss the first issue in the following remark.

   \begin{rem}{\label{R:gluing}}
       The source and target   maps of the groupoid constructed in the
       proof of the theorem are not expected to be LF in general.
       However, they are very close to be so. For instance,
       the source map is a map of the form $s \: X \to  Y\vee_L Z$
       where both restrictions $s|_Y \: s^{-1}(Y) \to Y$
       and  $s|_Z \: s^{-1}(Z) \to Z$ are LF. Although this in general
       does not imply that $s$ is LF, but it is quite likely that in
       specific examples one can prove by hand that this is
       the case. Especially, note
       that usually $L$ is a local retract in both $Y$ and $Z$.

       If we take $LF$ to be the class of local homeomorphisms, or locally
       cartesian maps (Example \ref{E:lf}.$\mathbf{3}$), then it {\em is}
       true that $s$ and  $t$ are again LF.
       So in this case the glued stack is again topological.
       Note also that in these cases every subspace is automatically
       neat (Example \ref{E:neat}).
   \end{rem}

 %


  Now, we turn to the universal property of $\Y\vee_{\LL} \Z$.
  It follows from the definition of $\Y\vee_{\LL} \Z$
  (see axioms $\mathbf{G1}$ and $\mathbf{G2}$) that there is a natural
  restriction functor \label{universal}

  $$\varrho\:\Hom_{\mathsf{St}}(\Y\vee_{\LL} \Z, \X) \to \Glue(\Y,\Z,\LL; \X),$$

  \vspace{1mm}
   \noindent where $\Glue(\Y,\Z,\LL; \X)$ is the groupoid of gluing data defined
   in page \ref{pushout}.
   We say $\Y\vee_{\LL} \Z$ satisfies the {\em universal property}
   if this functor is an equivalence of groupoids.

  We generalize Theorem \ref{T:gluing} as follows.

 \begin{thm}{\label{T:gluing3}}
   Let $\LL \to \Y$ and $\LL \to \Z$ be embeddings of pretopological
   stacks, and let $\Y\vee_{\LL} \Z$ be a gluing (as in Theorem
   \ref{T:gluing2}). Then:

   \begin{itemize}

          \item[$\mathbf{i.}$] The natural functor
  $$\varrho\:\Hom_{\mathsf{St}}(\Y\vee_{\LL} \Z, \X) \to \Glue(\Y,\Z,\LL; \X)$$
             is fully faithful.

          \item[$\mathbf{ii.}$] If there is a chart
             $p \: W \to \Y\vee_{\LL} \Z$ such that the corresponding embedding
             $i \: L \to Y$ is LTC ($L$ and $Y$ are invariant subspaces of $W$
             corresponding to $\LL$ and $\Y$),
             then $\varrho$ is an equivalence of groupoids for every topological
             stack $\X$.
       \end{itemize}
     \end{thm}

   \begin{proof}
      Pick a chart   $p \: W \to \Y\vee_{\LL} \Z$ and let
      $Y$, $Z$ and $L$ be the invariant subspaces of $W$
      corresponding to  $\Y$, $\Z$ and $\LL$, respectively.
      Using the description of maps coming out of a quotient stack
      (Section \ref{SS:maps}) the problem reduces to
      Theorem \ref{T:gluing}. Details left to the reader.
   \end{proof}

    \begin{cor}{\label{C:uniqueness}}
      If in part {\em ($\mathbf{ii}$)} of
      Theorem \ref{T:gluing3} the pretopological stack $\Y\vee_{\LL} \Z$
      is actually topological, then it is unique (up to an equivalence
      that is unique up to 2-isomorphism). This is automatically the
      case if we take $\mathbf{LF}$ to be the class of locally cartesian maps
      (Example \ref{E:lf}.$\mathbf{3}$) or local homeomorphisms.
    \end{cor}

   \begin{proof}
      The uniqueness follows from the universal property
      (Theorem \ref{T:gluing3}). The last statement follows from Remark
      \ref{R:gluing}.
   \end{proof}

   \begin{cor}{\label{C:gluing}}
       Let $\Y$ and $\Z$ be topological stacks, and let $\U \subseteq \Y$
       and $\V \subseteq \Z$ be open substacks. Suppose we have an
       equivalence $f \: \U \risom \V$. Then $\Y$ can be glued to
       $\Z$ along $f$ to form a topological stack.
       The glued stack satisfies the universal
       property of push-out.
    \end{cor}

  In fact, as long as we are gluing  along {\em open} substacks,
  any coherent gluing data can be used to glue topological stacks,
  and the resulting stack will again be topological
  (and will satisfy the relevant universal property).
  It should be clear how to
  formulate such a statement, but we will not do it here because it
  a bit messy and we will
  not be needing it later.

   \begin{thm}{\label{T:gluingDM}}
     Consider the diagram topological stacks
         $$\xymatrix@=12pt@M=10pt@M=10pt@=12pt{
                 \LL \ar@{^(->} [r]^i \ar@{^(->}[d]_j &    \Y  \\
                        \Z          &        }$$
    in which $i$ and $j$ are embeddings.
       \begin{itemize}

          \item[$\mathbf{i.}$] If $\Y$ and $\Z$ are
            weak Deligne-Mumford, then
            so is $\Y\vee_{\LL}\Z$. Furthermore, $\Y\vee_{\LL}\Z$
            satisfies the universal property
            of push-out (page \pageref{universal}).

          \item[$\mathbf{ii.}$] Assume $i$ and $j$ are locally-closed
            embeddings (that is, a closed embedding followed by an open
            embedding). If $\Y$ and $\Z$ are
           Deligne-Mumford, then
             so is $\Y\vee_{\LL}\Z$.

        \end{itemize}
     \end{thm}

   \begin{proof}[Proof of part {\em ($\mathbf{i}$)}]
   Let $\mathbf{LF}$ be the class of local homeomorphisms. Then
   every embedding $i \: L \to Y$ of topological spaces
   is LTC (Example \ref{E:LTC}.$\mathbf{6}$) and neat
   (Example \ref{E:neat}.$\mathbf{3}$).
   Furthermore, if we glue local homeomorphisms
   along subspaces we obtain a local homeomorphism (see Remark \ref{R:gluing}).
   So the chart $X \to \Y\vee_{\LL}\Z$ constructed in
   Theorem \ref{T:gluing2} is \'etale.

   \vspace{0.1in}
   \noindent{\em Proof of part} ($\mathbf{ii}$).
      Pick an \'{e}tale chart $X \to \Y\vee_{\LL}\Z$. Let $x \in X$
     be an arbitrary point. We have to find an open around $x$ that
     is invariant under $I_x$ and the action of $I_x$ on it is
     mild. Let $L\subseteq X$ be
     the subspace of $X$ corresponding to $\LL$. If $x$ is not in $L$
     the assertion is obvious. So assume $x$ is in $L$.
     Let $Y, Z \subseteq X$
     be the subspaces corresponding to $\Y$ and $\Z$. Since $\Y$ is
     Deligne-Mumford,
     there is a  subset
     $U \subseteq Y$   containing $x$ that is
     invariant under the action of $I_x$ and  such that the induced groupoid
     on $U$ is isomorphic to the action groupoid $[I_x \x U \sst{} U]$
     of a mild action of $I_x$ on $U$ (Proposition \ref{P:DMchart}).
     We can find $V \subseteq Z$ with the similar property. Since
     $L$ is locally closed in both $Y$ and $Z$, we may also assume that
     $U$ and $V$ are small enough so that  $U \cap L$ and $V \cap L$ are
     closed in $U$ and $V$, respectively. Let
     $U'=U\backslash\big(L\backslash(V\cap L) \big)$
     and   $V'=V\backslash\big(L\backslash(U\cap L) \big)$.
     These are open sets in $U$ and $V$,
     respectively, and we have $U'\cap V'=U\cap V$. Furthermore,
     both $U'$ and $V'$ are $I_x$ invariant. Set
     $W=U'\cup V' \subseteq X$. It is easy to see that $W$ is open in $X$,
     and the induced groupoid on $W$ is isomorphic to the action groupoid
     $[I_x \x W \sst{} W]$ of a mild action of $I_x$.
   \end{proof}

\section{Elementary homotopy theory}{\label{S:homotopy}}

   In this section we begin developing homotopy theory of topological stacks.
   Recall that our definition of a topological stack depends on a choice
   of a class $\mathbf{LF}$ of local fibrations (Section \ref{SS:fibrations}).
   We require that $\mathbf{LF}$ is so that CW inclusions of finite
   CW complexes are LTC. This is the case for all the classes
   considers in Example \ref{E:lf} except for ($\mathbf{4}$) and ($\mathbf{5}$).

    We use squiggly  arrows $\rsa$,
    instead of double arrows $\twomor{}$, to
    denote 2-isomorphisms between points.

    \begin{defn}{\label{D:pair}}
     By a {\bf pair} we mean a pair $(\X,\A)$ of topological stacks and a
     morphism (not necessarily an embedding) $i \: \A \to \X$ between them
     (we usualy drop $i$ in the notation). When $\A=*$ is a point, we call this
     a {\bf pointed stack}. A {\em morphism of pairs}
     $(\X,\A) \to (\Y,\B)$ consists of two morphisms $(f,f')$
     and a 2-isomorphism
     as in the following 2-cell:
        $$\xymatrix@=16pt@M=6pt{ \A \ar[r]^{f'}  \ar@{} @<-4pt> [r]| (0.4){}="a"
               \ar[d]_{i} \ar@{} @<4pt> [d]| (0.4){}="b"
                                   & \B \ar[d]^{j}     \ar @{=>}  "a";"b" \\
                                        \X \ar[r]_{f}        & \Y           }$$
     \noindent We usually abuse notation and denote such a morphism by $f$.
     Given two morphisms $f,g \: (\X,\A) \to (\Y,\B)$ of pairs, we define
     an {\bf identification} (or  {\em transformation}, {\em 2-isomorphism})
        from $f$ to $g$
     to be a pair of transformations $f\twomor{\varphi}g$ and
     $f'\twomor{\varphi'}g'$ making the following 2-cell commute:
     $$\xymatrix@=12pt@M=6pt@=30pt{     \A \ar[d]  \ar@/_/ @<-1pt> [r]_{g'} \ar@{} @<-5pt>    [r] | (0.5){}="b"
              \ar@/^/ @<1pt> [r]^{f'} \ar@{} @<5pt>    [r] | (0.5){}="a"        &    \B \ar[d]
                                          \ar @{=>}_{\varphi'}  "a";"b"  \\
               \X  \ar@/_/@<-1pt> [r]_{g}  \ar@{} @<-5pt>    [r]| (0.5){}="d"
                          \ar@/^/ @<1pt> [r]^{f} \ar@{} @<5pt>   [r] | (0.5){}="c"
                                      &      \Y  \ar @{=>}_{\varphi}  "c";"d"}$$
     \noindent We usually abuse notation and denote such a transformation by
     $\varphi$.
 \end{defn}

  For pointed stacks we use the (slightly different) notation
  $(\X,x)$, where $x \: * \to \X$ is the given point.

 A morphism of pointed stacks is a pair $(f,\phi) \: (\X,x) \to (\Y,y)$,
 where $f \: \X \to \Y$ is a morphism of stacks and $\phi \: y \rsa f(x)$
 is an identification.


\begin{defn}{\label{D:htpy}}
   Let $f,g \: (\X,\A) \to (\Y,\B)$ be maps of pairs as in the previous
   definition. A {\bf homotopy} from $f$ to $g$ is a triple
   $(H,\ep_0,\ep_1)$ as follows:

    \begin{itemize}

    \item  A map of pairs $H \: (I\times \X,I\times\A) \to (\Y, \B)$.

    \item  A pair of identifications $f \twomor{\ep_0} H_0$ and
           $H_1 \twomor{\ep_1} g$.  Here $H_0$ and $H_1$ stand for
           the morphisms of pairs $(\X,\A) \to (\Y,\B)$ obtained
           by restricting $H$ to $\{0\}\x\X$ and $\{1\}\x\X$, respectively.
    \end{itemize}
    We usually drop $\ep_0$ and $\ep_1$ in the notation.
    We denote the homotopy classes of maps of pairs $(\X,\A) \to (\Y,\B)$
    by $[(\X,\A), (\Y,\B)]$.
    If $x,y \: * \to \X$ are points in $\X$, a homotopy between $x$ and $y$
    is called a {\bf path} from $x$ to $y$.
\end{defn}

 \begin{rem}{\label{R:identification}}
   In view of the above definition, an identification between two maps
   of pairs can be regarded as a homotopy. For precisely, let
     $f,g \: (\X,\A) \to (\Y,\B)$ be maps of pairs and $\varphi \: f\twomor{} g$
     an identification. Then we have a homotopy $(H, \ep_0,\ep_1)$
     from $f$ to $g$ where $H$ is
       $$\xymatrix@C=-1pt@M=10pt{
          (I\times \X,I\times\A) \ar[rr]^(0.6){pr}   && (\X,\A) \ar[rr]^f  && (\Y,\B)
                        }$$
    and $\ep_0=\id$ and $\ep_1=\varphi$.
 \end{rem}

When the pair $(\X,\A)$ is reasonably well-behaved, we can compose
homotopies. More precisely, we have the following result.

  \begin{lem}{\label{L:composehtpy}}
     Let  $(\X,\A)$ be  a pair. Assume $\X$ and $\A$ admit
     charts $X \to \X$ and $A \to \A$ such that the $t=0$ inclusions
     $X \hra I\x X$ and $A \hra I\x A$ are LTC (e.g. $X$ and $A$
     both CW complexes). Then, for any pair $(\Y,\B)$ of topological stacks,
     homotopies between  maps from $(\X,\A)$ to $(\Y,\B)$ can be composed in a
     natural way. In other words, given maps of pairs
     $f_1,f_2,f_3 \: (\X,\A) \to (\Y,\B)$ , and homotopies $H_{12}$ from $f_1$
     to $f_2$ and $H_{23}$ from $f_2$ to $f_3$, we can compose
     $H_{12}$ and $H_{23}$ in a natural way
     to obtain a homotopy $H_{13}$ from $f_1$ to $f_3$.
  \end{lem}

 \begin{proof}
   Follows from Theorem \ref{T:gluing3}.
 \end{proof}

\begin{defn}
   Let $(\X,x)$ be a pointed topological stack. We define
   $\pi_n(\X,x)=[(S^n,\bullet),(\X,x)]$ for $n \geq 0$.
\end{defn}

As for topological spaces, there is a natural group structure on
$\pi_n(\X,x)$ for every $n \geq 1$. The multiplication is defined
as follows.

Let $f,g \: (S^n,\bullet) \to (\X,x)$ be pointed maps. From the
definition of a map of pairs (Definition \ref{D:pair}), there is a
natural identification $\al: f(\bullet) \rsa g(\bullet)$. So, by
Proposition \ref{P:gluing1}, there is a map $f\vee g \: S^n\vee
S^n \to \X$ whose restrictions to each  copy of $S^n$ is naturally
`identified' with  $f$ and $g$. The restrictions of these
identifications to the base points gives the following commutative
diagram of identifications:

  $$\xymatrix@R=14pt@C=-3pt{   & x \ar@{~}[ld] \ar@{~}[rd] &     \\
                  f(\bullet) \ar@{~}[rr]^{\al} \ar@{~}[rd] &
                                            & g(\bullet) \ar@{~}[dl] \\
                                       & (f\vee g)(\bullet) &           }$$

\noindent In particular, we have an identification   $x \rsa
(f\vee g)(\bullet)$, which makes $f\vee g \: (S^n\vee S^n,\bullet)
\to (\X,x)$ into a pointed map. It can be checked, using
Proposition \ref{P:gluing1},
 that this construction respects homotopies
of pointed maps (Definition \ref{D:htpy}).
 Precomposing $f\vee g$ with the $\frac{1}{2}+\frac{1}{2}$
map  $S^n \to S^n\vee S^n$  gives rise to a pointed map
$(S^n,\bullet) \to (\X,x)$, which we define to be the product $fg$
of $f$ and $g$ in $\pi_n(\X,x)$.

Similarly, given three maps $f,g,h \: (S^n,\bullet) \to (\X,x)$,
we can construct a map $f\vee g\vee h \: (S^n\vee S^n\vee S^n,
\bullet) \to (\X,x)$, unique up to a unique 2-isomorphism, whose
restriction to each copy of $S^n$ is naturally identified with the
corresponding map $f$, $g$ or $h$.
%

Pre-composing this map with the
$\frac{1}{4}+\frac{1}{4}+\frac{1}{2}$ map $S^n \to S^n\vee S^n\vee
S^n$ gives a pointed map $(S^n,\bullet) \to (\X,x)$ that is
naturally identified with $(fg)h$. Similarly, pre-composing with
the    $\frac{1}{2}+\frac{1}{4}+\frac{1}{4}$ map $S^n \to S^n\vee
S^n\vee S^n$   gives $f(gh)$. Since the
$\frac{1}{4}+\frac{1}{4}+\frac{1}{2}$ and
$\frac{1}{2}+\frac{1}{4}+\frac{1}{4}$ maps
 are homotopic, they remain homotopic after
composing with $f\vee g\vee h$, so we have $(fg)h=f(gh) \in
\pi_n(\X,x)$. This verifies associativity. The other axioms of a
group are verified in a similar fashion.

 The homotopy groups of topological stacks as defined above are
 functorial with respect to
 pointed maps of topological stacks. That is, given a pointed map
 $(f,\phi) \: (\X,x) \to (\Y,y)$, we get an induced map
 on the homotopy groups
        $$\pi_n(f,\phi) \: \pi_n(\X,x) \to \pi_n(\Y,y).$$
 \noindent Note that the above maps do depend on $\phi$.
 We sometime denote this map by  $\pi_n(f)$, $(f,\varphi)_*$, or $f_*$,
 if there is no fear of confusion.

 Homotopy groups of topological stacks  have all the
 natural properties of the homotopy groups of topological spaces. For instance,
 $\pi_n(\X,x)$ is abelian for $n \geq 2$,
 there is an action of $\pi_1(\X,x)$ on $\pi_n(\X,x)$, there is a Whitehead
 product which satisfies the Jacobi identity and so on.
 To see these, notice that
 in the classical case all these results reflect certain structures
 on the spheres $S^n$ (e.g., $\pi_n(\X,x)$ is an abelian group for
 $n \geq 2$ because
 $S^n$ is a homotopy abelian cogroup for $n \geq 2$). Therefore,  arguments,
 similar to the one use above to prove the associativity of the multiplication
 on $\pi_n(\X,x)$, can be used to prove other properties of
 homotopy groups of stacks.

As pointed out in the previous paragraph, the homotopy groups of a
topological stack behave pretty much like the homotopy groups of
topological spaces. There is, however, more structure on the
homotopy groups of a  topological stack. We discuss the case of
fundamental group in detail. The higher cases have been briefly
touched upon in Section \ref{SS:gerbes}.

For a pointed topological stack $(\X,x)$ there is natural group
homomorphism
        $$\ox \: I_x \to \pX$$
(see Remark \ref{R:identification}). For this reason, and also to
be consistent with the case of higher homotopy groups (see Section
\ref{SS:gerbes}), we sometimes refer to $I_x$ as the {\bf inertial
fundamental group} of $(\X,x)$. Recall that we have a natural
isomorphism
       $$I_x\cong\{\al \ | \ \al \: x \rsa x \ \text{self-transformation}\}.$$

The maps $\ox$  are functorial with respect to pointed maps. That
is, given a pointed map $(f,\phi) \: (\X,x) \to (\Y,y)$, we get an
induced map $(f,\phi)_* \: I_x \to I_y$, and the following square
commutes:
  $$\xymatrix@=12pt@M=10pt{I_x
           \ar[r]^(0.4){\ox} \ar[d]_{(f,\phi)_*} & \pX \ar[d]^{\pi_1(f,\phi)}   \\
           I_y   \ar[r]_(0.4){\oy}        &   \pY        }$$

    We will use inertial fundamental groups in the study of covering
    spaces of topological stacks.
\section{Covering theory for topological stacks}{\label{S:covering}}

   We develop a Galois theory of
   covering spaces for a topological stack. In particular,  we show that,
   under appropriate local assumptions,
   there is natural correspondence between
   subgroups of the fundamental group $\pX$ of a pointed topological
   stack $(\X,x)$ and the (isomorphism classes of) pointed covering spaces
   of $(\X,x)$.
   We then explain the role of inertial  fundamental groups as bookkeeping
   devices which keep track of the stacky structure of  covering spaces.
   Similar results have been proven in the algebraic setting in \cite{Noohi1},
   and most of the proofs can be easily adopted to the topological setting.

   We fix a class $\mathbf{LF}$ of local fibrations as in
   Section \ref{S:homotopy}.

  \subsection{Connectedness conditions}{\label{SS:connectivity}}

  In this subsection we generalize several notions of
  connectedness topological stacks.

  In what follows, whenever we talk about connected components,
  we can assume our class $\mathbf{LF}$ of local fibrations is arbitrary.
  So our discussion will be valid for pretopological stacks too.
  However, when we talk about {\em path} components, we have
  to assume $\mathbf{LF}$ is such that the inclusion   $\{0\} \hra I$
   is LTC. This is the case
  for every example in Example \ref{E:lf} except for ($\mathbf{4}$) and
  ($\mathbf{5}$).

  \begin{defn}{\label{D:pc}}
      Let $\X$ be a pretopological stack. We say $\X$ is {\bf
      connected} if it has no proper open-closed substack.
       We say $\X$ is
      {\bf path connected}, if for every two points $x$ and $y$ in $\X$,
      there is a path from $x$ to $y$ (Definition \ref{D:htpy}).
   \end{defn}

   \begin{defn}{\label{D:lpc}}
      Let $\X$ be a pretopological stack. We say $\X$ is {\bf
      locally connected} (resp.,
      {\bf locally path connected}, {\bf semilocally 1-connected}),
     if there is a chart $X \to \X$ such that  $X$ is so.
   \end{defn}

  This definition agrees with the usual definition when $\X$ is a topological
  space. This is because of the following lemma.

   \begin{lem}{\label{L:lpc}}
      Let $f \: Y \to X$ be an epimorphism of topological spaces.
      Assume $Y$ is locally connected (resp., locally path connected,
      semilocally 1-connected). Then so is $X$.
    \end{lem}

   \begin{proof}[Proof of locally connected.] Being locally connected
     is a local condition, so we may assume $f$ has a section
     $s \: X \to Y$. Recall that $X$ is
     locally connected if an only if for every open subset $U \subseteq X$
     every connected component of $U$ is open (hence also closed).
     So it is enough to show that every connected component of $X$ is
     open; same argument will work for any open subset $U$.
     Let $Y=\coprod Y_i$ be the decomposition of $Y$ into connected components.
     So every $Y_i$ is open-closed. For each $i$,
     $X_i :=s^{-1}(Y_i)$ is open-closed (possibly empty).
     It is also connected because, if $X_i$ is not empty, $f$ surjects $Y_i$
     onto $X_i$ (note that $f(Y_i)\subseteq X_i$ since $Y_i$ is connected and
     $X_i$ is open-closed).
      Also note that, if $i\neq j$,
     then $X_i$ and $X_j$ are disjoint. This shows that the connected
     components of $X$ are exactly those $X_i$ that are non-empty.
     So all the connected components of $X$ are open.

     \vspace{0.1in}
     \noindent{\em Proof of locally path connected.} Similar to the previous
     case.

    \vspace{0.1in}
    \noindent{\em Proof of semilocally 1-connected.}
     Let $x \in X$ be an
     arbitrary point, and assume $U$ is an open neighborhood of $X$
     over which $f$ has a section $s \: U \to Y$.
     Let $V$ be an open neighborhood of $y=s(x)$ such that
     $\pi_1(V,y) \to \pi_1(Y,y)$ is trivial. It is easily seen that
     $s^{-1}(V) \subseteq X$ is an open neighborhood of $x$ with the similar
     property.
   \end{proof}

    It is easy to see that, $\X$ is locally (path) connected
    if and only if for every map $f \: W \to \X$ from a topological space $W$
    to $\X$, and every point $w \in W$,
    there is an open neighborhood $U\subseteq W$ of $w$ in $W$ such that
    $f|_U$ factors through a map $V \to \X$, where $V$ is a locally (path)
    connected topological space.

 Similarly, $\X$ is semi-locally 1-connected
 if and only if for every map $f \: W \to \X$ from a topological space $W$
 to $\X$, and every point $w \in W$, there is an open neighborhood
 $U\subseteq W$
 of $w$ in $W$ such that
 induced map on fundamental groups
 $\pi_1(U,w) \to \pi_1(\X,f(w))$ is the zero map.

    \begin{lem}{\label{L:lpcmod}}\end{lem}
   {\em  \begin{itemize}

           \item[$\mathbf{i.}$] Let $\X$ be a pretopological stack.
              Then, $\X$ is connected if and only if $\Xm$ is so.
              If $\X$ is path connected, then so is $\Xm$.

           \item[$\mathbf{ii.}$] Let $\X$ be a pretopological stack.
              If $\X$ is locally connected
              (resp., locally path connected), then so is $\Xm$.

        \end{itemize}
    }

    \begin{proof}[Proof of part {\em ($\mathbf{i}$)}] Obvious.

 \vspace{0.1in}
    \noindent {\em Proof of part} ($\mathbf{ii}$).
    We only prove the locally connected case. The  locally
    path connected case is proved analogously.

     Let $X \to \X$ be a locally connected
       chart for $\X$. Then $f \: X \to \Xm$ is
       a quotient map (Example \ref{E:quotient}).
       So we have to show that, if $f \: X \to Y$ is a quotient map of
       topological spaces
       and $X$ is locally  connected, then so is $Y$. We show that,
       for every open $U \subseteq Y$, the connected components of $U$
       are open. Since the restriction of $f$ over $U$ is again a quotient map,
       we are reduced to the case $U=Y$. Let $Y_0$ be a connected component
       of $Y$. We have to show that $f^{-1}(Y_0)$ is open in $X$. But
        $f^{-1}(Y_0)$ is a union of connected components of $X$. The
         claim follows.
    \end{proof}



  It is possible to have a topological stack that is not locally
  connected and not path connected but whose coarse
  moduli space is path connected and locally path connected (see
  Example \ref{E:twopoint}).

    \begin{defn}{\label{D:components}}
      Let $\X$ be a topological stack (see the assumption
      on $\mathbf{LF}$ at the beginning of this subsection),
      and let $x$ be a point in $\X$.
      The ({\bf path}) {\bf component} of $x$ is defined to be the union
      (Section \ref{SS:operations})
      of all (path) connected embedded substacks of $\X$   containing
      $x$.
    \end{defn}

    It is an easy exercise that each  component is  connected.
    Similarly, each path component is path connected
    (Proposition \ref{P:gluing1} is needed for this).
    The connected  components of $\X$ correspond exactly to connected components
    of $\Xm$ (in the obvious way). The path components of $\X$ give rise
    to a partitioning of $\Xm$ that is  finer than (or equal to)
    the partitioning of $\Xm$
    by its path components.

  \begin{lem}{\label{L:components}}
     Let $\X$ be a topological stack and $p \: X \to \X$ a chart for it.
     Then, for every (path) component $\X_0$ of $\X$, the corresponding
     invariant subset $X_0 \subseteq X$ is a union of (path) components
     of $X$. In particular, if $\X$ is locally (path) connected,
     then every (path) component of $\X$ is an open-closed substack.
  \end{lem}

  \begin{proof}
     Let $A \subseteq X$ be a (path) component of $X$ that is not
     a subset of $X$ but $X_0\cap A \neq \emptyset$. Then it is easy
     to see that $\OO(X_0\cup A)$ is an invariant subspace of $X$
     whose corresponding embedded substack in $\X$ is (path) connected, and
     strictly bigger than $\X_0$. This is a contradiction.
  \end{proof}

   \begin{cor}{\label{C:pc1}}
      Let $\X$ be a connected, locally path connected topological stack.
      Then $\X$ is path connected.
   \end{cor}

  The following lemmas will be used later.

  \begin{lem}{\label{L:pc2}}
    Let $f \: \Y \to \X$ be a surjective local homeomorphism of pretopological
    stacks.
    Then $\X$ is locally (path) connected if and only if $\Y$ is so.
   \end{lem}

   \begin{proof}
      If $\X$ is locally (path) connected, the pull back of
      locally (path) connected chart for $\X$ is a locally (path) connected
      chart for $\Y$. Conversely, if $Y \to \Y$ is a locally (path) connected
      chart for $\Y$, the composition $Y \to \Y \to \X$ will be a locally
      (path) connected
      chart for $\X$. (Note that every surjective local homeomorphism is
      an epimorphism.)
  \end{proof}

   \begin{lem}{\label{L:lpcDM}}
     Let $\X$ be a locally (path) connected weak Deligne-Mumford topological
     stack, and let $p \: X \to \X$ be an \'{e}tale chart for it.
     Then $X$ is locally (path) connected. Same statement is true
     with semilocally 1-connected.
   \end{lem}

   \begin{proof}
      Let $q \: X' \to \X$ be a chart such that $X'$ is locally (path)
      connected. Then $Y:=X\x_{\X}X'$ is also locally (path) connected,
      because  the projection map $X\x_{\X}X' \to X'$ is a local homeomorphism
      and $X'$ is locally connected. The   map $Y \to X$
      is now an epimorphism from a locally (path) connected space
      to $X$. This implies that $X$ is also locally (path) connected
      (Lemma \ref{L:lpc}).
   \end{proof}

  \subsection{Galois theory of covering spaces of a topological stack }
  {\label{SS:Galois}}

   We begin with recalling
   the definition of a covering map of pretopological stacks.

   \begin{defn}{\label{D:covering}}
      A map $\Y \to \X$ of  pretopological stacks
     is a {\bf covering map}, if it is representable, and for every
     topological space $W$ and every map $W \to \X$, the base extension
    $W\x_{\X}\Y \to W$ is a covering map of topological spaces.
    \end{defn}

    \begin{prop}{\label{P:diagonal}}
       Let $p \: \Y \to \X$ be a covering map of pretopological stacks.
       Then the diagonal
       map $\De \: \Y \to \Y\x_{\X} \Y$ is an open-closed embedding.
    \end{prop}

    \begin{proof}
       Let $X \to \X$ be a chart for $\X$, and let $Y \to \Y$
       be the pull back chart for $\Y$. The base extension $q \: Y \to X$ of $p$
       is a covering map. We have a 2-cartesian diagram
           $$\xymatrix@=12pt@M=10pt{ Y \ar[r] \ar[d]_{\De} & \Y \ar[d]^{\De}  \\
                       Y\x_X Y \ar[r]        &   \Y\x_{\X}\Y        }$$
      Since $q \: Y \to X$ is a covering map of topological spaces,
      the left vertical
      map is an open-closed embedding. So we have shown that
      the base extension of $\De \: \Y \to \Y\x_{\X} \Y$ along the
      epimorphism $Y\x_X Y \to \Y\x_{\X}\Y$ is an open-closed embedding.
      So   $\De \: \Y \to \Y\x_{\X} \Y$ itself is also an open-closed embedding.
    \end{proof}

 Using the above proposition, we can prove the uniqueness of liftings
 for covering maps.

  \begin{lem}{\label{L:uniqueness}}
    Let $p \: \Y \to \X$ be a covering map of pretopological stacks.
    Let $\T$ be a connected pretopological stack,
    and let $t$ be a point in $\T$.
    Let $f, g \: \T \to \Y$ be maps, and let
    $\Phi \: p\circ f \twomor{} p\circ g$
    be a 2-isomorphism. Assume also  $f(t)$ and $g(t)$
    are equivalent points in $\Y$, and fix a 2-isomorphism
    $\psi  \: f(t) \rsa g(t)$. Then, there is a unique 2-isomorphism
    $\Psi \: f \twomor{} g$ such that its restriction to $t$ is equal to
    $\psi$ and $p\circ\Psi=\Phi$.
  \end{lem}

   \begin{proof}
     Consider the map $(f,g, \Phi) \: \T \to \Y\x_{\X}\Y$.
     Form the 2-fiber product
       $$\xymatrix@=16pt@M=10pt{    \Z \ar[r] \ar[d]_d &  \Y\ar[d]^{\De}  \\
                      \T \ar[r]_(0.4){(f,g,\Phi)} \ar @{..>} [ur]
                                                &    \Y\x_{\X}\Y   }$$
    The assertion is equivalent to saying that,
    $(f,g, \Phi) \: \T \to \Y\x_{\X}\Y$ can be lifted to $\Y$ as a pointed
    map. Note that this is a 2-cartesian
    diagram of pointed stacks,
    where we take $f(t)$ for the base point of $\Y$, and
    $\big(f(t),f(t),id\big)$
    for the base point of $\Y\x_{\X}\Y$). The base point of $\Z$ is obtained
    from $\psi$.

    The map $d \: \Z \to \T$ is a base extension of $\De$, so, by Proposition
    \ref{P:diagonal}, it is an open-closed embedding. Therefore, $d$ is
    an equivalence of pointed stacks. Composing an inverse of $d$ with
    the projection $\Z \to \Y$ gives us the desired lift.
   \end{proof}



  \begin{defn}{\label{D:fiber}}
      Let $f \: \Y \to \X$ be a map of stacks, and $x \in \X$ a point in
      $\X$. We define the {\bf fiber} $\Y_x$ of $\Y$ over $x$ to be the following
      groupoid:
\vspace{1mm} {\small \begin{itemize}
      \item $\Ob(\Y_x)= \left\{\begin{array}{rcl} (y, \phi) & \vert & y \: * \to \Y,
                                                     \ \text{a point in $\Y$} \\
                      & & \phi \: x \rsa f(y),
                  \ \text{an identification} \end{array} \right\}$
\vspace{3mm}

      \item $\Mor_{\Y_x}\big((y,\phi),(y',\phi')\big)=
           \left\{  \beta \: y \rsa y' \ \vert \
                    \phi  \cdot  f(\beta)  =\phi' \right\}$
    \end{itemize}
}
  \end{defn}

In other words, $\Y_x$ is the groupoid of points in the  stack
$*\x_{\X}\Y$. When $f$ is representable, this groupoid is
equivalent to a set.

  \vspace{0.1in}

 Now, we introduce the Galois category $\mathsf{Cov}_{\X}$
 of a path connected topological stack $\X$. 
 The Galois
  category $\mathsf{Cov}_{\X}$ of $\X$ is defined as follows: \label{Cov}
\vspace{1mm} {\small
\begin{itemize}
\item  $\Ob (\mathsf{Cov}_{\X}) = \left\{ \begin{array}{rcl}
         (\Y, f)   & \vert & \text{ $\Y$, \ a topological stack;} \\
            & &  f \: \Y \to \X, \   \text{a covering space.}
                                      \end{array} \right\}$
\vspace{3mm}

\item  $\Mor_{\mathsf{Cov}_{\X}}\big( (\Y,f), (\Z,g)
\big)=\left\{\begin{array}{rcl}
           (a,\Phi) & \vert & a \: \Y \to \Z \  \ \text{a morphism;} \\
        && \Phi \: f \Ra g \circ a, \   \text{a 2-morphism.}
                                              \end{array}\right\}_{/_\sim}$
\end{itemize}
}

\noindent Here, $\sim$ is defined by \vspace{-1mm}

$$(a, \Phi) \sim (b, \Psi) \  \text{if}\   \exists \Gamma  \: a \Ra b \
                         \text{such that} \  \Phi  \cdot g(\Gamma)  =\Psi.$$

\vspace{.1in}

\noindent  The fiber functor $F_x \: \mathsf{Cov}_{\X} \to
(\pX-\textsf{Set})$
           is defined by $F_x(\Y):=\pi_0(\Y_x)$, where $\pi_0(\Y_x)$
           stands for the set of isomorphism classes of $\Y_x$
           (Definition
           \ref{D:fiber}). Note that, since covering maps are representable,
           $\Y_x$ is equivalent to a set. Also note that, if
           we have a pointed map $(\Y,y) \to (\X,x)$, then $F_x(\Y)$ is
           naturally pointed.

 The  action of $\pX$ on $F_x(\Y)$ is defined as follows.
 Let $\gamma \: I \to \X$
be a map representing a loop at $x$ (in particular, we are given
identifications  $\ep_0 \: x \rsa \gamma(0)$ and $\ep_1\:
\gamma(1) \rsa x$).
 Since $\Y \to \X$ is a covering map,
 the base extension $E=I\x_{\X}\Y$ is isomorphic
to a disjoint union of copies $I$. There are natural bijections
$\varepsilon_0 \: E_0 \risom F_x(\Y)$ and $\varepsilon_1 \:
F_x(\Y)  \risom E_1$ corresponding to $\ep_0$ and $\ep_1$, where
$E_0$ and $E_1$ are fibers of $E$ over $0$ and $1$, respectively.
Since $E_0$ and $E_1$ are canonically identified, we get an action
$\varepsilon_0\circ\varepsilon_1 \: F_x(\Y) \to F_x(\Y)$. The same
base extension trick, this time applied to $I\x I$, shows that
this action is independent of the pointed homotopy class of
$\gamma$. It is also easy to check that this action respects
composition of loops, hence induces a group action of
$\pi_1(\X,x)$ on $F_x(\Y)$. This action can indeed be jazzed up to
give a local system on the fundamental groupoid $\Pi_1(\X)$ with
fibers   $F_x(\Y)$, $x\in \X$.

 An interesting special case is when $\gamma$ comes from $I_x$.
 More precisely, let $\gamma \in I_x$ be an element in the inertia group
 at $x$, and consider the loop $\ox(\gamma) \: I \to \X$ that is given by
 the constant map
 (i.e. the map that
 factors through $x$ on the nose) together with the identifications
 $\ep_0=id \: x \rsa x$ and $\ep_1=\gamma \: x \rsa x$ (Definition
 \ref{D:htpy}).
 It is easy to check that the  action of $\gamma$
  on $F_x(\Y)$ sends (the class of) $(y,\phi)$ in  $F_x$
 to (the class of) $(y,\phi\cdot\gamma)$, where $\cdot$ stands
 for composition of identifications.

 \begin{lem}{\label{L:lifting}}
    Let $p\: (\Y,y) \to (\X,x)$ be a pointed covering map of topological stacks.
    Let
    $\al \: (S^1,\bullet) \to (\X,x)$ be a loop. Then the following are
    equivalent:
     \begin{itemize}

    \item[$\mathbf{i.}$] The action of $\al$ on $F_x(\Y)$ leaves $y$ invariant.

    \item[$\mathbf{ii.}$] There is a (necessarily unique) lift
    $\tilde{\al} \: (S,*) \to (\Y,y)$.

    \item[$\mathbf{iii.}$] The class of $\al$ in $\pX$ is in $p_*(\pY)$.

 \end{itemize}
  \end{lem}

   \begin{proof}[Proof of $(\mathbf{i}) \Rightarrow (\mathbf{ii})$]
     Notations being as in  the preceding   paragraphs, let $y_0 \in E_0$
     and $y_1 \in E_1$ be the points corresponding to $y$ (via $\epsilon_0$
     and $\epsilon_1$, respectively). Since $y$ is invariant under the action
     of $\pi_1(\X,x)$, $y_0$ and $y_1$ lie in the same ``layer'' in $E$,
     that is, there is a section $s\: I \to E$ whose end points are $y_0$
     and $y_1$. Composing $s$ with the natural map $E \to \Y$ we obtain
     a path whose end points are canonically identified with $y$.
     By Proposition \ref{P:gluing1}, this gives a loop
     $\tilde{\al} \: (S,\bullet) \to (\Y,y)$, which is the desired lift.

     \vspace{0.1in}

    \noindent{\em Proof of} $(\mathbf{ii}) \Rightarrow (\mathbf{iii})$.
      Obvious.

    \vspace{0.1in}

    \noindent {\em Proof of} $(\mathbf{iii}) \Rightarrow (\mathbf{i})$.
     Since the action is independent of the choice the (pointed) homotopy
     class of the loop, we may assume that $\al$ lifts to $(\Y,y)$. Now
     run the proof of $(\mathbf{i}) \Rightarrow (\mathbf{ii})$ backwards.
   \end{proof}

  \begin{prop}{\label{P:injective}}
     Let $p\: (\Y,y) \to (\X,x)$ be a pointed covering map. Then
     the induced map $\pY \to \pX$ is injective.
  \end{prop}

  \begin{proof} Let $\al,\be \: (S,\bullet) \to (\Y,y)$ be loops
   whose images in $\pX$ are equal. Then, there exists a pointed
   homotopy $H \: (I\x S^1, I\x\{\bullet\}) \to (\X,x)$
   between $p\circ\al$ and $p\circ\be$. Let
   $A=(\{0\}\x S^1)\cup  (I\x\{\bullet\})\cup (\{1\}\x S^1) \subset I\x S^1$.
   By Proposition \ref{P:gluing1} (or Theorem \ref{T:gluing}),
   there is a natural map $f \: A \to \Y$ whose restrictions to
   $\{0\}\x S^1$,   $\{1\}\x S^1$ and  $I\x\{\bullet\}$ are naturally
   identified with $\al$, $\be$, and the constant map.
   By  Theorem \ref{T:gluing},
   the composition $p\circ f \: A \to \X$ is naturally identified with
   $H|_A$ (to see this, note that such an identification exists on each
   of the three pieces $A$ is made of, and
   then use Theorem \ref{T:gluing}
   to show that they glue to an identification defined on $A$).

     Set $Z=(I\x S^1)\x_{\X}\Y$. Consider the following diagram:

          $$\xymatrix@=20pt@M=6pt{ Z
           \ar[r]^g \ar[d]^(0.45)q & \Y \ar[d]^(0.45)p  \\
             I\x S^1 \ar[r]_(0.55)H \ar@{..>} [ru]_F
                            \ar@{..>} @/^/ [u]^(0.45)s     &   \X        }$$
     Note that $Z$ is a topological space
     and $q$ is a covering map. The map $f$ gives rise to a partial section
      $\bar{f} \: A \to Z$ of $q$ defined on $A$. Since $q$ is a covering map,
      $\bar{f}$ extends to a section $s \: I\x S^1 \to Z$. The composite
      $F=g\circ s$ is the desired homotopy between $\al $ and $\be$.
  \end{proof}


  \begin{lem}{\label{L:mainlemma}}
   Let $f \: \Y \to \X$ be a map of topological stacks. Let $y$ be a  point
   of $\Y$ and  $x=f(y)$ its image in $\X$. Consider the following
   commutative diagram, obtained from the functoriality of the inertial
   fundamental groups:

     $$\xymatrix@=12pt@M=10pt{ I_y \ar[r]^(0.4){\omega_y} \ar[d]_{f_*}
                                             & \pY \ar[d]^{\pi_1(f)}\\
                                    I_x \ar[r]_(0.4){\omega_x} & \pX } $$

  \noindent If $f$ is a covering map, then this diagram is cartesian.

 \end{lem}

 \begin{proof}
    Let $\gamma$ be in $I_x$. We want to show that, if $\omega_x(\gamma)$
    is in the image of $\pi_1(f)$, then there
    exists a unique $\alpha \in I_y$ such that $f_*(\al)=\gamma$.
    The uniqueness is obvious, because
    $f_* \: I_y \to I_x$ is injective (use
    Proposition \ref{P:rep}.$\mathbf{ii}$).
    On the other hand, $\omega_x(\gamma)$ being in the image of
    $\pi_1(f)$ exactly means that, under the action of $\gamma$ on $F_x(\Y)$,
    the point $(y,id)$ remains invariant (Lemma \ref{L:lifting}).
    That means, $(y,id) \sim (y, \gamma)$, where $\sim$ means
    an isomorphism in the groupoid $\Y_x$
    (see Definition \ref{D:fiber}), which in this case is indeed
    equivalent to a set. Therefore, there exists
    $\beta \in I_y$ such that
    $f(\beta)=\gamma$.
  \end{proof}

   \begin{lem}{\label{L:isom}}
      Let $p \: (\Y,y) \to (\X,x)$ be a pointed covering map of
       path connected
      topological stacks. If the induced map $\pY \to \pX$ is
      an isomorphism, then $p$ is an equivalence.
    \end{lem}

    \begin{proof}
      It is enough to show that the fiber $F_x(\Y)$ is the set with one
        element, that element being the class of  $(y,\phi)$, where
       $\phi \: x \rsa p(y)$ is part of the data of a pointed map.
       (More precisely, if this is proven, then path connectedness implies that
       the fiber over every point $x' \in \X$ is a set with one element.
       Therefore,
       for any map $X \to \X$
      from a topological space $X$ to $\X$, the base extension $Y \to X$
      of $p$ will be a covering map of degree one, hence an isomorphism.
      In particular, the base extension of $p$ along a chart $X \to \X$
      is an isomorphism. Hence $p$
      must be an equivalence by Lemma \ref{L:epimorphism}.)
      Assume $(y',\phi')$ is
      another point in $F_x(\Y)$. Since $\Y$ is path connected
      (Corollary \ref{C:pc1} and Lemma \ref{L:pc2}), we can find a
      path $\gamma$ connecting $y$ and $y'$. The image of $\gamma$ in $\X$
      is a map $I \to \X$ such that the images of the end points are
      canonically equivalent to $x$. So,
      by Proposition \ref{P:quotient1}, we obtain a loop
      $\al \: (S^1,\bullet) \to (\X,x)$.
      By hypothesis, $\al$ is in $p_*\big(\pY\big)$.
      By Lemma \ref{L:lifting}, we can lift $\al$ to a loop
      $\tilde{\al} \: (S^1,\bullet) \to (\Y,y)$ at $y$. This loop
      gives rise to a map
      $\gamma' \: I \to \Y$, sending both $0$ and $1$ to (points equivalent to)
      $y$. Now, we have
      two maps $\gamma, \gamma' \: I \to \Y$, both lifting $\al$, and both
      sending $0$ to (points equivalent to) $y$. So, by Uniqueness Lemma
      \ref{L:uniqueness}, $\gamma$ and  $\gamma'$ are equivalent.
      Therefore $\gamma(1)$
      and $\gamma'(1)$ are also equivalent. But the first one is (equivalent to)
      $y$ and the second one is (equivalent to) $y'$. Hence, we obtain an
      identification $\rho \: y \rsa y'$.

       We are not done yet, because existence of $\rho$ is not enough
       to guarantee that $(y,\phi)\sim(y',\phi')$. We have to verify
       that $\phi \cdot p(\rho)=\phi'$. In fact, this is not necessarily true.
       However,  we can adjust $\rho$ as follows. Set
       $\be=\phi^{-1}\cdot\phi'\cdot p(\rho)^{-1} \in I_{p(y)}$. By Lemma
       \ref{L:mainlemma}, there exists $\tilde{\be} \in I_x$ such that
       $f(\tilde{\be})=\be$. Now, if we replace $\rho$ by
       $\sigma=\tilde{\be}\cdot\rho$
       we have our desired identification $\sigma \: y \rsa y'$ (i.e
       $\phi \cdot p(\sigma)=\phi'$).
    \end{proof}

 \begin{prop}[{\bf The general lifting lemma}]{\label{P:lifting}}
   Let $(\Y,y)$ and $(\X,x)$ be pointed connected locally path connected
   topological stacks, and let $p \: (\Y,y) \to (\X,x)$ be a pointed covering
   map. Suppose $(\T,t)$ is a pointed connected locally path topological stack,
   and let $f \: (\T,t) \to (\X,x)$ be a pointed
   map. Then, $f$ can be lifted to a pointed
   map $\tilde{f} \: (\T,t) \to (\Y,y)$ if and only if
   $f_*\big(\pi_1(\T,t)\big) \subseteq p_*\big(\pY\big)$.
   Furthermore, if such a lifting exists,
   it is unique.
   (Of course, the terms `lifting' and `uniquness' are to be interpreted as
   `up to unique
   2-isomorphism'.)
 \end{prop}

  \begin{proof}
      One implication it trivial. Assume now that
      $f_*(\pi_1(\T,t)) \subseteq p_*(\pY)$.
      Let $(\Z,z)$ be the connected
      component of the 2-fiber product $\T\x_{\X}\Y$ containing the  canonical
      base point $z$. The map $q\: (\Z,z) \to (\T,t)$ is a pointed
      covering map. It follows from the hypothesis that
      the induced map $\pi_1(\Z,z) \to \pi_1(\T,t)$ is surjective (Lemma
      \ref{L:lifting}),
      hence an isomorphism (Proposition \ref{P:injective}).
      Therefore, by Lemma \ref{L:isom},  $q$ is  an equivalence.
      Pick a (pointed) inverse for $q$ and compose it with $(\Z,z) \to (\Y,y)$
      to obtain the desired lifting.
  \end{proof}


The main theorem of this section is

 \begin{thm}{\label{T:Galois}}
   Let $(\X,x)$ be a pointed connected topological stack. Assume $\X$
      is locally path connected  and semilocally 1-connected
      (Definition \ref{D:lpc}).
      Then, the functor $F_x \: \mathsf{Cov}_{\X} \to (\pX-\mathsf{Set})$ is
      an equivalence of categories.
 \end{thm}

 Before proving Theorem \ref{T:Galois}, we remark that, under the correspondence
stated in the theorem, the connected covering spaces correspond to
transitive $\pX$-sets. Hence, we have the following

 \begin{cor}{\label{C:universalcover}}
      Let $(\X,x)$ be a pointed connected topological stack. Assume $\X$
      is locally path connected  and semilocally 1-connected. Then, there is
      a natural bijection between the subgroups $H \subseteq \pX$
      and the (isomorphism classes of) pointed connected covering maps
      $p \: (\Y,y) \to (\X,x)$. In particular, $\X$ has a
      universal cover.
 \end{cor}

 The group $H$ in the above corollary is nothing but $p_*\big(\pY\big)$
 (which is
 isomorphic to $\pY$ by Proposition \ref{P:injective}).

 \vspace{0.1in}

 Before proving Theorem \ref{T:Galois}, we need some preliminaries.

 \begin{defn}{\label{D:Aequivalence}}
     Let $(\X,x)$ be a topological stack, and let $A$ be a $\pX$-set.
      Let $f \: W \to \X$ be a map from a path connected
     topological spaces $W$  (unpointed) such that the induced map
     on fundamental groups is trivial. By an {\em $A$-equivalence class of
     a path from $x$ to $W$ (via $f$)} we mean a triple $(w,\gamma, a)$,
     where $w$ is a point
     in $W$, $a$ a point in $A$, and $\gamma$ is a path from $x$ to $f(w)$,
     modulo the following equivalence relation:
     \begin{quote}
       $(w,\gamma, a) \sim (w',\gamma', a')$ if there is a path $\al$ in
      $X$ from $w$ to $w'$ such that  the action of
      $\gamma f(\al)\gamma'^{-1} \in \pX$ sends $a$ to $a'$.
      \end{quote}
   \end{defn}
 Note that the choice of $\al$ is immaterial since $f$ induces the trivial
 map on fundamental groups. Also note that if $f,g \: W\to \X$ are
 2-isomorphic maps, then the $A$-equivalence classes of paths from
 $x$ to $W$ via $f$ are in natural bijection with the $A$-equivalence
 classes of paths from $x$ to $W$ via $g$.

\begin{proof}[Proof of Theorem \ref{T:Galois}]
  Fully faithfulness follows from Proposition \ref{P:lifting}.
  We have to prove essential surjectivity.

  Take a  $\pX$-Set $A$. We  construct the covering space $\Y \to \X$
  associated to $A$ as a sheaf (of sets) over the
  comma category $\mathsf{Top}_{\X}$ of topological spaces over $\X$.
  Recall that the comma category $\mathsf{Top}_{\X}$ is defined
  in a similar way to $\mathsf{Cov}_{\X}$ (see page \pageref{Cov}),
  but we take the objects to be {\em all} maps $X \to \X$, where $X$
  is a topological space.

  Define a presheaf $\mathcal{F}$
  on $\mathsf{Top}_{\X}$ as follows. Let $f \: W \to \X$ be an object in
  $\mathsf{Top}_{\X}$. Suppose $W$ is connected and locally path connected, and
  that $f$ induces the trivial map on the fundamental groups.
  Define
      $$\mathcal{F}(f)=\{\text{ $A$-equiv. classes of  paths
                                    from $x$ to $W$.}\}$$
  For an arbitrary object $g \: V \to \X$ in $\mathsf{Top}_{\X}$,
  assume there exists an object $W \to \X$   as above and
  a morphism $\varphi \: g \to f$ in $\mathsf{Top}_{\X}$.
  Then we set $\mathcal{F}(g):=\mathcal{F}(f)$.  Observe that, for
   another choice of $f' \: W' \to \X$ and $\varphi' \: g \to f'$, there
   is a canonical bijection between $\mathcal{F}(f)$ and $\mathcal{F}(f')$
   (exercise).
   So $\mathcal{F}(g)$ is well-defined.
   If such a $W$ does not exists, we set  $\mathcal{F}(V)=\emptyset$.

 Now, define $\Y$ to be the sheafification $\mathcal{F}^a$.
 We claim that $\Y$, viewed as a stack over $\X$, is the desired
 covering stack. To see this, observe that, if $f \: W \to \X$ is an object
 in $\mathsf{Top}_{\X}$ such that $W$ is connected and locally path connected,
 and $f$ induces the trivial map on the fundamental groups, then the
 restriction of $\mathcal{F}^a$ to $\mathsf{Top}_W$ is a constant sheaf.
 In fact, upon fixing a point $w$ in $W$ and a path $\gamma$ in $\X$ from
  $x$ to $f(w)$, this sheaf becomes canonically isomorphic to the
  constant sheaf associated to $A$.

 The above observation translates as saying that, for every such
 $f \: W \to \X$, the fiber product $W\x_{\X}\Y$ is equivalent (as a stack)
 to the topological space $W\x A$.
 Since $\X$ is locally path connected and semilocally
 1-connected, we can find a family of map $W_i \to \X$
 with the above property such that the induced map $\coprod W_i \to \X$
 is an epimorphism. Lemma \ref{L:representable} now  implies that
 $\Y \to \X$ is representable. It is also a covering map by Definition
 \ref{D:representable2} (note that being a covering map is invariant
 under base extension and local on the target).

  It is straightforward to check that the fiber of $\Y \to \X$
  is  canonically isomorphic to $A$ as a $\pX$-set.
 \end{proof}

 \begin{rem}
   The Galois category $\mathsf{Cov}_{\X}$ is in fact defined for every
   connected {\em pretopological} stack $\X$. Also, for a base point $x \in \X$,
    there is a fiber functor $F'_x \: \mathsf{Cov}_{\X} \to \mathsf{Set}$.
    If we define $\pi_1'(\X,x)$
    to be the group of automorphisms of the functor $F'_x$, then it
    can be shown that    $\mathsf{Cov}_{\X}$ is equivalent to
    $\pi_1'(\X,x)-\mathsf{Set}$. When $\X$ is topological, there is a natural
    group homomorphism $\pX \to \pi_1'(\X,x)$. This map is an isomorphism if
    $\X$ is locally path connected and semilocally 1-connected.

  \end{rem}
\subsection{Role of the inertial fundamental groups}{\label{SS:hidden}}

 Having the apparatus of Galois categories at hand, we can
 extend the results of \cite{Noohi1} to the topological setting.
 The role played by inertial fundamental groups is that they control the
 stacky structure of the covering spaces (see Lemma \ref{L:mainlemma}
 and Theorem \ref{T:uniformization}). \footnote{Inertial fundamental groups
  can also be used to compute the fundamental group of the coarse moduli
 space (\oldcite{Noohi2}).}


 \begin{prop}{\label{P:kernelofomega}}
    Let $(\X,x)$ be a pointed topological stack, and let $\tilde{\X}$
    be a universal cover for $\X$. Then, for any point $x_0 \in \tilde{\X}$
    lying above $x$, we have an isomorphism $I_{x_0}\cong \ker\ox$.
\end{prop}

\begin{proof}
  Follows from Lemma \ref{L:mainlemma}.
\end{proof}

 \begin{thm}{\label{T:uniformization}}
   Let $\X$ be a connected locally path connected
   semilocally 1-connected topological stack, and let
   $\tilde\X$ be its universal cover.
   Then the following conditions are equivalent:
    \begin{itemize}

    \item[$\mathbf{i.}$] The maps $\ox \: I_x \to \pX$
       are injective for every point $x$.

    \item[$\mathbf{ii.}$] $\tilde\X$ is a quasitopological space;

    \item[$\mathbf{iii.}$] $\X$ is the quotient stack of the action of a
        discrete group on a quasitopological space.

 \end{itemize}
    \noindent Furthermore, if $\X$ is a  Deligne-Mumford topological stack,
    we can drop `quasi' in
       $(\mathbf{ii})$ and $(\mathbf{iii})$.
 \end{thm}

   \begin{proof}[Proof of  $(\mathbf{i}) \Rightarrow (\mathbf{ii})$]
          This follows from Proposition \ref{P:kernelofomega}.

   \vspace{0.1in}
   \noindent {\em Proof of} $(\mathbf{ii}) \Rightarrow (\mathbf{iii})$.
          Let $R=\tilde{\X}\x_{\X}\tilde{\X}$, and consider
          the groupoid $[R \sst{} \tilde{\X}]$ (in the category
          of sheaves of sets over $\mathsf{Top}$). The source
          and target maps are covering maps, and $\tilde{\X}$ is
          connected and
          simply connected. So $R$ is a disjoint union of copies
          of $\tilde{\X}$ (indexed by $\pi_1\X$), each of which mapping
          homeomorphically to $\tilde{X}$ via the source
          and target maps. Lemma \ref{L:actiongpd} is easily adapted
          to this situation and it implies that
          $[R \sst{} \tilde{\X}]$ is the action groupoid of an action
          of $\pi_1\X$ on $\tilde{\X}$. (Having chosen a base point
          in $\tilde{\X}$ this action of $\pi_1\X$ on $\tilde{\X}$
          could also be recovered from the general Lifting
          Lemma \ref{P:lifting}.)

    \vspace{0.1in}
    \noindent {\em Proof of} $(\mathbf{iii}) \Rightarrow (\mathbf{i})$.
          In this case, $\X$ has a covering stack in which all the
          inertial fundamental groups are zero. The result follows from
          Lemma \ref{L:mainlemma}.

     The final statement of the theorem follows from
     Proposition \ref{P:strongcovering}
     below and the obvious
     fact that
     a quasirepresentable Deligne-Mumford topological stack is
     representable.
 \end{proof}

  \begin{prop}{\label{P:strongcovering}}
        Let $p \: \Y \to \X$ be a covering map of topological
        stacks, and assume that $\X$ is a locally  connected
        Deligne-Mumford topological stack. Then, for every
        point $x \in \X$, there exists an open substack $\U \subseteq \X$
        containing $x$ with the following properties:
           \begin{itemize}

              \item[$\mathbf{i.}$] $\U\cong[U/I_x]$, for some
                topological space $U$ acted on by $I_x$ (mildly at $x$);

              \item[$\mathbf{ii.}$] $p^{-1}(\U)\cong \coprod_{k\in K} [U/G_k]$,
                  where $G_k$, for $k$ ranging in an index set $K$,
                  are subgroups of $I_x$ acting on $U$ via $I_x$.
           \end{itemize}
        In particular,  $\Y$ is a Deligne-Mumford topological stack.
 \end{prop}

  \begin{proof}

    By shrinking $\X$ around $x$ we may assume that $\X=[X/I_x]$, for some
    locally  connected topological space $X$ acted on by $I_x$
    (Lemma \ref{L:lpcDM}).
    Consider the corresponding \'{e}tale chart $X \to \X$, and let
    $Y \to \Y$ be the pull back
    chart for $\Y$. The map $q \: Y \to X$, being a pull
    back of $p$, is again a covering map. Let $x'$ be the (unique)
    point in $X$ lying
    over $x$. Note that $x'$ is a fixed point of the action of
    $I_x$ and $I_x$ acts mildly at $x'$.
    There is an open set $U \subseteq X$ containing $x'$ over which
    $q$ trivializes. After replacing $U$ with a smaller open set
    containing $x$ (say, by the connected component of $x$ in $U$),
    we may assume that $U$ is $I_x$-invariant and connected.
    Set $\U=[U/I_x]$. We claim that $\U$ has the desired property.
    Let $\V=p^{-1}(\U)\subseteq \Y$,
    and  $V=q^{-1}(U)\subseteq Y$. Then,  $V$ is an \'{e}tale chart for $\V$
    and it is of the form
                  $$V=\coprod_{j\in J} U_j, \ U_j=U,$$
    for some
    index set $J$. We can think of $V$ also
    as an \'etale chart for $\U$ (by first
    projecting it down to $U$ via $q$, and then composing with $U \to \U$). The
    corresponding groupoid look as follows:
       $$\big[V\x_{\U}V=\coprod_{J\x J\x I_x} U_{j_1,j_2,g}
             { \xymatrix@=12pt@M=8pt@C=70pt{
                 \ar @<4pt> [r]^{g \: U_{j_1,j_2,g} \to  U_{j_1}}
                 \ar @<-4pt> [r]_{g \: U_{j_1,j_2,g} \to  U_{j_2}} & }}
                                              \coprod_{j\in J} U_j=V\big],$$
    where $U_{j_1,j_2,g}=U$.
    We can think of $V$ also as a chart for $\V$. In this case,
    the corresponding groupoid looks as follows:
         $$R:=[V\x_{\V}V \sst{} V].$$
    Now, the key observation is that, the latter groupoid is an open-closed
    subgroupoid of the former.
    This follows by applying Proposition \ref{P:diagonal}
    to the top left square of
    the following cartesian
    diagram:
        $$\xymatrix@=12pt@M=10pt{ V\x_{\V}V \ar[r] \ar[d] &
                                   V\x_{\U}V \ar[d] \ar[r]& V\x V \ar[d] \\
          \V \ar_(0.4){\De} [r] &  \V\x_{\U}\V \ar[r] \ar[d] & \V\x \V \ar[d] \\
                   & \U \ar[r] & \U\x\U}$$
    Since $U$ is connected, this implies that
             $$V\x_{\V}V=\coprod_A U_{j_1,j_2,g},$$
    where $A$  is a subset of $J\x J \x I_x$. An immediate conclusion is that,
    for each slice $U_{0}$ in $V=\coprod_{j\in J} U_j$,
    the orbit $\OO(U_0)$ is a disjoint
    union of copies of some $U_j$, $j \in J$.
    By considering orbits of different $U_j$, $j\in J$,
    we obtain a partitioning $J=\coprod_{k \in K} J_k$. This gives
    rise to a decomposition $\V = \coprod_{k \in K} \V_k$.
    Let us fix a $k \in K$ and see how each individual $\V_k$ looks like.
    For this, we have to look at
    the restriction $[R|_{\mathbb{U}_k} \sst{} \mathbb{U}_k]$ of the groupoid
    $[R \sst{} V]$ to the invariant
    open
    $\mathbb{U}_k=\coprod_{j\in J_k} U_j$.  Let us fix a
    representative slice $U_{j_k}$  in $\coprod_{j\in J_k} U_j$.
    Since all the slices of $U$ in $\coprod_{j\in J_k} U_j$ are permuted  around
    transitively under $R$,
    we have     $[\mathbb{U}_k/R|_{\mathbb{U}_k}]\cong[U_{j_k}/R|_{U_{j_k}}]$.
    On the other hand,
    $R|_{U_{j_k}}=\coprod_{G_k}U_g$, where $G_k \subseteq I_x$ is a subset
     (we are actually
    thinking of $G_k$ as the subset
    $\{j_k\}\x \{j_k\}\x G_k \subseteq J\x J\x I_x$).
    The fact that $[R|_{U_{j_k}} \sst{} U_{j_k}]$
    is a subgroupoid of $[R \sst{} V]$
    implies that $G_k$ is in fact a subgroup of $I_x$, and that the groupoid
     $[R|_{U_{j_k}} \sst{} U_{j_k}]$ is
     naturally identified with the action groupoid
     $[G_k \x U \sst{} U]$ of $G_k$ acting (via $I_x$) on $U$.
     Therefore, we have $\V_k\cong [U_{j_k}/R|_{U_{j_k}}]\cong[U/G_k]$.
     This completes
     the proof.
  \end{proof}

The example at the end of Section \ref{S:DMFiber} shows that the
above proposition fails if we  do not assume that $\X$ is locally
connected.

\section{Examples}{\label{S:examples}}

In this section we supply some examples of topological stacks.
There are five subsections. The first one collects a few
pathological examples that can be used here and there as
counterexamples. In the next four subsections we consider four
general classes of example, namely: gerbes, orbifolds, weighted
projective lines, and graphs of groups.

\subsection{Some pathological examples}{\label{SS:pathological}}
\begin{ex}{\label{E:twopoint}}
   Let $X=\mathbb{Z}\cup\{\infty\}$ be the one-point compactification
   of $\bbZ$, and
   let $\mathbb{Z}$
   act on $X$ by  translation on $\mathbb{Z} \subset X$ and fixing $\infty$.
   The corresponding quotient stack $[X/\mathbb{Z}]=\X$
   is a connected non locally connected weak Deligne-Mumford topological stack
   that is not Deligne-Mumford. It has a  dense open substack $\U$ that
   is equivalent to a point. The residue gerbe at $\infty$ is the closed
   substack that is complement to $\U$, and it is equivalent
   to $\B\mathbb{Z}$. The coarse moduli space of $\X$ is a 2 element set
   $\{0,\infty\}$ with $\{\infty\}$ its only non trivial open set.
   Note that, although $\X$ is not path connected, it has
   a path connected coarse moduli space.
\end{ex}

 \begin{ex}{\label{E:torus}}
  Let $c$ be a non-zero real number. Consider the action
  of  $\mathbb{Z}$   on $\bbR^2$ given by multiplication by $c$,
  and let $\X$ be the quotient stack. Then $\X$ is a path-connected
  locally path connected weak Deligne-Mumford
  topological  stack that is not Deligne-Mumford. It has an open dense
  substack $\U$ that is equivalent to a (2-dimensional) torus.
  The complement of $\U$ is a closed substack equivalent to $\B\bbZ$
  which can be identified as the residue gerbe at $0$. The coarse moduli
  space of $\X$ is homeomorphic to $\mathbb{T}^2\coprod\{0\}$
  where the topology induced on $\mathbb{T}^2$ is the usual topology
  of a torus,
  and $0$ is in the closure of every point in $\mathbb{T}^2$ (i.e.
  the only open set containing $0$ is the whole space).
  We have $\pi_1\X\cong\bbZ$. All the higher homotopy groups vanish.
  \end{ex}

  \begin{ex}{\label{E:circle}}
 A variation of Example \ref{E:torus} is constructed as follows.
 Consider the vector field $(x,y) \mapsto (x,y)$ on $\bbR^2$
 and look at the
 corresponding flow. This gives an action of $\bbR$ on $\bbR^2$,
 having $0$ as its unique fixed point. The quotient stack $\X$
 of this action is a path connected locally path connected
 topological stack (that is not  weak Deligne-Mumford,
 since the inertia group at $0$ is $\bbR$, which is not discrete).
 It has an open dense
  substack $\U$ that is equivalent to the circle $S^1$.
  The complement of $\U$ is a closed substack equivalent to $\B\bbR$
  which can be identified as the residue gerbe at $0$.
  The coarse moduli
  space of $\X$ is homeomorphic to $S^1\coprod\{0\}$
  where the topology induced on $S^1$ is the usual topology
  and $0$ is in the closure of every point in $\mathbb{T}^2$ (i.e.
  the only open set containing $0$ is the whole space).
  All homotopy groups of $\X$ vanish.
  \end{ex}

\subsection{Topological gerbes}{\label{SS:gerbes}}

Any trivial pretopological gerbe is of the form $B_XG$, where $X$
is a topological space and $G$ is a relative topological group
over $X$. This pretopological  gerbe is topological (resp., weak
Deligne-Mumford, Deligne-Mumford), if the structure map $G \to X$
of the group $G$ is LF (resp., local homeomorphism, resp. covering
space).

Recall that any gerbe $\X$ can be covered by open substacks that
are trivial gerbes, and $\X$ is topological (resp., weak
Deligne-Mumford, Deligne-Mumford), if each of these open substacks
are so.

\begin{ex}{\label{E:trivialgerbe}}
   Let $G$ be a topological group, acting trivially on a topological
   space $X$. Then $[X/G]\cong X\x \B G$ is a trivial gerbe (over $X$).
   This gerbe is topological for any choice of LF in Example
   \ref{E:lf} except for ($\mathbf{6}$).
   If $G$ is a  discrete group, then $X/G$ is Deligne-Mumford.
  \end{ex}

    \begin{ex}{\label{E:pointgerbe}}
      Let $\X$ be a topological stack whose coarse moduli space is a single
      point. Then $\X$ is {\em not} necessarily a   gerbe (Examples
      \ref{E:Q} and \ref{E:Z2}).
      If $\X$ is Deligne-Mumford, however, then $\X$ is a trivial gerbe
      (over $*$).
    \end{ex}

   \begin{ex}{\label{E:nontrivialgerbe1}}
     Let $G$ be a topological group, and   $H \subseteq G$   a closed
     normal subgroup. Let $X$ be a topological space, and  $Y$
     a $G/H$-torsor over $X$. Let $G$ act on $Y$ via $G/H$. Then
     $[Y/G]$ is a topological gerbe (over $X$). This gerbe is trivial if and only
     $Y$ can be extended to a $G$-torsor, that is, if
     there is a $G$-torsor $Z$ over $X$ such that $Y=Z/H$. In particular,
     $[Y/G]$ is trivial when  the extension
     $$1 \to H \to G \to G/H \to 1$$
     is split. (Choose a splitting, and define $Z=Y\x_{G/H}G$.)
   \end{ex}

    \begin{ex}{\label{E:nontrivialgerbe2}}
       In the previous example, take $X=*$ and $Y=G/H$.
       Then $[Y/G]$ is trivial if and only if  the extension
     $$1 \to H \to G \to G/H \to 1$$
     is split.
    \end{ex}

    Let $G$ be an arbitrary topological group. The source
    and target maps of the trivial groupoid $[G\sst{}*]$
    are $\mathbf{LF}$ for, say, $\mathbf{LF}$=locally cartesian maps (Example
    \ref{E:lf}.$\mathbf{3}$). So we can talk about homotopy groups
    of $\B G=[*/G]$.
    Using the torsor description
    of $\B G=[*/G]$, it is easy to show that there are  natural
    isomorphisms:
      $$\pi_{n-1}(G) \risom \pi_n(\B G), \ n=1,2,\cdots.$$
    Here, $\pi_0(G)$ stands for the group of {\em path}
    components of $G$ (we do not need any assumptions  on $G$);
    it is isomorphic to $G/G_0$, where $G_0$ is the path component
    of  identity.
    Take $x \: * \to \B G$ to be the base point of $\B G$. Then, the natural
    map $\ox \: I_x \to \pi_1(\B G,x)$ is simply the quotient
    map $G \to G/G_0$.

    Using this, we obtain another interpretation of the maps
    $\ox \: I_x \to \pX$, for an arbitrary pointed topological stack $(\X,x)$,
    as follows.

    Let $\Gamma_x$ be the residue gerbe at $x$.
    Then, $\Gamma_x\cong \B \mathbb{I}_x$,
    where $\mathbb{I}_x$ is naturally a topological group since $\X$
    is (pre)topological.
      The natural pointed map
    $\Gamma_x \to \X$ induces maps on the homotopy groups
      $$\pi_{n-1}(\mathbb{I}_x)=\pi_n(\Gamma_X) \to \pi_n(\X), \ n=1,2,\cdots.$$
    Setting $n=1$ and pre-composing with
    $\mathbb{I}_x \to \mathbb{I}_x/(\mathbb{I}_x)_0$, we obtain a group
    homomorphism
       $$I_x \to  \pi_0(G) \to \pi_1(\X).\footnote{Recall that
       $I_x$ stands for the underlying discrete group of the topological group
       $\mathbb{I}_x$}$$
     This homomorphism can be canonically identified with $\ox$.

     This shows that each homotopy group $\pi_n(\X)$ comes with a natural
     extra structure, namely, the group homomorphism
     $\pi_{n-1}(\mathbb{I}_x) \to \pi_n(\X)$. For this reason,
     we call $\pi_{n-1}(\mathbb{I}_x)$ the $n^{\text{th}}$ {\em inertial
     homotopy group}. Inertial homotopy groups
     appear in the study of  the loop stack of a topological
     stack.

  Now, let us compute inertia stacks of topological gerbes.
  We begin with the following proposition.

  \begin{prop}{\label{P:doublecoset1}}
  Consider a diagram of discrete groups
   $$\xymatrix@=12pt@M=10pt{
                     &  H\ar[d]^f  \\
            K \ar[r]_g       &  G        &          }$$
  and set $H'=f(H)$ and $K'=g(K)$.
  Let $A$ be the the set of double cosets $H'aK'$ (so we have a partitioning
  $G=\coprod_A H'aK'$). For each double coset   $H'aK'$, choose
  a representative $a$ and let $C_a$
  be the subgroup of $H\x K$ defined by
      $$C_a=\{(h,k)\, | \; f(h)\, a\, g(k)^{-1} =a\}.$$
  Then, we have a natural 2-cartesian diagram
  $$\xymatrix@=12pt@M=10pt{
     \coprod_A \B C_a\ar[r]\ar[d] & \B H\ar[d]    \\
         BK    \ar[r]       &    \B G      }$$
  \end{prop}

  This proposition is a special case of the following.

  \begin{prop}{\label{P:doublecoset2}}
     Consider a diagram of topological groups
       $$\xymatrix@=12pt@M=10pt{
                                  &  H\ar[d]^f  \\
                 K \ar[r]_g       &  G        &          }$$
     Consider the action of $H\x K$ on $G$ in which
     the effect of $(h,k) \in H\x K$  sends
     $a \in G$ to $f(h)\, a\, g(k)^{-1}$. Then we have
     a natural 2-cartesian diagram
         $$\xymatrix@=12pt@M=10pt{
           [G/(H\x K)]  \ar[r]\ar[d] & \B H\ar[d]    \\
              BK    \ar[r]       &    \B G      }$$
  \end{prop}

  \begin{proof}
    This follows from the general construction
    of Section \ref{S:fiber}. Details left to the reader.
  \end{proof}

  \begin{cor}{\label{C:inertiaBG}}
    Let $G$ be a topological group. Then, the inertia stack of $\B G$
    is equivalent to the quotient stack of the conjugation action of
    $G$ on itself. In particular, when $G$ is discrete,
    we have
      $$\mathcal{I}\B G\cong\coprod_A \B C_a,$$
    where $A$ is the set of
    conjugacy classes in $G$, and $C_a$ is the stabilizer group of
    (a representative $a$ of) a conjugacy class.
  \end{cor}

 \begin{proof}
     The inertia stack $\mathcal{I}\B G$ of $\B G$ is equivalent
     $\B G\x_{B(G\x G)} \B G$,
     where the maps $\B G \to B(G\x G)$ are induced from the diagonal map
     $\De \: G \to G\x G$.
     So we have
     an action of $G\x G$ on $X:=G\x G$, where an element $(g_1,g_2)$ acts on
      $(x,y)$ as follows:
       $$(x,y) \mapsto (g_1\, x\, g_2^{-1}, g_1\, y\, g_2^{-1}).$$
     (We use the notation $X$ for $G\x G$ when it is viewed as a space
     acted on by $G\x G$.) Hence, $\mathcal{I}\B G$  is equivalent
     to the quotient stack $[X/(G\x G)]$. Consider the subspace
     $Y=G\x \{1\}$ of $X$, and let $[R \sst{} Y]$ be the restriction
     of the action groupoid $[(G\x G)\x X \sst{} X]$. The quotient
     stack $[Y/R]$ is a substack of $\mathcal{I}\B G$. On the other hand,
     from the fact that $s^{-1}(Y) \to X$ is an epimorphism (in
     fact, it has a section), it follows that
     $Y \to \mathcal{I}\B G$ is an epimorphism.
     This implies that $[Y/R]\cong\mathcal{I}\B G$
     (Proposition \ref{P:quotient1}). We leave it to the reader
     to verify that $[R \sst{} Y]$ is indeed th action groupoid
     of the conjugation action of $G$ on $Y$ (where $Y$ is identified
     with $G$ is the obvious way).
 \end{proof}

\subsection{Orbifolds}{\label{SS:orbifolds}}

 An {\em orbifold} in the sense of Thurston \cite{Thurston} is
 a Deligne-Mumford topological stack whose isotropy groups are finite
 and  admits an \'{e}tale
 chart $p \: X \to \X$ with $X$ a manifold. It is also assumed that
 there is an open dense substack of $\X$ that is a
 manifold.\footnote{Thurston also assumes that
 $\Xm$ is Hausdorff. I am not sure why we need this.} The 2-category
 $\mathbf{Orb}$ of orbifolds is a full sub 2-category of the 2-category
 of Deligne-Mumford topological stacks. It is customary in the literature to
 identify 2-isomorphic morphisms in $\mathbf{Orb}$ and
 work with the resulting 1-category.

 The covering theory developed by
 Thurston in {\em loc. cit.} is a special case of
 the covering theory of topological stacks   developed in
 Section \ref{S:covering}. This is because of
 Proposition \ref{P:strongcovering}. In particular, the existence of a
 universal cover for an orbifold (\cite{Thurston}, Proposition 13.2.4)
 follows from Corollary \ref{C:universalcover}.
 Our definition of the fundamental group of a topological
 stack generalizes Thurston's definition.

 What Thurston calls a {\em good} orbifold is what we have called
 a {\em unifromizable} orbifold. We have the following.

  \begin{thm}{\label{T:orbifoldunif}}
    Let $\X$ be a connected orbifold, and let $I_x$ be the inertia group
    attached
    to the point $x \in \X$. Then, we have natural  group
    homomorphisms $I_x \to \pX$.
    Furthermore,
    $\X$ is a
    good orbifold if and only if for every $x \in \X$
     the map $I_x \to \pX$ is injective.
  \end{thm}

  \begin{proof}
    This follows from Theorem \ref{T:uniformization}.
  \end{proof}

  In   \cite{Noohi2} we give a formula for computing the
  fundamental group of the underlying space (read, coarse moduli space)
  of a connected orbifold $\X$.\footnote{In fact, the formula is valid
  for a quite general class of topological stacks, including Deligne-Mumford
  topological stacks. Under some mild conditions, from this we can derive
  a formula for the fundamental group of the (naive) quotient of a topological
  group acting on topological space (possibly with fixed points).}
  Roughly speaking, the formula can be
  stated as follows. Consider the normal
  subgroup $N \subseteq \pX$
  generated by images in $\pX$
  of  $I_{x'}$, for all  $x' \in \X$ (for
  this we have to identify $\pi_1(\X,x')$ with $\pX$ by choosing a path from
  $x'$ to $x$). Then, $\pX/N$ is naturally isomorphic to the fundamental group
  of the underlying space of $\X$.

\subsection{Weighted projective lines}{\label{SS:weighted}}

  \begin{ex}{\label{E:weightedPL}}
    Let
    $\mathbb{C}^{\times}$ act on $\mathbb{C}^2\backslash \{0\}$   by
    $t\cdot(x,y)=(t^mx,t^ny)$, where $m$ and $n$ are fixed
    positive integers. The corresponding quotient stack is
    called the {\em weighted projective line of weight} $(m,n)$, and
    is denoted by $\PP(m,n)$. It is a Deligne-Mumford topological stack
    (indeed, it is the topological stack associated to
    a Deligne-Mumford   algebraic stack with the same name).
    Using the theory of fibrations of topological stacks (which will
    appear in a sequel to this paper), we can compute homotopy
    groups of $\PP(m,n)$. We have  a  $\mathbb{C}^{\times}$-fibration
    over $\PP(m,n)$ whose total space is $\mathbb{C}^2\backslash\{0\}$.
    A fiber homotopy exact sequence argument shows that
    $\pi_k\PP(m,n)$ is isomorphic to $\pi_kS^2$, for $k \geq 1$.
    In particular, $\PP(m,n)$ is simply connected, for all $m,n \geq 1$.
    The coarse moduli space of $\PP(m,n)$ is homeomorphic to $S^2$.
    It has two distinguished points, corresponding to  $[1:0]$ and $[0:1]$
    in $\mathbb{C}^2$. The residue gerbe at these points
    are isomorphic to $\B\mathbb{Z}_m$ and $\B\mathbb{Z}_n$, respectively.
    At every other point the residue gerbe is isomorphic to $\B\mathbb{Z}_d$,
    where $d=\gcd(m,n)$.

    \vspace{0.1in}

    \noindent{\bf Exercise}. Compute the effect of the moduli map
    $\PP(m,n) \to \PP(m,n)_{mod}\cong S^2$ on the homotopy groups.

    \vspace{0.1in}

    Let us now remove a point (other than the two special points) from
    $\PP(m,n)$. Denote the resulting stack by $\U$.
    An easy van Kampen argument shows that $\pi_1\U\cong\bbZ_m*_{\bbZ_d}\bbZ_n$.
    All other homotopy groups of $\U$ vanish. This follows from the fact
    that, using Theorem \ref{T:uniformization}, $\U$ has a universal
    cover that is a 2-dimensional simply connected manifold.  It
    can be seen  that the only possibility is that the universal cover
    be homeomorphic to $\mathbb{R}^2$.

  \end{ex}

   More details and more examples of this type can be found in \cite{BN}.

\subsection{Graphs of groups}{\label{SS:graphs}}

 In this subsection, we briefly indicate how the theory of
 {\em graphs of groups} (for example see \cite{Serre}, \cite{Bass})
 can be embedded in the   theory of Deligne-Mumford topological stacks.
 Many of the basic results of the theory of graphs of groups are more or less
 immediate consequences of the existence of a homotopy theory for
 Deligne-Mumford topological stacks.\footnote{
    A covering theory for graphs of groups has been developed by Bass
    in \oldcite{Bass}
    in which he points out:
    {\em The definition, and verifications of its properties,
    is regrettably more technical than one might anticipate.}
    Our theory of covering spaces of topological stacks generalizes
    Bass's theory. }

 \vspace{0.1in}
 We start with the simplest case, namely, the graph of groups that looks
 as follows:

 $$\xymatrix@M=0pt{G_0 \, \bullet \ar@{-}[r]^A & \bullet \, G_1}$$

 Here, $G_0$, $G_1$ and $A$ are discrete groups, and the incidence of the
 edge labeled $A$ with the vertex labeled $G_i$ indicates an
 injective group homomorphism $A \hra G_i$.

 We would like to think of this
 as a Deligne-Mumford topological stack $\X$ whose  coarse moduli
 space is the interval $I=[0,1]$. The inertia groups at the end points
   of this interval are $G_0$ and $G_1$, respectively, and the inertia group
 of any other point in the open interval $(0,1)$ is $A$.
 Here is how to make this
 precise.

    Consider the star shaped topological space
    \vspace{-0.1mm}

             $$Y:= (G_0/A) \x I\Big/(G_0/A)\x\{0\}.$$

    \noindent  This is a single point, with
    a bunch of rays coming out of it, each labeled by  a coset of $A$ in $G_0$.
    There is a natural action of $G_0$ on this space, fixing the
    vertex $0$ and permuting the rays emanating from it. This action
    is mild at $0$ (Definition \ref{D:mild}). Set $\Y=[Y/G_0]$; $\Y$
    is a Deligne-Mumford topological stack. Do the similar
    thing with $G_1$, namely, take
     \vspace{-0.1mm}

      $$Z:=(G_1/A) \x I\Big/(G_1/A)\x\{1\}$$

    \noindent and set $\Z=[Z/G_1]$. The complement of the point $0$ in $\Y$
    is an open substack of $\Y$ that is
    equivalent to $B_{(0,1)}A\cong BA \x (0,1)$.
    Similarly, the complement of the point $1$ in $\Z$
    is an open substack of $\Z$ that is equivalent to $BA \x (0,1)$.
    We can glue $\Y$ and $\Z$ along these open substacks
    (Corollary \ref{C:gluing})
    to obtain a Deligne-Mumford topological  stack $\X$
    (Theorem \ref{T:gluingDM}). This Deligne-Mumford  topological
    stack is what we want
    to think of as the stacky incarnation of our graphs of groups.

    Let us play around a  bit with this stack. Using van Kampen
    theorem (whose proof will appear elsewhere), the fundamental group
    of $\X$ is easily computed to be isomorphic to $G=G_1*_AG_2$.
    Note that,  the possible inertial fundamental groups
    of $\X$ at various points are one of $A$, $G_1$ and $G_2$, and these all
    inject into $G$. So by  Theorem \ref{T:uniformization}, the
    universal cover  of $\X$ is an honest topological space.
    That is, there is a topological space $X=\tilde{\X}$ with a properly
    discontinuous action (Definition \ref{D:pd})
    of $G$ such that $[X/G]\cong\X$. By restricting $X$ over
    the open substacks $\Y$ and $\Z$, we see that $X$ locally looks
    like either $Y$ or $Z$; thus, $X$ is a graph. In fact, $X$
    is a tree, because it is simply connected. The upshot is
    that, there is a tree $X$ such that our graph of groups
    is isomorphic to $[X/G]$, where $G=G_1*_AG_2$.

    This tree is easy to construct explicitly (see \cite{Serre}).
      Take the set of vertices of $X$ to be $V=G/G_1\coprod G/G_2$, and
     the set of edges to be  $E=G/A$. The end points of edges are determined by
     maps
     $G/A \to G/G_0$ and $G/A \to G/G_1$. There is a natural action
     of $G$ on this graph that has the following properties:
      \begin{itemize}
         \item All the edges are permuted transitively. The vertices are
            partitioned  into two orbits $P_i=G/G_i$, $i=0,1$.
         \item The inertia group of an interior point of an edge is $A$.
            The inertia
            group of a point in $P_i$ is $G_i$, $i=1,2$.
     \end{itemize}

    \noindent Therefore, the  quotient stack $[X/G]$ is equivalent to $\X$.

      \vspace{0.1in}

      Another basic graph of groups looks like this:

    
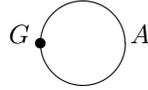
\begin{center}
      \begin{picture}(0,60)
        \put(20,30){\circle{30}}
        \put(4,30){\circle*{4}}
        \put(-8,30){$G$}
        \put(38,30){$A$}
      \end{picture}
    \end{center}

    \noindent Here, $G$ and $A$ are discrete groups. The incidence of the two
    ends of the loop with the vertex correspond to injective group
    homomorphisms $i,j\: A \hra G$. We can identify $A$ with a subgroup
    of $G$ via $i$, and think of $j$ as an injective group homomorphisms
    $\theta \: A \to G$.

    We will explain shortly how to turn this graph of groups
    into a Deligne-Mumford topological stack. van Kampen theorem implies
    that the fundamental group of this graph of groups is isomorphic
    to the HNN-extension  associated to the data $(G,A,\theta)$.

    In general, a graph of groups is  defined to be a graph $\mathcal{G}=(E,V)$,
    (possibly with multiple edges and loops) equipped with a family
    $G_v$ of groups, one for each vertex $v \in V$, and a family of groups
    $G_e$, one for each edge $e \in E$.
    We are also given an inclusion $G_v \hra G_e$
    for an incidence of a vertex $v$ with an edge $e$.

     Using  Theorem \ref{T:gluingDM},
     we can turn a graph of groups into a Deligne-Mumford topological stacks.
     Here is one way of doing it. First, a star-shaped
     graph of groups (i.e. one
     having a single vertex in the center with rays coming
     out of it) can be turned into a Deligne-Mumford topological stack in a way
     analogous to the very first example of this subsection. For a general
      graph of groups we do the following.
     For every edge $e$, pick two distinct points
     $e_{\frac{1}{3}}$ and  $e_{\frac{2}{3}}$, different from the end points.
     Let
         $$\LL=\coprod_{e \in E} \B G_e\x \{\frac{1}{3},\frac{2}{3}\},$$
         $$\Y=\coprod_{e\in E} \B G_e\x [\frac{1}{3},\frac{2}{3}],$$
         $$\Z=\coprod_{v \in V} \mathcal{G}_p,$$
     where $\mathcal{G}_p$ is the star-shaped graphs of groups whose center
     is the vertex $p$
     (think of $\Z$ as the graph of groups that remains when we remove the
     middle third of every edge).
     There are natural embeddings $\LL \hra \Y$ and
     $\LL \hra \Z$ that can be used to glue $\Y$ and $\Z$ along $\LL$
     (Theorem \ref{T:gluingDM}).
     The glued Deligne-Mumford topological stack is a stacky model for
     our graph of groups $\mathcal{G}$.

     The structure theory for graphs of groups developed by Serre
      \cite{Serre} is best interpreted
     from the stacky point of view.

     The fundamental group
     of a graph of groups $\mathcal{G}$ (\cite{Serre}, Section 5.1) is simply
     the fundamental
     group of the topological stack associated to $\mathcal{G}$.
     The existence of the universal cover (\cite{Serre}, Section 5.3)
     follows from Corollary \ref{C:universalcover} as follows.

     Let $\tilde{\mathcal{G}}$ be the universal cover of  $\mathcal{G}$.
     By Proposition \ref{P:strongcovering}.$\mathbf{ii}$, $\tilde{\mathcal{G}}$
     is again a graph of groups. It is easy to check that, for any
     Deligne-Mumford topological stack $\X$, the induced map $\pi_1(\X) \to
     \pi_1(\Xm)$ is surjective (in fact, $\pi_1(\Xm)$ is explicitly computed
     in \cite{Noohi2}). This implies that the coarse moduli space of
     $\tilde{\mathcal{G}}$ is a tree. Now, an easy van Kampen argument shows
     that, for $\tilde{\mathcal{G}}$ to have trivial fundamental
     group it is necessary and sufficient that all the groups $G_v$ and $G_e$
     be trivial; that is, $\tilde{\mathcal{G}}\cong X$ for some {\em tree} $X$.
     So $\mathcal{G}$ is equivalent to $[X/\pi_1(\mathcal{G})]$, where $X$
     is a tree.

     The Structure Theorem (\cite{Serre}, Section 5.4, Theorem 13)
     is completely obvious from the stacky point of view.

     The injectivity of $G_v \to \pi_1(\mathcal{G})$ (\cite{Serre},
     Section 5.2, Corollary 1, whose proof therein is quite tedious and long),
     is an immediate
     consequence of Theorem \ref{T:uniformization}, because we just showed
     that the universal cover of $\mathcal{G}$ is an honest topological
     space (namely, a tree).

\vspace{0.1mm}

    The residue gerbes of a graph of groups are easy to figure out:
    if $x$ is a vertex  of $\mathcal{G}$, then the residue gerbe
    at $x$ is $\B G_x$; if $x$ lies on the interior of an edge $e$,
    then the residue gerbe at $x$ is $\B G_e$.

    More interesting is the inertia stack $\mathcal{I}\mathcal{G}$
    of a graph of groups
    $\mathcal{G}$. This is again a graph of group whose underlying graph
    has its set of vertices
       $$\{v_{[g]} \,|\, v \in V, \ [g]\ \text{a conj. class in} \ G_v\}.$$
    The set of edges is given by
       $$\{e_{[g]} \,|\, e \in E, \ [g]\ \text{a conj. class in} \ G_e\}.$$

    A vertex $v_{[g]}$ is incident with an edge $e_{[h]}$,
    if $v$ is incident with $e$, and, under the inclusion
    $G_v \hra G_e$, the conjugacy class $[g]$ maps to the conjugacy
    class $[h]$.

   The group associated to the vertex $v_{[g]}$ is the centralizer
   $C_g \subseteq G_v$ (we have to make a choice of a representative
   $g \in [g]$). Similarly, the group associated to the edge $e_{[h]}$ is
   the centralizer
   $C_h \subseteq G_e$ (again, we have to make a choice of a representative
   $h \in [h]$).

   For an incident pair $v_{[g]}$ and $e_{[h]}$, the inclusion map
   $C_h \hra C_g$ is defined by the composition
        $$\xymatrix@=12pt@M=10pt{
          C_h \subseteq G_e \ar@{^{(}->} [r] &
                 G_v \ar[rr]^{\text{conj. by $x$}}   &&  G_v     }$$
   where $x \in G_v$ is an element such that conjugation by $x$
   sends (the image in $G_v$ of) $h$ to $g$; again, we have to make
   a choice for $x$. It is easy to see that under this composition $C_h$ lands
   inside $C_g$.

    Up to equivalence, the resulting stack will be independent of
    all the choices made.

\section{The topological stack associated to an algebraic stack}
 {\label{S:algebraic}}

Take the class $\mathbf{LF}$ to be any of the Example \ref{E:lf}
except for ($\mathbf{6}$).

  All  algebraic stacks and schemes considered in this section are assumed
  to be locally of finite type over $\mathbb{C}$.

  In this section we show how to associate a topological stack to
  an algebraic stack (locally of finite type)  over $\mathbb{C}$.
  That is, we describe how to construct  a functor of 2-categories
   $$-^{top} \: \mathsf{AlgSt}_{\mathbb{C}} \to \mathsf{TopSt}.$$
  \noindent (Of course, this functor factors through the 2-category
  of analytic stacks, but we will not discuss this here.) We will
  then prove the stacky version of the {\em Riemann Existence Theorem}.

   \begin{thm}[\oldcite{Milne},
  Chapter III, Lemma 3.14]{\label{T:Riemann}}
     Let $X$ be an algebraic space that is
      locally of finite type over $\mathbb{C}$, and let
     $X^{top}$ be the associated topological space. The functor
     $Y \mapsto Y^{top}$ defines an equivalence between the category of
     finite \'{e}tale maps  $Y \to X$ and the category of finite covering
     spaces of $X^{top}$.
   \end{thm}

  \noindent We remark that the above theorem is stated in \cite{Milne}   for
    passage from schemes to  {\em analytic} spaces. Passage from analytic
    spaces to topological spaces is straightforward. Also, extending the result
    from {\em schemes} to algebraic spaces is easy (and is implicit in what
    follows).

  \vspace{0.1in}
  \noindent{\em Sketch of the construction of $\X^{top}$.} Assume $\X$ is an
    algebraic stack that is
    locally of finite type over $\mathbb{C}$. Let $\X\cong[X/R]$ be a
    presentation of $\X$ as the quotient of a smooth groupoid $[R \sst{} X]$
    (so, $X$ and $R$ are  locally of finite type over $\mathbb{C}$).
    We can now consider the topological groupoid $[R^{top} \sst{} X^{top}]$.
    Every smooth map of  $f  \: X \to Y$ schemes (or algebraic spaces)
    looks locally (in the \'{e}tale topology) like a Euclidean projection
    $\mathbb{A}^m \to \mathbb{A}^n$. That is, for any $x \in X$, after replacing
    $X$ and $Y$ by suitable Zariski open sets, the map $f$ fits in a commutative
    diagram
           $$\xymatrix@=14pt@M=10pt{ X \ar[r]^{\text{\'{e}tale}} \ar[d]_f &
                                             \mathbb{A}^m \ar[d]^{pr}  \\
                        Y \ar[r]_{\text{\'{e}tale}}        &    \mathbb{A}^n  }$$
   \noindent Here $m \geq n$ and $\mathbb{A}^n$ is identified as the
   first $n$ coordinates of $\mathbb{A}^m$.

     Using the fact that
      \'{e}tale maps between schemes induces local homeomorphisms on the
     associated topological spaces, the above diagram implies that
      smooth maps  induce locally cartesian maps with Euclidean fibers
      (see Example \ref{E:lf}.$\mathbf{3'}$)
     on the associated topological
     spaces. Therefore, the groupoid $[R^{top} \sst{} X^{top}]$ has
     $\mathbf{LF}$ source and target maps, where $\mathbf{LF}$ is
     any of classes of maps considered in Example \ref{E:lf},
     except for ($\mathbf{6}$). We denote the quotient (topological) stack
     $[X^{top}/R^{top}]$ by $\X^{top}$.

     If $X_1 \to \X$ and $X_2 \to \X$ are two different
     smooth charts for $\X$, with $[R_1 \sst{} X_1]$ and $[R_2 \sst{} X_2]$
     the corresponding groupoids, there is always a third
     one $X_3 \to \X$ dominating
     both (say $X_3=X_1\x_{\X}X_2$, also see
     Section \ref{S:Morita}), which gives rise to canonical
     equivalences $[X_3^{top}/R_3^{top}] \to [X_1^{top}/R_1^{top}]$
     and $[X_3^{top}/R_3^{top}] \to [X_2^{top}/R_2^{top}]$.  This way we can
     construct an equivalence $[X_2^{top}/R_2^{top}] \to [X_1^{top}/R_1^{top}]$
     which, using a similar argument, can be shown to be independent of the
     choice of $X_3$, up to a 2-isomorphism. In other words,
     $\X^{top}$ is well-defined
     up to an isomorphism that is unique up to 2-isomorphism.

      More generally,
     given a morphism $f \: \Y \to \X$ of algebraic stacks, one can realize it
     as a morphism of groupoids from  $[T \sst{} Y]$ to $[R \sst{} X]$.
     This gives a map
     $f^{top} \: \Y^{top} \to \X^{top}$ of
     topological stacks and again it is easy to show that this map is
     well-defined up to a unique 2-isomorphism.  This, modulo a lot of
     choices that need to be made
     (and hence lots of set theoretical problems that are swept under
     the carpet), will give us the effect of the functor
     $-^{top} \: \mathsf{AlgSt}_{\mathbb{C}} \to \mathsf{TopSt}$ on objects and
     morphisms. Once the choices are made, and the effect of $-^{top}$
     on objects and 1-morphisms is determined,
     the effect of this functor on 2-morphisms can be traced
     through the construction and be easily seen to be well-defined
     and uniquely determined. Therefore, we get a functor of 2-categories.

    %
    %
    %

      \begin{prop}{\label{P:fiberproduct}}
        The functor $-^{top} \: \mathsf{AlgSt}_{\mathbb{C}} \to \mathsf{TopSt}$
        commutes with fiber products (more generally, with all finite
        limits).
      \end{prop}

      \begin{proof}
        When all the stacks involved are affine schemes the result is easy
        to verify because an affine scheme of finite type over $\mathbb{C}$
        can be realized as the zero set of a finitely many polynomial in some
        affine space $\mathbb{A}^n$. The case of fiber products of general
        schemes can be reduced to the affine case, since both on the algebraic
        side and the topological side the   fiber products
        can be constructed locally (that is, by constructing fiber products of
        affine patches and then gluing them together).

        The general case of stacks now follows from the construction of fiber
        products of stacks using their groupoid presentations (
        Section \ref{S:fiber}).
      \end{proof}

     The functor $-^{top}$ has other nice properties: it sends representable
     morphisms to
     representable morphisms (use Lemma \ref{L:representable}),
      smooth morphisms to local fibrations,
     \'{e}tale morphisms to local homeomorphisms, and finite \'{e}tale morphisms
     to covering spaces. It also sends Deligne-Mumford stacks
     (with finite stabilizer) to
     Deligne-Mumford topological
      stacks. The latter is due to the fact that a locally
     Noetherian  Deligne-Mumford stack, with finite stabilizer, has a
     coarse moduli space and is locally, in the \'{e}tale topology
     of its coarse moduli space, a quotient stack of a finite group
     action.

\vspace{0.1in}

     The following result will be used in the proof of Riemann existence theorem
     for stacks.

      \begin{lem}{\label{L:finiteetale}}
        Let $[R \sst{} Y]$ and $[T \sst{} Y]$ be smooth groupoids (in the
        category of schemes), and let
        $f \: T \to R$  be a map of groupoids (the map
        $Y \to Y$ being identity)
        such that the induced map
        $[Y/T] \to [Y/R]$ on the quotients is (representable) finite \'{e}tale.
        Then
        $f$ is an open-closed embedding. The topological version of the
        statement is also true.
      \end{lem}

      \begin{proof}
        Let $\Y=[Y/T]$ and $\X=[Y/R]$. Then $T=Y\x_{\Y}Y$ and $R=Y\x_{\X}Y$.
        So we have a 2-cartesian diagram
           $$\xymatrix@=12pt@M=10pt{  T \ar[r]^f \ar[d]     & R \ar[d]  \\
                        \Y \ar[r]_(0.35){\De}      &  \Y\x_{\X}\Y         }$$
       \noindent The lower map is an open-closed embedding since $\Y \to \X$
       is finite \'{e}tale, therefore so is $f$. The topological version
       is proved similarly (also see Proposition \ref{P:diagonal}).
      \end{proof}

     Now, we are ready to prove the Riemann Existence Theorem for stacks.

    \begin{thm}[Riemann Existence Theorem for Stacks]{\label{T:StackyRiemann}}
     Let $\X$ be an algebraic stack that is
     locally of finite type over $\mathbb{C}$, and let
     $\X^{top}$ be the associated topological stack. The functor
     $\Y \mapsto \Y^{top}$ defines an equivalence between the category of
     (representable)
     finite \'{e}tale maps  $\Y \to \X$ and the category of finite covering
     stacks of $\X^{top}$. (Note that these are honest 1-categories, because we
     are only considering {\em representable} maps.)\footnote{
     More precisely, 2-categories in which there is at most one
     2-isomorphism between any two 1-morphisms.}
    \end{thm}

    \begin{proof}
      First we prove essential surjectivity. Let $\Y' \to \X^{top}$
      a covering map (we use the notation $\Y'$, instead of $\Y$,
      to emphasis that $\Y'$ is a {\em topological} stack -- in fact,
      as a matter of notational convention, in this proof
      every thing that has a superscript is topological, otherwise algebraic).
      Let $X \to \X$ be a smooth chart for $\X$, with $X$ a scheme.
      So $X^{top} \to \X^{top}$ will be an LF chart for $\X^{top}$.
      Let $Y'=\Y'\x_{\X^{top}}X^{top}$ be the pull-back chart for $\Y'$.
      We can think of the composition $Y' \to X^{top} \to \X^{top}$
      as an LF chart for $\X^{top}$.
      Set $T'=Y'\x_{\Y'}Y'$ and  $R'=Y'\x_{\X^{top}}Y'$. By (the topological
      version of) Lemma \ref{L:finiteetale}, $T'$ is an open-closed subspace of
      $R'$. Since $Y' \to X^{top}$ is a covering space,
      by the classical version of the Riemann Existence Theorem,
      there exists a scheme $Y$, locally of finite type over $\mathbb{C}$,
      such that $Y^{top}=Y'$. Let $R=Y\x_{\X}Y$.
      By Proposition \ref{P:fiberproduct}, we have $R^{top}\cong R'$.
      It is obvious that the decomposition of $R$ into its connected components
      is preserved under $-^{top}$. Therefore, $T'=T^{top}$ for some open-closed
      subspace of $R$. It is easy to see that $[T \sst{} Y]$ is a subgroupoid
      of $[R \sst{} Y]$. If we let $\Y=[Y/T]$, then $\Y \to \X$ is finite
      \'{e}tale,  and $\Y^{top}\to \X^{top}$ is isomorphic,
      as a covering space of
      $\X^{top}$, to $\Y' \to \X^{top}$. This proves
      the essential surjectivity.

      Next we prove fully-faithfulness. Let $\Y$ and $\Z$ be finite \'{e}tale
      covers of $\X$. We may assume $\Y$ is connected. To give a map
      from $\Y$ to $\Z$ relative to $\X$ is the same
      is to specify a connected component of $\Y\x_{\X}\Z$ that maps
      isomorphically  to $\Y$.  Using Proposition \ref{P:fiberproduct},
      this is the same as to specify a connected component of
      $\Y^{top}\x_{\X^{top}}\Z^{top}$ that maps isomorphically  to $\Y^{top}$,
      and these are of course in bijection with maps
      from $\Y^{top}$ to $\Z^{top}$ relative to $\X^{top}$. This proves
      the fully-faithfulness.
    \end{proof}

     \begin{cor}{\label{C:profinitecompletion}}
       Let $\X$ be a connected algebraic stack that is
       locally of finite type over $\mathbb{C}$. Assume $\X^{top}$ is
       locally path connected and semilocally 1-connected. Then
       we have an isomorphism
         $$\pi_1^{alg}(\X) \to \widehat{\pi_1(\X^{top})}.$$
     \end{cor}


\bibliographystyle{amsplain}
\bibliography{homotopy}
\end{document}